\magnification=1200
\font\Large=cmbx10 scaled \magstep2
\def\res{\mathop{\rm res}}
\def\deg{\mathop{\rm deg}}
\hfill{Preprint SISSA 24/98/FM}
\vskip 3cm
\centerline{\bf\Large Painlev\'e transcendents}
\smallskip
\centerline{\bf\Large in two-dimensional topological field theory}
\bigskip
\centerline{\bf Boris DUBROVIN}
\medskip
\centerline{SISSA, Trieste}
\bigskip
\centerline{\bf Contents}
\medskip
Introduction
\smallskip
Lecture 1. Algebraic properties of correlators in 2D topological field
theory. Moduli of a 2D TFT and WDVV equations of associativity.
\smallskip
Lecture 2. Equations of associativity and Frobenius manifolds. Deformed
flat connection and its monodromy at the origin.
\smallskip
Lecture 3. Semisimplicity and canonical coordinates.
\smallskip
Lecture 4. Classification of semisimple Frobenius manifolds.
\smallskip
Lecture 5. Monodromy group and mirror construction for semisimple
Frobenius manifolds.
\smallskip
References.
\vfill\eject

\centerline{\bf Introduction}
\medskip
This paper is devoted to the theory of WDVV equations of associativity.
This remarkable system of nonlinear differential equations was discovered
by E.Witten [Wi1] and R.Dijkgraaf, E.Verlinde and H.Verlinde [DVV] in the
beginning of '90s. It was first derived as equations for the so-called
primary free energy of a family of two-dimensional topological field
theories.  Later it proved to be an efficient tool in solution of problems
of the theory of Gromov - Witten invariants, reflection groups and
singularities, integrable hierarchies.

Here we mainly consider the relationships of WDVV to the theory of
Painlev\'e equations. This is a two-way connection. First, any solution to
WDVV satisfying certain semisimplicity conditions, can be expressed via
Painlev\'e-type transcendents. Conversely, theory of WDVV works as a
source of remarkble particular solutions of the Painlev\'e equations.

The paper is an extended version of the lecture notes of a course given
at 1996 Carg\`ese summer school ``The Painlev\'e property: one century
later". It is organized as follows.

In Lecture 1 we give a sketch of the ideas of two-dimesnional topological
field theory, we formulate WDVV and give main examples of solutions coming
from quantum cohomology and from singularity theory. In Lecture 2 we give
a coordinate-free reformulation of WDVV introducing the notion of
Frobenius manifold. We also construct the first main geometrical object,
namely, the deformed affine connection on a Frobenius manifold. The
monodromy at the origin of the deformed connection gives us the first
set of important invariants of Frobenius manifolds. In Lecture 3 we define
the class of semisimple Frobenius manifolds. In physics they correspond to
two-dimensional topological field theories with all relevant
perturbations. We construct the so-called canonical coordinates on such
manifolds. In Lecture 4 we complete the classification of semisimple
Frobenius manifolds in terms of monodromy data of certain universal
linear differential operator with rational coefficients. We give a
nontrivial example of computation of the monodromy data in quantum
cohomology. In the last Lecture we develop a ``mirror construction''
representing the principal geometrical objects on a semisimple Frobenius
manifold by residues and oscillatory integrals of a family of analytic
functions on Riemann surfaces.
\smallskip
{\bf Acknowledgment.} I would like to thank the organizers of the
Carg\`ese summer school for the invitation and generous support.
I thank A.B.Givental for fruitful discussion of Theorem 3.2.

\vfill\eject
\centerline{Lecture 1.}
\medskip
{\bf Algebraic properties of correlators in 2D topological field
theories.}

{\bf Moduli of a 2D TFT and WDVV equations of associativity}
\medskip
By definition, a quantum field theory (QFT) on a D-dimensional oriented
manifold $\Sigma$ (in our case $D=2$) consists of:

1) Local fields $\phi_\alpha (x)$, $x\in \Sigma$. The metric $g_{ij}(x)$ 
on
$\Sigma$ could be one of the fields. It is called gravity.

2) Lagrangian
$$L=L(\phi, \partial_x \phi, \dots).
$$
The equations of motion of the classical field theory have the form
$${\delta S\over \delta \phi_\alpha (x)} =0
$$
where
$$S[\phi] = \int_\Sigma L(\phi, \partial_x \phi, \dots)
$$
is the classical action.

3) In the path-integral quantization we are interested in the partition
function
$$Z_\Sigma = \int [d\phi ] e^{-S[\phi]}
$$
and, more generally, in the (non-normalized) correlation functions
$$<\phi_\alpha (x)\phi_\beta (y) \dots >_\Sigma 
= \int [d\phi] \phi_\alpha (x) \phi_\beta (y) \dots  e^{-S[\phi]}.
$$
The integration in both cases is over the space of local fields $\phi$ on
$\Sigma$ with an appropriate measure $[d\phi]$. In the full theory we are
also to take an integration over the space of manifolds $\Sigma$.

4) The theory admits topological invariance if an arbitrary change of the
metric on $\Sigma$ preserves the action
$${\delta S\over \delta g_{ij} (x)} \equiv 0.
$$
In $D=2$ case such a theory will be called 2D topological field theory
(TFT). For example, in the 2D case the total curvature functional
$$S={1\over 2\pi} \int_\Sigma R \sqrt{g} d^2 x
$$
is topologically invariant. Indeed, due to Gauss - Bonnet theorem it is
equal to the Euler character of the surface $\Sigma$.

For a topological field theory the partition function gives a topological
invariant of $\Sigma$. The correlation functions depend only on the
topology of $\Sigma$ and on the fields (but not on their positions).
Particularly, in the 2D case we have
$$<\phi_\alpha (x) \phi_\beta (y)\dots >_\Sigma  
\equiv <\phi_\alpha  \phi_\beta \dots >_g.
$$
In the r.h.s. there are just numbers depending on the genus $g$ of the
surface $\Sigma$ and on the labels $\alpha, ~ \beta, \dots$ of the
fields.

5) In the {\it matter sector} of the QFT we integrate over the space
of all fields but the metric $\left( g_{ij}(x)\right)$. For a TFT the
correlators of the matter sector have a nice algebraic description to
be presented in a moment. To describe coupling of the QFT to gravity one
is
to integrate over the space of metrics. In TFTs coupling to gravity can be
reduced to integration over the space of conformal classes of the metrics
on $\Sigma$, i.e., over the moduli space of Riemann surfaces of the genus
$g=g\left(\Sigma \right)$. This is a much more complicated procedure
by now fixed only for the genera $g=0, ~1$.
\medskip
We describe now the algebraic properties of the matter sector correators
in a 2D TFT. We will consider simple theories having a finite number of
observables in the matter sector
$$\phi_1, \dots, \phi_n
$$
(the so-called primary chiral fields). One can easily derive all algebraic
properties of the correlators using the general Atiyah axioms of a
topological field theory. We present here only the summary of the
properties.
\smallskip
{\bf Definition 1.1.} A {\it Frobenius algebra} is a pair $\left( A,
~<~,~>\right)$ where $A$ is a commutative associative algebra (over {\bf
C}) with a unity and $<~,~>$ stands for a symmetric
non-degenerate {\it invariant} bilinear form on $A$. The invariance
means validity of the following identity
$$<a\, b, c> = < a, b\, c>
\eqno(1.1)
$$
for arbitrary 3 vectors $a, ~b, ~c \in A$.
\smallskip
{\bf Theorem 1.1} (see [Dij1, Dij1, Du7]). {\it The matter sector
correlators of any
2D
TFT with $n$ observables can be encoded by a Frobenius algebra $\left( A,
<~,~>\right)$ with a marked basis $e_1, \dots, e_n$. The genus $g$
correlators of the observables have the form
$$<\phi_{\alpha_1} \phi_{\alpha_2}\dots \phi_{\alpha_k}>_g = <e_{\alpha_1}
\cdot e_{\alpha_2} \cdot 
\dots
\cdot \phi_{\alpha_k}, H^g>
$$
where
$$H=\eta^{\alpha\beta} e_\alpha \cdot e_\beta \in A
$$
$$\left( \eta^{\alpha\beta}\right) = \left(
\eta_{\alpha\beta}\right)^{-1},
~~~\eta_{\alpha\beta} := <e_\alpha, e_\beta>.
$$
}
Physicists call   $\left( A,
<~,~>\right)$ the primary chiral algebra of the TFT. Observe that the
structure of the Frobenius algebra is uniquely determined by the genus
zero two- and three-point correlators
$$<e_\alpha, e_\beta> = <\phi_\alpha \phi_\beta>_0, ~~~ <e_\alpha \cdot
e_\beta, e_\gamma> = <\phi_\alpha \phi_\beta \phi_\gamma>_0.
$$
Usually the observables are chosen in such a way that the vector $e_1$
coincides with the unity of the algebra $A$. Then
$$<e_\alpha, e_\beta> = <\phi_1 \phi_\alpha \phi_\beta>_0.
$$
\medskip
We give now the two main ``physical'' examples of Frobenius algebras.

{\bf Example 1.1.} Let $X$ be a 2d-dimensional closed oriented manifold
without odd-dimensional cohomologies. Take the full cohomology algebra
$$A= H^*(X)
$$ 
with the bilinear form
$$<\omega_1, \omega_2> = \int_X \omega_1 \wedge \omega_2, ~~~ \omega_1,
~\omega_2 \in H^*(X)
\eqno(1.2)
$$
(we realize cohomologies by classes of closed differential forms).
Symmetry and invariance of this bilinear form are obvious. Nondegeneracy
follows from the Poincar\'e duality theorem. This Frobenius algebra
describes the matter sector of the topological sigma model ($X$ is the
target space).

Actually, in this example we have a certain graded structure on  $\left(
A,
<~,~>\right)$. Generalizing, we give
\smallskip

{\bf Definition 1.2.} The Frobenius algebra is called {\it graded} if a
linear operator $Q: A\to A$ and a number $d$ are defined such that
$$
Q(a \, b) = Q(a) b + a Q(b),
\eqno(1.3a)
$$
$$
<Q(a),b> + <a, Q(b)> = d <a,b>
\eqno(1.3b)
$$ 
for any $a, ~b \in A$. The operator $Q$ is called {\it grading operator}
and the number $d$ is called {\it charge} of the Frobenius algebra. We
will consider only the case of diagonalizable grading operators. Then we
may
assign degrees to the eigenvectors $e_\alpha$ of $Q$
$$
\deg (e_\alpha) =q_\alpha ~~{\rm if} ~~ Q(e_\alpha) =q_\alpha e_\alpha.
\eqno(1.4a)
$$
For the topological example the vectors of a homogeneous basis are chosen
in such a way that
$$
e_\alpha \in H^{2 q_\alpha} \left(X\right), ~~\deg (e_\alpha) =
q_\alpha.
\eqno(1.4b)
$$
The charge $d$ is equal to the half of the dimension of $X$.

Particular example: $X={\bf CP}^d$. The full cohomology space has
dimension
$n=d+1$. The natural basis in $A=H^*\left( {\bf CP}^d\right)$ is
$$1, ~ \omega, ~ \omega^2, \dots, \omega^d
$$
where $\omega$ is the standard K\"ahler form on the projective space. We
normalize it by the condition
$$\int_{{\bf CP}^d}\omega^d =1.
$$
Then $\left( A, <~,~>\right)$ is isomorphic to the quotient of the
polynomial algebra
$$A= {\bf C}[\omega]/\left(\omega^{d+1}\right)
$$
with the bilinear form
$$<\omega^k, \omega^l> = \delta_{k+l, d}.
$$
\smallskip
{\bf Remark 1.1.} We will consider below also graded Frobenius algebras
$\left( A, <~,~>\right)$ over
graded commutative associative rings $R$. In this case we have two grading
operators $Q_R: R\to R$ and $Q_A:A\to A$
satisfying the properties
$$
\eqalignno{Q_R(\alpha\beta) & = Q_R (\alpha)\beta + \alpha Q_R(\beta),~
\alpha, ~\beta \in R & (1.5a)
\cr
Q_A(ab) & = Q_A(a)b+aQ_A(b), ~a, ~b \in A & (1.5b)
\cr
Q_A(\alpha a) &= Q_R(\alpha)a + \alpha Q_A(a), ~\alpha\in R, ~a\in A
& (1.5c)
\cr
Q_R<a,b> +d\, <a,b> &=<Q_A(a),b> + <a, Q_A(b)>, ~a, ~b\in A.
& (1.5d)
\cr}
$$
The number $d$ is called {\it the charge} of the graded Frobenius algebra.
\medskip
{\bf Example 1.2.} of Frobenius algebra. Let $f(x)$ be a polynomial of
$x\in
{\bf C}^N$ with an isolated singularity at $x=0$. This means that 
$$
df(x)|_{x=0} =0
$$
(we may also assume that $f(0)=0$),
$$
df(x)|_{x\neq 0} \neq 0
$$
for $x$ sufficiently close to the origin. Take the quotient of the
polynomial algebra
$$
A={\bf C}[x]\big/ _{\left({\partial f\over \partial x_1}, \dots,
{\partial f\over \partial x_N}\right)}.
\eqno(1.6)
$$
This is called the Jacobi ring, or the local algebra of the singularity.
This is a finite-dimensional algebra if the singularity has finite
multiplicity $n$. (The number $n={\rm dim}\, A$ is also called Milnor
number of the singularity.) We define bilinear form on $A$ taking the
residue
$$<p, q> ={1\over \left(2 \pi i\right)^N} \int_{\cap_i |{\partial f\over
\partial x_i}|=\epsilon}
{p(x) q(x) d^N x\over {\partial f\over
\partial x_1}\dots {\partial f\over
\partial x_N}}.
\eqno(1.7)
$$
Here $\epsilon$ is sufficiently small positive number. Again, symmetry
and invariance of the bilinear form are trivial. Nodegeneracy is less
trivial; see the proof in [AGV], Volume 1, Section 5. To obtain a graded
Frobenius algebra one is to take a quasihomogeneous polynomial $f(x)$.
This Frobenius algebra describes the matter sector of a topological Landau
- Ginsburg model. The function $f(x)$ is called superpotential of the
theory.

Particular example: the simple singularity of the $A_n$ type. Here $N=1$,
$f(x)=x^{N+1}$. The local algebra
$$A={\bf C}[x]\big/ \left( x^{N+1}\right) =
{\rm span}\, \left(1, x, x^2, \dots, x^{n-1}\right)
$$
$$x^k\cdot x^l =\cases{ x^{k+l}, ~k+l<n\cr
0, ~k+l\geq n}$$
$$<x^k, x^l> =\res {x^{k+l}\over (n+1) x^n} =\cases{0, ~k+l\neq n-1\cr
{1\over n+1}, ~k+l=n-1.\cr}
$$
The grading operator is determined by 
$$Q(x) = {1\over n+1}x,
$$
the charge is 
$$d={n-1\over n+1}.
$$
\medskip
We have already said that the procedure of coupling to gravity of a 2D TFT
is more complicated (not settled in full generality). For the genus zero
case it can be done still in an axiomatic way. It turns out that the
axioms of coupling to gravity can be reduced to WDVV equations of
associativity. Here WDVV stands for Witten - Dijkgraaf - E.Verlinde -
H.Verlinde. In the paper [Wi1] the equations of associativity were derived
in the setting of topological sigma models. In [DVV] they were derived in
a more general class of TFTs obtained by the so-called twisting from $N=2$
supersymmetric QFTs. Basically, the idea was to consider correlators of
a particular $n$-dimensional family of TFTs
$$S\mapsto S-\sum_{\alpha=1}^n \int_\Sigma \phi_\alpha^{(2)}
$$
as functions of the coupling constants $t=\left( t^1, \dots, t^n\right)$.
Here $\phi_1^{(2)}, \dots, \phi_n^{(2)}$ are certain two-forms on $\Sigma$
being in one-to-one correspondence with the observables $\phi_1, \dots,
\phi_n$. The deformation preserves the topological invariance (not the
grading!). So one obtains a $n$-dimensional deformation  $\left( A_t,
<~,~>_t\right)$
of $n$-dimensional
Frobenius algebra $\left( A, <~,~>\right) =  \left( A_0, <~,~>_0\right)$.
A basis $e_1=1, e_2, \dots, e_n$ corresponding to the chosen system of
observables $\phi_1,  \dots, \phi_n$ is marked in all of the algebras
$A_t$. The following properties of the family of Frobenius algebras 
 $\left( A_t, <~,~>_t\right)$ were proved by WDVV:
$$<e_\alpha, e_\beta>_t\equiv <e_\alpha,e_\beta>
\eqno(WDVV1)
$$
$$c_{\alpha\beta\gamma}(t):=<e_\alpha\cdot e_\beta, e_\gamma>_t
={\partial^3 F(t)\over \partial t^\alpha \partial t^\beta \partial
t^\gamma}.
\eqno(WDVV2)
$$
Here $F(t) =\log Z_0(t)$ is the genus zero free energy of the family of
TFTs (the so-called primary free energy).

The last is the quasihomogeneity condition: the structure constants
$c_{\alpha\beta\gamma}(t)$ of the algebras $A_t$ are weighted homogeneous
functions of the degree $q_\alpha+q_\beta+q_\gamma-d$ where we assign
degree $1-q_\alpha$ to the variable $t^\alpha$ for each $\alpha=1, \dots,
n$:
$$c_{\alpha\beta\gamma}\left(\lambda^{1-q_1} t^1, \dots, \lambda^{1-q_n}
t^n\right) =
\lambda^{q_\alpha+q_\beta+q_\gamma-d}c_{\alpha\beta\gamma}\left(t^1,
\dots, t^n\right)
$$
for an arbitrary $\lambda\neq 0$. Observe that $q_1=0$ if $e_1=1$. All the
quasihomogeneity equations can be written as a one for the primary free
energy
$$F\left(\lambda^{1-q_1} t^1, \dots, \lambda^{1-q_n}
t^n\right) 
=\lambda^{3-d} F\left( t^1, \dots, t^n\right) + {\rm quadratic}
\eqno(WDVV3)
$$
where ``quadratic'' stands for at most degree two polynomial in $t^1,
\dots, t^n$. (Later we will slightly modify the quasihomogeneity
requirement for those $t^\alpha$ where $q_\alpha =1$ - see the beginning
of Lecture 2.)

The WDVV equations of associativity is the problem of classification of
$n$-dimensional
families of $n$-dimensional Frobenius algebras satisfying the above
properties WDVV1 - WDVV3. One can consider this problem as the first
approximation to the problem of classification of 2D TFTs, at least
of those obtained by twisting from N=2 supersymmetric theories. We do not
present here other stores of the whole building of a 2D TFT (coupling to
gravity [Wi2, Du3, Du7], Zamolodchikov-type Hermitean metric on the space
of
parameters $t$ [CV1, Du4]). Probably, the upper stores can be put not on
an
arbitrary solution of WDVV as on the basement. However, before proceeding
to the upper stores we will study the structure of the eventual basement,
i.e., of a solution of WDVV. These my lectures are devoted just to this
problem of classification of solutions of WDVV equations of associativity.

We finish this section with a sketch of construction of the deformed
2D TFTs for the two above examples. Observe first that for a graded
Frobenius algebra $(A_0, <~,~>, Q,d)$ one can construct a trivial
cubic solution of WDVV
$$
F_0 = {1\over 6} <1, (t)^3>, ~~t=t^\alpha e_\alpha \in A_0.
\eqno(1.8)
$$
In all the physical examples the free energy $F(t)$ is constructed as an
analytic perturbation of a cubic $F_0$.
\medskip
{\bf Example 1.3.} We will additionally assume the $2d$-dimensional target
space $X$ to be K\"ahler. The deformation of $F_0$ is defined as the
generating function of Gromov - Witten invariants. Let us consider the
moduli space of instantons
$$X_{[\beta ], l}:= \left\{
{\rm holomorphic}~ \beta:
\left( S^2, p_1,\dots, p_l\right) 
\to X, ~
{\rm given ~homotopy~class}~
[\beta ]\in H_2 \left( X;{\bf Z}\right)\right\}.
\eqno(1.9)
$$
The holomorphic maps $\beta$ of the Riemann sphere $S^2$ with marked
points $p_1, \dots, p_l$ are considered up to holomorphic change of
parameter. Under certain assumptions about the manifold $X$
(see [KM, RT]) it can be shown that $X_{[\beta],l}$ can be compactified
to produce an orbifold of the complex dimension
$${\rm dim}_{\bf C} X_{[\beta ],l} =d + \int_{S^2} \beta^{*} \left(
 c_1(X)\right) + l-3.
$$
Here $c_1(X)\in H^2(X)$ is the first Chern class of $X$.

Observe that any of the marked points $p_i$ defines the evaluation map
that we denote by the same symbol
$$
p_i: ~X_{[\beta ],l} \to X, ~~\left(\beta, p_1, \dots, p_l\right)
\mapsto
\beta(p_i).
\eqno(1.10)
$$
For an element
$$a_1 \otimes a_2 \otimes \dots \otimes a_k \in \left( H^{*}
 (X)\right)^{\otimes k}
$$
define the number
$$
<a_1\otimes \dots \otimes a_k >_{[\beta ],l} = \cases{ 0, ~k\neq l\cr
\int_{X_{[\beta ],l}} p_1^{*} (a_1) \wedge \dots \wedge p_l^{*} (a_l),
~k=l\cr}
\eqno(1.11)
$$
We extend this symbol linearly onto the infinite direct sum
$${\bf C} \oplus H^{*} \oplus H^{*} \otimes H^{*} \oplus \left(
H^{*} \right)^{\otimes 3} \oplus \dots
$$
with $H^{*} := H^{*} (X)$.

Define now the function $F(t)$,
$$
t=\left( t', t''\right) \in H^{*} (X)
\eqno(1.12a)
$$
$$
t'\in H^2(X)/2 \pi i H^2 (X,{\bf Z}), ~~t''\in H^{{*} \neq 2}(X),
\eqno(1.12b)
$$
$$
F(t) =F_0(t) +\sum_{[\beta ]\neq 0, l} \left< e^{t''}\right>_{[\beta ],
l} e^{\int_{S^2} \beta^{*} (t')}.
\eqno(1.13)
$$
Here $F_0(t)$ is the cubic (1.8) for the Frobenius algebra $A_0
=H^{*}(X)$.
The exponential
$$e^t := 1 +{t\over 1!} + {1\over 2!} t\otimes t + \dots
$$
is considered as an element of the infinite direct sum.

The numbers $<a_1 \otimes \dots \otimes a_k>_{[\beta ], l}$ can be nonzero
only if the following dimension condition holds true
$${\rm deg} a_1 + \dots {\rm deg} a_l = {\rm dim} \, X_{[\beta ],l} 
= d+\int_{S^2} \beta^{*} \left( c_1(X)\right) + l-3.
$$
This can be written in the form
$$
\sum_{i=1}^l \left( 1 - {\rm deg}\, a_i\right) = 3 -d -
\int_{S^2}\beta^{*} \left( c_1(X)\right).
\eqno(1.14)
$$
We see from this dimension condition that for any $[\beta ], l$ the
coefficient 
$$\left< e^{t''}\right>_{[\beta ], l}
$$
is a polynomial in $t''\in H^{{*}\neq 2}(X)$. The coefficients of these
polynomials proved to be independent on the complex structure on $X$
[Gr, MS, RT] but only on the homotopy class of the symplectic structure
on $X$
given by the imaginary part $\Omega$ of the K\"ahler metric. They are
called
{\it Gromov - Witten invariants} of $\left(X, \Omega\right)$. (Actually,
one can start with more general situation to define GW invariants of a
compact symplectic manifold $\left( X, \Omega\right)$. To this end one is
to consider pseudoholomorphic maps $\beta : S^2 \to X$ w.r.t. an
appropriate almost complex structure on $X$. See details in [Gr, MS, RT].)

The family of algebras $A_t$ with the parameter 
$$t\in H^{*}(X)/2 \pi i
H^2(X,{\bf Z})
$$
is called {\it quantum cohomology} of $X$.
Sometimes they considered quantum cohomology in the restricted sense where
the parameter $t=t'$ of the deformation belongs to 
$$
t\in H^{1,1}(X)/2 \pi i
H^2(X,{\bf Z}).
$$ 
This restricted quantum cohomology is closely related
to Floer symplectic cohomology of $\left(X, \Omega\right)$ (see [Sad,
Pi, MS]). In
the {\it point of classical limit} $t'\to -\infty$ (i.e.,
$\int_{S^2}\beta^{*} (t') \to -\infty$ for any $[\beta ]\neq 0$)
$F(t)\to f_0(t)$, so the quantum cohomology goes to the classical ones.

Particular example. Quantum cohomology of the projective plane
${\bf CP^2}$. For $t=t^1 + t^2 \omega + t^3 \omega^2 \in H^{*} ({\bf
CP^2})$
the cubic function $F_0(t)$ is
$$F_0(t) = {1\over 2} \left( t^1\right)^2 t^3 + {1\over 2} t^1 \left(
t^2\right)^2.
$$
Here $t' = t^2 \omega$, $t''=t^1 + t^3 \omega^2$. The series $F(t)$
has the form [KM1]
$$
F(t)=F_0(t) + \sum_{k=1}^\infty {N_k\over (3k-1)!}
\left( t^3\right)^{3k-1} e^{k t^2}.
\eqno(1.15)
$$
Here
$$N_k = \# \left\{ {\rm rational~curves~of~degree}~k~{\rm on}~{\bf CP^2}
\right.
$$
$$\left. {\rm passing~through}~3k-1~{\rm generic~points.}
\right\}
$$
E.g., $N_1=1$ (one line through 2 points), $N_2=1$ (one conic through 5
points). One can see that the quasihomogeneity condition WDVV3 must be
modified: the function $F(t)$ has degree $1=3-2$ (up to quadratic terms)
if $t^1$ has degree 1, $t^3$ has degree -1, $t^2$ has degree 0 but
$\exp t^2$ has degree 3. This quasihomogeneity anomaly comes from the term
$$\int_{S^2} \beta^{*} \left( c_1\left( {\bf CP^2}\right)\right)
$$
in the dimension condition (1.14).

The series $F(t)$ has nonempty domain of convergence
$$
{\rm Re}\, \left( t^2 + 3 \log t^3 \right) < R
\eqno(1.16)
$$
for some positive $R$. Numerical estimation for $R$ was obtained by [DI]
$$R\simeq 1.981.
$$
Actually, the following asymptotic ansatz was proposed in 
[DI]
$${N_k\over (3k-1)!} \simeq a^k b k^{-{7\over2}}, ~~k->\infty
$$
with $a\simeq 0.138$, $b\simeq 6.1$. The exact values of the constants
$a$, $b$ are not known.

The structure constants of the restricted quantum cohomology ring are
obtained by triple differentiation of $F(t)$ and setting $t^1=t^3=0$.
The resulting ring has very simple structure: this is the quotient
of the polynomial ring
$$QH^{*} \left({\bf CP^2}\right) ={\bf C}[e_2]/\left(e_2^3=q\right)
$$
with
$$q=e^{t^2}.
$$
Clearly, at the point of classical limit $q\to 0$ one obtains the
classical cohomology ring of the projective plane.

The function $F(t)$ proves to solve the WDVV equatons of associativity
[KM1].
It was observed by Kontsevich that, plugging the ansatz (1.15) with $N_1=1$
into the equations of associativity one can compute recursively all the
coefficients $N_k$. We leave as an exercise for the reader to derive these
recursion relations for $N_k$.
\smallskip
{\bf Remark 1.2.} Denote 
$$
\phi(x)=\sum_{k=1}^\infty {N_k\over (3k-1)!} e^{kx}
$$
and
$$
\psi(x) ={\phi'''-27\over 8(27+2 \phi'-3\phi'')}
$$
(the prime stands for the $x$-derivative). Then the coefficients
$N_K^{(1)}$ of the expansion
$$
\psi(x)=-{1\over 8} + \sum_{k=1}^\infty {k\,N_k^{(1)}\over (3 k)!} e^{kx}
$$
are the elliptic Gromov - Witten invariants of ${\bf CP}^2$, i.e., they
are the numbers of elliptic curves of the degree $k$ passing through $3k$
generic points on ${\bf CP}^2$. This was proved in [DZ2].
\medskip

Also in the general situation of quantum cohomology of a manifold $X$ one
can prove validity of WDVV for a vast class of manifolds $X$ [KM, MS, RT].
The quasihomogeneity conditions have the form WDVV3 for the dependence of
$F(t)$ on the coordinates of the component $t''\in H^{{*}\neq
2}\left(X\right)$. For the other component $t'=\sum {t'}^\alpha
e'_\alpha\in H^2\left( X\right)$
of $t=\left( t', t''\right)$ the coordinates ${t'}^\alpha$ are
dimensionless. We assign then the degrees to the exponentials
$$
\deg e^{{t'}^\alpha} = r_\alpha
\eqno(1.17)
$$
if
$$
c_1\left( X\right) = \sum r_\alpha e'_\alpha.
\eqno(1.18)
$$
Clearly, for $X={\bf CP^2}$ we obtain the above condition $\deg \exp t^2
=3$. For Calabi - Yau (CY) varieties $X$ also the exponentials $\exp
{t'}^\alpha$ are dimensionless since $c_1 \left( X\right) =0$.
Particularly, for CY 3-folds all the GW polynomials
$$\left < e^{t''}\right >_{[\beta ], l}
$$
are just numbers, as it follows from the dimension condition (1.14). That
means that, essentially, the full quantum cohomology of a CY 3-fold is
reduced to the restricted one (we do not consider here the contributions
from the odd-dimensional classes of th CY). According to mirror
conjecture [COGP], the free energy of a CY 3-fold $X$ can be expressed via
certain generalized hypergeometric functions. These hypergeometric
functions are periods of the holomorphic three-form on the so-called dual
CY 3-fold $X^{*}$. The mirror conjecture has been proved in [Gi2 - Gi4]
for CY
complete intersections in projective spaces. A general geometrical setting
justifying mirror conjecture was proposed in [Wi4].

In the opposite case of Fano varieties, where $c_1\left( X\right) >0$,
nothing is known about the analytic structure of the free energy
(besides the trivial example of the projective line where the full quantum
cohomology is reduced to the restricted one). The restricted quantum
cohomology can be often computed (actually, they are computed for all Fano
complete intersections in [Beau]). For many examples of Fano varieties 
it was shown that, like in the above example of ${\bf CP^2}$, one can
reconstruct all the GW invariants from the restricted quantum cohomology
just solving recursively the WDVV equations of associativity. The
restricted quantum cohomology serve as the initial data to specify
uniquely the solution of WDVV.

We suggest that the success of this reconstruction of GW invariants of
Fano varieties, unlikely CY varieties, where WDVV gives essentially no
information about the GW invariants, is based on the following conjectural
property [TX] of quantum cohomology of Fano varieties: the deformed
Frobenius algebra $A_t$ is semisimple for generic value of the parameter
$t$. In these lectures we describe the general solution of WDVV satisfying
the semisimplicity condition. We will show that they can be expressed via
certain Painlev\'e-type transcendents. We will also discuss the problem of
selection of the particular solutions of WDVV corresponding to free
energies of physically motivated models of 2D TFT.
\medskip

{\bf Example 1.4.} In the topological Landau - Ginsburg models with the
superpotential $f(x)$ the deformed Frobenius algebra is given by the
formulae similar to (1.6), (1.7) where one is to use the versal deformation 
[AGV, Ar2]
$$
f_s(x) = f(x) + \sum_{i=1}^n s^i p_i(x)
\eqno(1.19)
$$
of the singularity. Here $p_1(x)=1$, $p_2(x)$, \dots, $p_n(x)$ is a basis
of the local algebra of the singularity. (Actually, one is to choose
properly the volume form $d^Nx$ in (1.7). The construction of the needed
volume form is given in [Sai2].) The metric
$$\sum \eta_{ij}(s)\, ds^i ds^j
$$
on the space of the parameters $s=\left( s^1, \dots, s^n\right)$ has the
form
$$
\eta_{ij}(s) ={1\over (2\pi i )^N} \int_{\cap_j |{\partial
f_s(x)\over \partial x_j}|=\epsilon} {p_i(x) p_j(x) d^N(x)\over
{\partial f_s(x)\over \partial x_1} \dots {\partial f_s(x)\over
\partial x_n}}.
\eqno(1.20)
$$
Under certain assumptions [Sai2] one can prove that this metric has zero
curvature. Thus one
can introduce new coordinates $\left( t^1, \dots, t^n\right)$ on the space
of parameters such that
$$\eta_{ij} ds^i ds^j =\eta_{\alpha\beta} dt^\alpha dt^\beta
$$
with a constant matrix $\eta_{\alpha\beta}$. In these coordinates
$$
c_{\alpha\beta\gamma} (t) ={1\over (2\pi i)^N}
 \int_{\cap_j |{\partial
f_s(x)\over \partial x_j}|=\epsilon} {{\partial f_s(x)\over \partial
t^\alpha}{\partial f_s(x)\over \partial
t^\beta}{\partial f_s(x)\over \partial
t^\gamma} d^N(x)\over
{\partial f_s(x)\over \partial x_1} \dots {\partial f_s(x)\over
\partial x_n}}.
\eqno(1.21)
$$
The explicit formulae for the A-D-E simple singularities see in [BV].

Particular case. Simple singularity of $A_3$ type. Here $p_1=1$, $p_2=x$,
$p_3=x^2$ is a basis in the local algebra. So
$$f_s = x^4 + s_1 + s_2 x + s_3 x^2
$$
(I use all lower indices in concrete examples). The metric (1.7) in the
coordinates $s_1$, $s_2$, $s_3$ has the matrix depending on $s$
$$<p_i, p_j>_s = - 4 \res_{x=\infty} {p_i(x) p_j(x)\over 4 x^3 + 2 s_3 x +
s_2}.
$$
We obtain the following matrix of the metric
$$\eta_{ij}(s) = \left( \matrix{0 & 0 & 1\cr
0 & 1 & 0\cr
1 & 0 & -{1\over 2} s_3\cr}\right).
$$
Introducing the new coordinates 
$$\eqalign{s_1 &=t_1 +{1\over 8} t_3^2\cr
s_2 &= t_2\cr
s_3 &=t_3\cr}
$$
we obtain the constant matrix
$$\eta_{\alpha\beta} = \left( \matrix{0 & 0 & 1\cr 0 & 1 & 0\cr 1 & 0 &
0\cr} \right).
$$
The new parametrization of the versal deformation has the form
$$P_t(x) \equiv f_s(x) = x^4 + t_1 + {1\over 8} t_3^2 + t_2 x + t_3
x^2.
$$
The only nontrivial ``three-point functions''
$$c_{\alpha\beta\gamma} =- \res_\infty {\partial_\alpha P_t \partial_\beta
P_t \partial_\gamma P_t\over \partial_x P_t} dx
$$
are
$$c_{113}=c_{122}=1, ~~c_{223}=-{1\over 4} t_3, ~ c_{233} =-{1\over 4}
t_2, ~ c_{333} = {1\over 16} t_3^2.
$$
This gives a polynomial solution of WDVV
$$
F(t_1, t_2, t_3) = {1\over 2} t_1^2 t_3 + {1\over 2} t_1 t_2^2 -{1\over
16} t_2^2 t_3^2 + {1\over 960} t_3^5.
\eqno(1.22)
$$
We can continue our experiments with WDVV and try to find {\it all}
polynomial solutions $F(t_1, t_2, t_3)$. This simple exercise gives only
4 polynomial solutions [Du6, Du7]! Besides (1.22) they are
$$
\eqalignno{F & = 
{1\over 2} t_1^2 t_3 + {1\over 2} t_1 t_2^2 +{1\over
6} t_2^3 t_3 + {1\over 6} t_2 ^2 t_3^3 + {1\over 210} t_3^7
& (1.23) \cr
F & = 
{1\over 2} t_1^2 t_3 + {1\over 2} t_1 t_2^2 +{1\over 6} t_2^3 t_3^2
+{1\over 20} t_2^2 t_3^5 +{1\over 3960} t_3^{11}
& (1.24) \cr
F & = 
{1\over 2} t_1^2 t_3 + {1\over 2} t_1 t_2^2 +t_2^4. \cr
}
$$
The last polynomial does not satisfy the semisimplicity condition. It
turns out that the first two can be described in terms of singularities 
of the type $B_3$ and $H_3$ respectively. In Lecture 5 I will explain the
construction of the polynomials (1.22) - (1.24) in terms of invariants of 
the Coxeter
groups of the type $A_3$, $B_3$, $H_3$ resp. and the generalization of
this construction to higher dimensions. Observe that the Coxeter groups
$A_3$, $B_3$, $H_3$ are just all the groups of symmetries  of Platonic
solids (of tetrahedron, octahedron, and icosahedron resp.). So, WDVV
equations of associativity ``know'' not only enumeration of rational plane
curve, but they also ``know'' the list of Platonic solids! See also
Conjecture of Lecture 5 below regarding polynomial solutions of WDVV.

\vfill\eject
\def\diag{{\rm diag}\,}
\def\res{\mathop{\rm res}}
\def\deg{\mathop{\rm deg}}
\def\wdvv{WDVV equations of asociativity}
\def\V{{\cal V}}
\def\C{{\bf C}}
\def\U{{\cal U}}
\def\L{{{\cal L}_E}}
\centerline{Lecture 2}
\medskip
\centerline{\bf Equations of associativity and Frobenius manifolds}
\smallskip
\centerline{\bf Deformed flat connection and its monodromy at the origin}
\bigskip
We give first the precise formulation of \wdvv . Next, we will
reformulate them in a coordinate-free form.

We look for a function
$F\left( t^1, \dots, t^n\right) \equiv F(t)$, a constant symmetric
nondegenerate matrix $\left( \eta^{\alpha\beta}\right)$, numbers $q_1,
\dots, q_n$, $r_1, \dots, r_n$, $d$ such that
$$\partial_\alpha\partial_\beta\partial_\lambda F(t) \eta^{\lambda\mu}
\partial_\mu\partial_\gamma\partial_\delta F(t) = 
\partial_\delta\partial_\beta\partial_\lambda F(t) \eta^{\lambda\mu}
\partial_\mu\partial_\gamma\partial_\alpha F(t)
\eqno(WDVV1)
$$
for any $\alpha, ~\beta, ~\gamma, ~\delta \, = \, 1, \dots, n.$ (We denote
$$\partial_\alpha := {\partial\over\partial t^\alpha} 
$$
etc., summation over repeated indices is assumed.) Equivalently, the
algebra
$$A_t = {\rm span}\, \left( e_1, \dots, e_n\right)
$$
with the multiplication law
$$
\eqalign{e_\alpha\cdot e_\beta &=  c_{\alpha\beta}^\gamma (t)
e_\gamma\cr
c_{\alpha\beta}^\gamma(t) :&= \eta^{\gamma\epsilon}
\partial_\epsilon\partial_\alpha\partial_\beta F(t)\cr}
\eqno(2.1)
$$
is to be associative for any $t$. The algebra will automatically be
commutative. 

The symmetric nondegenerate bilinear form $<~,~>$ on $A_t$ defined by
$$
<e_\alpha,e_\beta > := \eta_{\alpha\beta}
\eqno(2.2)
$$
where the matrix $\left( \eta_{\alpha\beta}\right)$ is the inverse one to 
 $\left( \eta^{\alpha\beta}\right)$, is invariant (in the sense of (1.1))
since the expression
$$
<e_\alpha\cdot e_\beta, e_\gamma > =
\partial_\alpha\partial_\beta\partial_\gamma F(t)
\eqno(2.3)
$$
is symmetric w.r.t. any permutation of $\alpha$, $\beta$, $\gamma$.

The variable $t^1$ will be marked and we require that
$$\partial_\alpha\partial_\beta\partial_1 F(t)\equiv \eta_{\alpha\beta}.
\eqno(WDVV2)
$$
This means that the first basic vector $e_1$ will be the unity of all the
algebras $A_t$. From (WDVV1,2) we conclude that $\left( A_t, <~,~>\right)$
is a Frobenius algebra for any $t$.

The last one is the quasihomogeneity condition that we write down in the
infinitesimal form using the Euler identity for the quasihomogeneous
functions. Introducing the Euler vector field
$$
E =\sum_{\alpha=1}^n \left[ (1-q_\alpha ) t^\alpha + r_\alpha\right]
\partial_\alpha
\eqno(2.4)
$$
we require the function $F(t)$ to satisfy 
$$
{\cal L}_E F(t) :=\sum_{\alpha=1}^n \left[ (1-q_\alpha ) t^\alpha +
r_\alpha\right]  
\partial_\alpha F(t) = (3-d) F(t) + {1\over 2} A_{\alpha\beta} t^\alpha
t^\beta + B_\alpha t^\alpha +C
\eqno(WDVV3)
$$
for some constants $A_{\alpha\beta}$, $B_\alpha$, $C$.
The numbers $q_\alpha$, $r_\alpha$, $d$ must satisfy the following
normalization conditions
$$
q_1=0, ~~r_\alpha\neq 0 ~{\rm only ~if}~ q_\alpha=1.
\eqno(2.5)
$$
Loosely speaking, we assign the degree $1-q_\alpha$ to the variable
$t^\alpha$. But if $q_\alpha =1$ then the degree $r_\alpha$ is assigned to
$\exp t^\alpha$. With respect to this assignment the function $F(t)$ has
degree $3-d$ up to quadratic terms.

We will consider the class of equivalence of solutions modulo additions of
quadratic polynomials in $t$.
\smallskip
{\bf Exercise 2.1.} For any $\alpha$, $\beta$ prove that
$$
(q_\alpha+q_\beta -d)\eta_{\alpha\beta} =0.
\eqno(2.6)
$$
\smallskip
{\bf Exercise 2.2.} Prove that, by adding a quadratic polynomial to
$F(t)$,
the coefficients $A_{\alpha\beta}$, $B_\alpha$, $C$ in (WDVV3) can be
normalized in such a way that
$$
\eqalign{A_{\alpha\beta} &\neq 0 ~{\rm only ~if}~q_\alpha+q_\beta =
d-1\cr
A_{1\alpha} &=\sum_\alpha \eta_{\alpha\epsilon}r_\epsilon\cr
B_\alpha &\neq 0  ~{\rm only ~if}~q_\alpha=d-2\cr
B_1 &=0\cr
C &\neq 0  ~{\rm only ~if}~d=3.\cr}
\eqno(2.7)
$$
\medskip
Normalized in such a way coefficients $A_{\alpha\beta}$, $B_\alpha$, $C$
must also be considered as the unknown parameters of the WDVV problem.

Trivial solutions are cubics corresponding to graded Frobenius algebras
$\left( A_0, <~,~>\right)$
$$
{\rm cubic}\, ={1\over 6} c_{\alpha\beta\gamma} t^\alpha t^\beta
t^\gamma = {1\over 6} <1, (t)^3>.
\eqno(2.8)
$$
Solutions needed are analytic perturbations of cubics. That means that
$$
F(t) = {\rm cubic} + \sum_{k,l\geq 0} a_{k,l} \left( t''\right)^l e^{k
t'}
\eqno(2.9)
$$
where the vector argument $t$ is subdivided in two parts
$$
t=\left( t', t'' \right), ~\deg t'=0, ~\deg t'' \neq 0,
\eqno(2.10)
$$
$k$, $l$ are multiindices with all nonnegative coordinates. For
$$
t''\to 0, ~~t'\to -\infty
\eqno(2.11)
$$
$F(t)$ goes to the cubic. In quantum cohomology this is called the point
of classical limit.

There are two main approaches in the theory WDVV.

\smallskip
1. Algebraic approach. We study the formal series solutions (2.9) to WDVV
analyzing, say, the recursion relations for the coefficients $a_{k, l}$. 
An example of the algebraic approach are Kontsevich's recursion relations
for the numbers of plane curves, and also our discovery of Platonic solids
when classifying polynomial solutions to WDVV. A general approach to
construct solutions of WDVV in the class of formal series was recently
proposed in [BS]. Some family of formal power series solutions (not
satisfying the quasihomogeneity (WDVV3)) was very recently constructed in
[Los]. 
\smallskip
2. Analytic approach. First to describe {\it all} solutions to WDVV and
then to select the solutions of the needed class (2.9).
\medskip
In my lectures I will speak on the analytic approach to WDVV. This can be
applied to the solutions of the form (2.9) only if the series converges near
the point of classical limit (2.11). The convergence can be easily verified
in concrete examples of quantum cohomologies. However, the general proof
of convergence is still missing.

Let me first be more specific about the explicit form of WDVV.
\smallskip
{\bf Exercise 2.3.} Let $Q$ be the grading operator in ${\bf C}^n={\rm
span}\, (e_1, \dots, e_n)$ defined by 
$$Q(e_\alpha ) = q_\alpha e_\alpha, ~\alpha = 1, \dots, n.
$$
Show that WDVV remain invariant under the linear transformations
of the variables $t$
$$
\eqalign{t &\mapsto Mt, ~t= \left(t^1, \dots, t^n\right)^T\cr
(e_1, \dots, e_n) &\mapsto (e_1, \dots, e_n)M^{-1}\cr}
$$
(the upper label ${}^T$ stands for the transposition)
where the matrix $M$ satisfies the two conditions
$$\eqalign{Me_1 &=e_1\cr
MQ &= QM.\cr}
$$
Prove that, if $d\neq 0$ then the matrix
$\eta=\left(\eta_{\alpha\beta}\right)$
by a transformation of the above form can be reduced to the antidiagonal
form
$$
\eta_{\alpha\beta}=\delta_{\alpha+\beta, \, n+1}.
\eqno(2.12)
$$
Derive from (WDVV2) that, in these coordinates, the function $F(t)$ can
be represented in the form
$$
F(t) = {1\over 2} {t^1}^2 t^n +{1\over 2} t^1 \sum_{\alpha=2}^{n-1}
t^\alpha t^{n-\alpha +1} + f\left( t^2, \dots, t^n\right)
\eqno(2.13)
$$
for some function $f$ of $n-1$ variables. WDVV can be written as a system
of differential equations for this function. 
\medskip
We will usually consider only the case $d\neq 0$, although this is not
important for the mathematical theory of WDVV.
\smallskip
{\bf Example 2.1.} $n=2$. Here $f=f(t_2)$. The equations WDVV1 are empty.
The quasihomogeneity condition WDVV3 gives that
$$\eqalign{f &= {t_2}^{{3-d\over 1-d}}, ~d\neq 1\cr
f &= e^{{2 t_2 \over r}}, ~ d=1, ~ E=t_1 \partial_1 + r \partial_2,\cr
f &= -{1\over 2} c \log t_2, ~ d=3, ~ E= t_1 \partial_1 - 2 t_2
\partial_2, ~ {\cal L}_E F = F+c.\cr}
$$
\medskip
{\bf Example 2.2.} For $n=3$ the function $f=f(x,y)$, $x=t_2$, $y=t_3$
must satisfy the following PDE
$$f^2_{xxy}=f_{yyy}+f_{xxx} f_{xyy}.
$$
For generic $d$ the variables $t_1$, $t_2$, $t_3$ have the scaling
dimensions $1$, $1-{d\over 2}$, $1-d$ respectively, and the scaling
dimension of the function $f$ is $3-d$. So
$$f(x,y) = {x^4\over y} \phi \left( \log (y x^q)\right), ~~q=2{d-1\over
2-d}.
$$
Plugging this to the above PDE one obtains the following complicated 3d order
ODE for the function $\phi$
$$-6\,\phi  + 48\,{{\phi }^2} + 11\,\phi ' +
  88\,q\,\phi \,\phi ' - (144 +
  144\,q - 3\,{q^2})\,{{\phi '}^2} -
  6\,\phi '' +48 (2 +
  2\,q +
  {q^2})\,\phi \,\phi '' 
$$
$$-4\,q(16+
16\,{q} +
  {q^2})\,\phi '\,\phi '' -
  (13\,{q^2} + 13\,{q^3} +
  {q^4})\,{{\phi ''}^2}$$
$$ + \phi '''+ 
  8q(3 + 3q +q^2)\,\phi \,\phi ''' +
  2\,{q^2}(1+q+q^2)\,\phi '\,\phi ''' -
  {q^3}(1 + q)\,\phi ''\,\phi '''=0.
$$
The nongeneric are the integer values $-2 \leq d \leq 4$. In this case
the ansatz for $f$ is to be modified. For example, for $d=2$, $r_2=r$
$$
f(x,y) = {1\over y} \phi(x+r\log y)
$$
where the function $f$ satisfies the ODE
$$
\phi''' [r^3 + 2\phi' - r \phi''] -(\phi'')^2 - 6r^2 \phi''
+ 11 r \phi' - 6\phi = 0.
\eqno(2.14)
$$
The case of quantum cohomology of ${\bf CP^2}$ corresponds to $r=3$
(see Lecture 1 above). In this case the equation (2.14) has a unique solution
$\phi = \phi(x)$ of the form
$$
\phi = \sum_{k\geq 1}A_k e^{k\, x}
\eqno(2.15)
$$
normalized by the condition $A_1={1\over 2}$. Plugging the series (2.15) 
into the equation (2.14) one obtains the recursion relations for the numbers
$$
N_k = (3k-1)! A_k
$$ 
of rational curves of degree $k$ on  ${\bf CP^2}$ passing through $3k-1$ 
generic points.
\smallskip
{\bf Exercise 2.4.} Derive from the recursion relations that the series
(2.15)
converges if
$$
{
\rm Re}\, x< \log {6\over 5}.
$$
\medskip
Recall that the numerical estimate of [DI] guarantees convergence for 
$$
{\rm Re}\, x< 1.981.
$$

We conclude that in the first nontrivial case $n=3$ the general solution of 
WDVV depends on 3 arbitrary parameters. However, this parametrization does 
not say anything about the analytic properties of solutions. For the next 
case $n=4$
the situation looks to be even worse: the function $f=f(t_2, t_3, t_4)$ is to 
be found from an overdetermined system of 6 PDEs. With $n$ growing the 
overdeterminancy of the system WDVV1 grows rapidly.

In these lectures we will give complete classification of the solutions of 
WDVV satisfying the semisimplicity condition. Recall that this condition
means that the algebra $A_t$ is semisimple for generic $t$. The solution will 
be expressed via certain Painlev\'e-type transcendents (via particular 
transcendents of the Painlev\'e-VI type in the first nontrivial case $n=3$).

In this lecture we will develop some preliminary geometric constructions of 
the theory of WDVV. First we will give a coordinate-free reformulation of \wdvv.
The basic idea is to identify the algebra $A_t$ with the tangent space $T_tM$ 
to the space of the parameters $t\in M$,
$$
A_t \ni e_\alpha \leftrightarrow \partial_\alpha \in T_t M, ~\alpha=1,
\dots, n.
$$
The space of parameters $M$ acquires a new geometrical structure: the tangent 
spaces $T_t M$ are Frobenius algebras w.r.t. the multiplication
$$
\partial_\alpha \cdot \partial_\beta = c_{\alpha\beta}^\gamma(t) 
\partial_\gamma
\eqno(2.16)
$$
and metric
$$
<\partial_\alpha, \partial_\beta> =\eta_{\alpha\beta}
\eqno(2.17)
$$
We arrive [Du5] at the following main
\smallskip
{\bf Definition 2.1.} (Smooth, analytic) {\it Frobenius structure} 
on the manifold $M$ is a structure of Frobenius algebra on the tangent spaces 
$T_tM=\left(A_t, <~,~>_t\right)$
depending (smoothly, analytically) on the point $t$. This structure must 
satisfy the following axioms.
\smallskip
\noindent{\bf FM1}. The metric on $M$ induced by the invariant bilinear form
$<~,~>_t$ is flat. Denote $\nabla$ the Levi-Civita connection for the metric
$<~,~>_t$. The unity vector field $e$ must be covariantly constant,
$$\nabla \, e=0.
\eqno(2.18)
$$
We use here the word `metric' as a synonim of a symmetric nondegenerate
bilinear form on $TM$, not necessarily a positive one. Flatness of the metric,
i.e., vanishing of the Riemann curvature tensor, means that locally a
{\it system of flat coordinates} $(t^1, \dots, t^n)$ exists such that the 
matrix $<\partial_\alpha , \partial_\beta>$
of the metric in these coordinates becomes constant.
\medskip
\noindent{\bf FM2}. Let $c$ be the following symmetric trilinear form on $TM$
$$c(u,v,w) := <u\cdot v, w>.
\eqno(2.19)
$$
The four-linear form
$$
(\nabla_z c)(u,v,w), ~u,v,w,z\in TM
$$
must be also symmetric.
\medskip
Before formulating the last axiom we observe that the space $Vect(M)$ of 
vector fields on $M$ acquires a structure of a Frobenius algebra over the 
algebra
$Func(M)$ of (smooth, analytic) functions  on $M$. 
\smallskip
\noindent{\bf FM3}. A linear vector field $E\in Vect(M)$ must be fixed on $M$,
i.e.,
$$
\nabla\nabla \, E =0.
\eqno(2.20)
$$
The operators
$$
\eqalign{
Q_{Func(M)}&:=E\cr
Q_{Vect(M)}&:={\rm id} + {\rm ad}_E\cr}
\eqno(2.21)
$$
introduce in $Vect(M)$ a structure of graded Frobenius algebra  of a
given charge $d$ over the graded
ring $Func(M)$ (see above Remark 1.1 after Definition 1.2).
\medskip
{\bf Lemma 2.1.} {\it Locally a Frobenius manifold with diagonalizable 
$\nabla\,E$  is described by a solution of WDVV and vice versa.}
\smallskip
{\bf Proof.} 1. Starting from a solution of WDVV define the multiplication
(2.16) and the metric (2.17) on the tangent planes to the parameter space. 
In the original coordinates $(t^1 \dots, t^n)$ the metric is manifestly flat. 
In these coordinates the covariant derivatives coincide with the partial ones
$$
\nabla_\alpha =\partial_\alpha.
$$
Since $e=\partial_1$ we have $\nabla\,e \equiv 0$. 
The first axiom FM1 is proved. The tensor $c$ in (2.19) has the components
$$
c_{\alpha\beta\gamma} (t) \equiv 
c(\partial_\alpha, \partial_\beta, \partial_\gamma ) 
=\partial_\alpha \partial_\beta \partial_\gamma F(t).
$$
So
$$
\left( \nabla_{{\partial_\delta}} c\right) 
(\partial_\alpha, \partial_\beta, \partial_\gamma) = 
\partial_\alpha \partial_\beta \partial_\gamma \partial_\delta
F(t)
$$
is totally symmetric. This proves FM2.

Let us now prove FM3. The equations
$$
\eqalign{Q_{Vect(M)} (a\cdot b) &=Q_{Vect(M)}(a)\cdot b + a\cdot 
Q_{Vect(M)}(b)
\cr
Q_{Func(M)}<a,b> + d<a,b>
&=<Q_{Vect(M)}(a),b>+<a,Q_{Vect(M)}(b)> \cr}
$$
can be recasted in the form
$$
\eqalignno{\L (a\cdot b) - \L (a) \cdot b - a\cdot \L (b) &= a\cdot b
& {(2.22)}\cr
\L <a,b> -<\L a,b> -<a, \L b> &= (2-d) <a,b>.
& {(2.23)}
\cr}
$$

We will prove the last two equations.

The Euler vector field is clearly a linear one. The gradient $\nabla\, E$
is a diagonal constant matrix
$$
\nabla\, E= {\rm diag}\, \left (1-q_1, \dots, 1-q_n\right).
\eqno(2.24)
$$
Triple differentiating of the quasihomogeneity equation (WDVV3) w.r.t. 
$t^\alpha$, $t^\beta$, $t^\gamma$ gives $$
\sum_\epsilon \left[ (1-q_\epsilon)t^\epsilon + r_\epsilon \right]
\partial_\epsilon \left( c_{\alpha\beta\gamma} (t)\right) = (q_\alpha + q_\beta
+q_\gamma -d) c_{\alpha\beta\gamma} (t).
\eqno(2.25)
$$
From this and from (2.6)
easily follow the identities (2.22), (2.23) of the definition
of the graded
Frobenius algebra over graded ring of functions. 

2. Choose locally flat coordinates $(t^1, \dots, t^n)$ on a Frobenius manifold.
We can choose them in a particular way such that $\partial_1$, \dots,
$\partial_n$ are the eigenvectors of the linear operator $\nabla\,E:TM\to TM$
$$
(\nabla\,E) \partial_\alpha = \lambda_\alpha \partial_\alpha
$$
for some constant $\lambda_\alpha$ (in the flat coordinates the matrix of
the 
covariantly constant tensor $\nabla\, E$ is constant). This will be the
homogeneous basis for the grading operator $Q_{Vect(M)}$ 
$$
Q_{Vect(M)} \partial_\alpha =(1-\lambda_\alpha) \partial_\alpha.
$$
So
$$
E=\sum_{\alpha=1}^n (\lambda_\alpha t^\alpha + r_\alpha)\partial_\alpha
$$
for some constants $r_\alpha$. We can kill by a shift all these constants
but those for which $\lambda_\alpha=0$. This gives the form (2.4) of the
Euler vector field with 
$$
q_\alpha:= 1-\lambda_\alpha.
$$
From the obvious equation
$$
Q_{Vect(M)} e =0
$$
we immediately obtain that $\lambda_1=1$, i.e., $q_1=0$.

From the symmetry w.r.t. $\alpha$, $\beta$, $\gamma$, $\delta$ of
partial derivatives
$$
\partial_\delta c_{\alpha\beta\gamma} (t) 
$$
of the symmetric tensor
$$
c_{\alpha\beta\gamma} (t) = <\partial_\alpha \cdot \partial_\beta, \partial_\gamma>
$$
we conclude that locally a function $F(t)$ exists such that
$$
c_{\alpha\beta\gamma}(t) =\partial_\alpha\partial_\beta\partial_\gamma F(t).
$$
For the invariant metric we obtain
$$
\eta_{\alpha\beta} =<\partial_\alpha, \partial_\beta> = <\partial_\alpha\cdot
\partial_\beta, \partial_1> =\partial_1\partial_\alpha\partial_\beta F(t).
$$
We have proved WDVV1 and WDVV2.

Spelling the last axiom FM3 out we obtain the following two formulae
$$
\eqalign{ {\cal L}_E \eta_{\alpha\beta} &= \partial_\alpha E^\epsilon
\eta_{\epsilon\beta} 
+ \partial_\beta E^\epsilon \eta_{\alpha\epsilon}
=(2-d)\eta_{\alpha\beta}\cr
{\cal L}_E c_{\alpha\beta}^\gamma &=
E^\epsilon \partial_\epsilon c_{\alpha\beta}^\gamma - \partial_\epsilon 
E^\gamma c_{\alpha\beta}^\epsilon +\partial_\alpha E^\epsilon c_{\epsilon\beta}^\gamma
+\partial_\beta E^\epsilon c_{\alpha\epsilon}^\gamma = c_{\alpha\beta}^\gamma
\cr}
$$
From this it follows
$$
(q_\alpha+q_\beta-d) \eta_{\alpha\beta}=0
\eqno(2.26)
$$
$$
E^\epsilon \partial_\epsilon c_{\alpha\beta}^\gamma
=(q_\alpha+q_\beta-q_\gamma) c_{\alpha\beta}^\gamma.
\eqno(2.27)
$$
Using (2.6) we lower the index $\gamma$ in the last equation to obtain
$$
E^\epsilon \partial_\epsilon
c_{\alpha\beta\gamma}=(q_\alpha+q_\beta+q_\gamma -d)c_{\alpha\beta\gamma}, 
~~\alpha, \, \beta, \, \gamma =1, \dots, n.
$$
Triple integration gives 
$$
E^\epsilon \partial_\epsilon F = (3-d) F + ~{\rm quadratic}.
$$
Lemma is proved.
\smallskip
{\bf Remark 2.1}. The definition of Frobenius manifold can be easily
translated into algebraic language as a graded Frobenius algebra structure
on the module of derivations of a graded commutative associative algebra.
An important extension of this definition for the case of ${\bf
Z}_2$-graded algebras was done by Kontsevich and Manin [KM1]. Such {\it
Frobenius supermanifolds}
are necessary to deal with Gromov - Witten invariants of
manifolds with nontrivial odd-dimensional cohomologies. In this paper we
will not discuss this extension. 
\medskip
{\bf Exercise 2.5.} Prove that the direct product $M'\times M''$ of 
two Frobenius manifolds of {\it the same} charge $d$ carries a natural
structure of a Frobenius manifold of the charge $d$, the unity vector 
field $e'\oplus e''$, and the Euler vector field $E'\oplus E''$.
\medskip
We now address the problem of (local) classification of Frobenius manifolds
coinciding with local clasification of solutions of WDVV. To be more specific 
we give
\smallskip
{\bf Definition 2.2.} A (local) diffeomorphism 
$$
\phi: M\to \tilde M
$$
of two Frobenius manifolds is called {\it (local) equivalence} if
the differential
$$
\phi_* : T_t M\to T_{\phi(t)}\tilde M
$$
is an isomorphism of algebras for any $t\in M$ and
$$
\phi^* <~,~>_{\tilde M} = c^2 <~,~>_M
$$
where $c$ is a nonzero constant not depending on the point of $M$. 
\medskip
The corresponding free energies $F$ and $\tilde F$ are related by
$$
\tilde F\left( \phi(t)\right) = c^2 F(t) + {\rm quadratic}.
$$
\smallskip
{\bf Definition 2.3.} A Frobenius manifold is called {\it reducible} if it 
is equivalent to the direct product of two Frobenius manifolds (see
Exercise 2.5 above).
\medskip
The first main tool in dealing with Frobenius manifolds is a deformation
of the Levi-Civita connection $\nabla$. We put
$$
\tilde\nabla_u v:= \nabla_u v + z\, u\cdot v.
\eqno(2.28a)
$$
Here $u$, $v$ are two vector fields on $M$, $z$ is the parameter of the deformation. We extend this up to a meromorphic connection on the direct product
$M\times {\bf C}$ by the formulae
$$
\eqalign{
\tilde\nabla_u {d\over dz} &=0\cr
\tilde\nabla_{{d\over dz}} {d\over dz} &=0\cr
\tilde\nabla_{{d\over dz}} v &= \partial_z v +E\cdot v - {1\over z}\mu\,
v\cr}
\eqno(2.28b)
$$
where
$$
\mu:= {2-d\over 2} - \nabla E ={\rm diag}\, (\mu_1, \dots, \mu_n)
\eqno(2.29a)
$$
$$
\mu_\alpha := q_\alpha -{d\over 2}.
\eqno(2.29b)
$$ 
Here $u$, $v$ are tangent vector fields on $M\times {\bf C}$ having zero
component along ${\bf C}$. Observe that $\tilde\nabla$ is a symmetric connection.

\smallskip
{\bf Proposition 2.1.} {\it For a Frobenius manifold $M$ the curvature of
the
conection $\tilde\nabla$ equals zero. Conversely, if on the tangent spaces
to $M$ a structure of a Frobenius algebra is defined satisfying FM1,
and the curvature of the connection $\tilde\nabla$ vanishes and the Euler
vector field $E$ satisfies
$$
{\cal L}_E <~,~> = (2-d) <~,~>
\eqno(2.30)
$$
with a constant $d$, then $M$ is a Frobenius manifold.}

{\bf Proof.} For a covector
$$
\xi =\xi_\alpha dt^\alpha + 0\, dz
$$
one has
$$
\eqalign{\tilde\nabla _\alpha \xi_\beta &= \partial_\alpha\xi_\beta -z\,
c_{\alpha\beta}^\gamma \xi_\gamma\cr
\tilde\nabla_{{d\over dz}} \xi_\beta &= \partial_z\xi_\beta -E^\gamma
c_{\gamma\beta}^\alpha\xi_\alpha +{1\over z} M_\beta^\epsilon
\xi_\epsilon\cr}
$$
where we 
denote
$M_\beta^\epsilon$
the matrix entries of the linear operator $\mu= {1\over 2}(2-d) -
\nabla\,E$. Any solution of
the system $\tilde\nabla \xi =0$ is a (local) horizontal section of
$T^* \left( M\times {\bf C}\right)$ for the connection $\tilde\nabla$. A
basis of horizontal sections is given by $dz$ and by $n$ linearly
independent solutions of the system
$$
\eqalignno{ \partial_\alpha\xi_\beta &=
z\,c_{\alpha\beta}^\gamma \xi_\gamma
& {(2.31)}\cr
\partial_z\xi_\beta &=
E^\gamma
c_{\gamma\beta}^\alpha\xi_\alpha -{1\over z} M_\beta^\epsilon
\xi_\epsilon.
& {(2.32)} \cr}
$$
Such a basis exists {\it iff} the compatibility conditions 
$$
\partial_\alpha\partial_\gamma = \partial_\gamma\partial_\alpha, ~
\partial_\alpha\partial_z = \partial_z\partial_\alpha
$$
hold true. Differentiating (2.31) w.r.t. $t^\gamma$ and subtracting the same
expression with $\alpha$ and $\gamma$ permuted we obtain the first
compatibility condition in the form
$$
z\, \left( \partial_\gamma c_{\alpha\beta}^\epsilon -\partial_\alpha
c_{\gamma\beta}^\epsilon \right) \xi_\epsilon
+ z^2 \left( c_{\alpha\beta}^\lambda c_{\lambda\gamma}^\epsilon -
c_{\gamma\beta}^\lambda c_{\lambda\alpha}^\epsilon\right) \xi_\epsilon =0.
$$
This must vanish for arbitrary $\xi$. We obtain
$$ c_{\alpha\beta}^\lambda c_{\lambda\gamma}^\epsilon =
c_{\gamma\beta}^\lambda c_{\lambda\alpha}^\epsilon
$$
(associativity) and
$$ \partial_\gamma c_{\alpha\beta}^\epsilon =\partial_\alpha
c_{\gamma\beta}^\epsilon
$$
(local existence of $F(t)$). Similarly, from compatibility of (2.31) and 
(2.32) we obtain, first, that 
$$
\partial_\alpha M_\beta^\epsilon =0.
$$
So $\nabla\,E$ is a constant matrix. Assume, for simplicity,  $\nabla\,E$
to be diagonal,  $\nabla\,E = {\rm diag}\, (1-q_1, \dots, 1-q_n)$. Then
we further obtain 
$$
E^\epsilon \partial_\epsilon c_{\alpha\beta}^\gamma =
(q_\alpha+q_\beta-q_\gamma) c_{\alpha\beta}^\gamma.
$$
As we already know, this together with (2.30) (i.e., with (2.6))  is 
equivalent to FM3. Proposition is proved.
\medskip
{\bf Remark 2.2.} Due to Proposition, one can alternatively define
Frobenius manifolds as those carrying a metric and a linear pencil
of affine connections (2.28a) satisfying the above conditions
(such a definition was explicitly used in [Du4]). It is interesting that
manifolds with a metric and a linear pencil of affine connections
deforming the Levi-Civita connection are known also in mathematical
statistics - see the book [Ama]. This structure appeares in 
the parametric statistics that studies families of probabilistic measures
depending on a finite number of parameters. The metric was introduced by
Rao about 1945 using the classical Fischer matrix of the family. The
deformed Levi-Civita connection was discovered by N.N.Chentsov in 1972.
It has the form (2.28a). However, the curvature of the deformed connection
does not vanish identically but it vanishes for two values of the
parameter $z$.

\medskip
{\bf Exercise 2.6.} Prove that the solutions of the linear system 
(2.31), (2.32) are all
closed differential forms
$$
\xi_\alpha dt^\alpha =d\, \tilde t
$$
(the differential along $M$ only).
\smallskip
Choosing a basis of $n$ linearly independent solutions $\xi_\alpha^{(1)}, 
\dots, \xi_\alpha^{(n)}$ of the system
we obtain $n$ functions $\tilde t_1(t,z), \dots, \tilde t_n(t,z)$.
Together with $z$ they give a system of flat coordinates for the
connection $\tilde \nabla$ on a domain in $M\times {\bf C}$. This means
that, in these coordinates, the covariant derivatives coincide with
partial ones.
\medskip
How to choose a basis of the deformed flat coordinates  $\tilde t_1(t,z),
\dots, \tilde t_n(t,z)$? Let us forget first about the last component
(2.28b) of the connection $\tilde\nabla$. The first part (2.28a) can be 
considered
as a deformation of the affine structure on $M$ with $z$ being the
parameter
of the deformation. We can look for the deformed flat coordinates in the
form of the series
$$
\tilde t_\alpha = \sum_{p=0}^\infty h_{\alpha,p}(t) z^p =:h_\alpha(t;z),
~\alpha=1, \dots,n.
\eqno(2.33)
$$
\smallskip
{\bf Lemma 2.2.} {\it The coefficients $h_{\alpha,p}(t)$ can be determined
recursively from the relations
$$
\eqalign{h_{\alpha,0} &=t_\alpha \equiv \eta_{\alpha\epsilon}t^\epsilon\cr
\partial_\beta\partial_\gamma h_{\alpha, p+1} &= c_{\beta\gamma}^\epsilon
\partial_\epsilon h_{\alpha, p}, ~p=0, \, 1, \, 2, \dots\cr}
\eqno(2.34)
$$
uniquely up to a transformation of the form
$$
\tilde t_\alpha \mapsto \sum_{\beta=1}^n \tilde t_\beta G_\alpha^\beta (z)
$$
where the coefficients $G_1$, $G_2$, \dots of the matrix-valued series
$$
G(z) =\left( G_\alpha^\beta (z)\right) = 1 + z\, G_1 + z^2 G_2 + \dots
$$
do not depend on $t$.}

{\bf Proof.} We are only to show that the right-hand sides of (2.34) are
second derivatives along $t^\beta$ and $t^\gamma$. This can be proved
inductively using the identity
$$
\partial_\alpha\left( c_{\beta\gamma}^\epsilon \partial_\epsilon 
h_{\lambda,p}\right) -\partial_\beta\left( c_{\alpha\gamma}^\epsilon
\partial_\epsilon
h_{\lambda,p}\right)=
\left( c_{\beta\gamma}^\epsilon c_{\epsilon\alpha}^\rho -
 c_{\alpha\gamma}^\epsilon c_{\epsilon\beta}^\rho\right) \partial_\rho
h_{\lambda, p-1} =0
$$
due to associativity.
\medskip
The gradients $\nabla\,h_{\alpha,p}(t)$ and their inner products
$<\nabla\, h_{\alpha, p}, \nabla\, h_{\beta, q}>$ play a very important
role in the theory and applications of Frobenius manifolds. Before
discussing how to normalize them uniquely I will give here two important
identities for these coefficients.
\smallskip
{\bf Exercise 2.7.} Prove that
$$
t_\alpha = < \nabla\, h_{\alpha,0}, \nabla \, h_{1,1}>
\eqno(2.35)
$$
$$
F(t)
 ={1\over 2} \left\{ 
<\nabla\, h_{\alpha,1},\nabla\, h_{1,1}>\eta^{\alpha\beta}
<\nabla\, h_{\beta,0}, \nabla\, h_{1,1}>\right.
$$
$$
\left. -<\nabla\, h_{1,1},\nabla\,h_{1,2}> 
 -<\nabla\, h_{1,3},\nabla\,h_{1,0}>\right\} .
\eqno(2.36)
$$
\smallskip
{\bf Exercise 2.8.} Prove the identity
$$
\nabla \left< \nabla \, h_\alpha (t;z), \nabla\, h_\beta (t;w)\right>
=(z+w) \nabla\, h_\alpha (t;z) \cdot \nabla\, h_\beta(t;w).
$$

Observe that $\left< \nabla\, h_\alpha (t;z), \nabla\,
h_\beta(t;-z)\right>$ does not depend on $t$.
\medskip
To choose canonically the system of deformed flat coordinates
$\tilde t_1(t;z), \dots, \tilde t_n(t;z)$ we will now use the last
equation (2.32) of the horizontality of the gradients $\xi_\alpha
=\partial_\alpha \tilde t(t;z)$
$$
\partial_z\xi_\alpha =\U_\alpha^\beta \xi_\beta -{1\over z} \mu_\alpha
\xi_\alpha.
\eqno(2.37)
$$
Here
$$
\U_\alpha^\beta (t) := E^\epsilon c_{\epsilon\alpha}^\beta
\eqno(2.38)
$$
is the matrix of multiplication by the Euler vector field. The choice of
the basis can be done by carefully looking at the behaviour of the
solutions at $z=0$ where the connection $\tilde\nabla$ has logarithmic
singularity. The analysis of this behaviour will provide us with some
numerical invariants of the Frobenius manifold.

Let us introduce the numbers
$$
\mu_\alpha =q_\alpha -{d\over 2}, ~\alpha=1,\dots, n.
$$
Recall that they are the entries of the diagonal matrix
$$
\mu={2-d\over 2} 1 -\nabla\, E ={\rm diag}\, (\mu_1, \dots, \mu_n).
\eqno(2.39)
$$
The operator $\mu$ is antisymmetric w.r.t. the metric $<~,~>$
$$
<\mu\, a, b> +<a,\mu\,b> =0.
\eqno(2.40)
$$
We say that the operator $\mu$ is {\it resonant} if some of the
differences $\mu_\alpha-\mu_\beta$ is a nonzero integer. Otherwise it is
called {\it nonresonant}. We will also use expressions
resonant/nonresonant Frobenius manifold if the correspondent operator
$\mu$ is resonant/nonresonant. For example, any Frobenius manifold 
related to quantum cohomology is resonant (all the numbers $q_\alpha$ are
integers). The Frobenius manifold on the space of versal deformations
of $A_3$ singularity is nonresonant one
$$
\mu={\rm diag}\, \left(-{1\over 4}, 0, {1\over 4} \right).
$$

We consider first the nonresonant case. In this case the system of the
deformed flat coordinates can be uniquely chosen in such a way that
$$
\tilde t_\alpha (t;z) =\left[ t_\alpha +O(z)\right]\, z^{\mu_\alpha}, ~
\alpha=1, \dots, n
$$
(here, as usual, $t_\alpha = \eta_{\alpha\epsilon} t^\epsilon$). The
coordinates are multivalued analytic functions of $z$ defined for
sufficiently small
$z\neq
0$. Going along a closed loop around $z=0$ one obtains the monodromy
transformation
$$
\left( \tilde t_1, \dots, \tilde t_n\right) \mapsto
\left( \tilde t_1, \dots, \tilde t_n\right)\, M_0,
$$
$$
M_0 ={\rm diag}\, \left( e^{2\pi i\mu_1}, \dots, e^{2\pi i \mu_n }\right).
$$

For  the resonant case such a choice is not possible. The monodromy matrix
$M_0$ cannot be diagonalized. 

We will first rewrite the equation (2.32) in
the matrix form. Doing the linear change
$$
\xi=(\xi_1, \dots, \xi_n) \mapsto \eta^{-1} \xi^T
$$
rewrite (2.31), (2.32) as follows
$$
\eqalignno{
\partial_\alpha \xi & = z\, C_\alpha \xi
& {(2.41a)}
\cr
\partial_z \xi & =\left( \U+{1\over z} \mu \right) \xi
& {(2.41b)}
\cr}
$$
where
$$
\left( C_\alpha\right) _\beta^\gamma = c_{\alpha\beta}^\gamma.
$$
The matrix $\U$ is $\eta$-symmetric
$$
\U^T\eta = \eta\, \U
\eqno(2.42)
$$
and $\mu$ is $\eta$-antisymmetric
$$
\mu\,\eta +\eta\,\mu=0.
\eqno(2.43)
$$
The solutions of the system (2.41) are gradients of the deformed flat
coordinates $\xi = \nabla\, \tilde t$.
\smallskip
{\bf Lemma 2.3.} {\it The bilinear form
$$
\left< \xi_1, \xi_2\right> := \xi_1^T (-z) \, \eta\, \xi_2(z)
\eqno(2.44)
$$
on the space of solutions of (2.41b) does not depend on $z$.}

Proof is obvious.
\medskip
Let us study the classes of equivalence of the system (2.41b) under gauge
transforms
$$
\xi \mapsto G(z) \, \xi
\eqno(2.45a)
$$
of the form
$$
G(z) = 1+z\, G_1 + z^2 G_2 + \dots
\eqno(2.45b)
$$
$$
G^T(-z) \eta G(z) \equiv \eta.
\eqno(2.45c)
$$
\smallskip
{\bf Lemma 2.4.} {\it After an arbitrary gauge transform (2.45) the
vector-function
$$
\tilde \xi =G(z)\xi
\eqno(2.46)
$$
satisfies the system
$$
\partial_z \tilde \xi =\left({1\over z} \mu + \tilde U_1 + z\, \tilde U_2
+ \dots \right) \tilde\xi
$$
where the matrices $\tilde U_{2k+1}$ are $\eta$-symmetric and the matrices
$\tilde U_{2k}$ are $\eta$-antisymmetric.}

{\bf Proof.} After the gauge transform the vector function $\tilde \xi$
satisfies
$$
\partial_z \tilde\xi =A(z)\tilde\xi
$$
with
$$
\eqalign{
A(z) &= G(z) \left( {1\over z} \mu + \U\right) G^{-1}(z) + G'(z)
G^{-1}(z)\cr
 &=: {1\over z}\mu + \tilde U_1 + z\, \tilde U_2 +\dots\cr}
$$
where the matrix coefficients $\tilde U_k$ are defined by this equation.
Using (2.42), (2.43) and (2.45c) one obtains
$$
A^T(-z) =\eta\, A(z)\,\eta^{-1}.
$$
This gives
$$
\tilde U_k^T =(-1)^{k+1} \eta\, \tilde U_k \, \eta^{-1}, ~k=1,\, 2, \dots.
$$
Lemma is proved.
\medskip
In the nonresonant case one can choose the gauge transform in such a way
that all $\tilde U_1 = \tilde U_2 =\dots =0$. The only gauge invariant of
the system (2.41b) near the logarithmic singularity $z=0$ is the diagonal
matrix $\mu$. 

Let us consider slightly more general system
$$
\partial_z \xi =\left( {1\over z} \mu +U_1 +z\, U_2 +z^2\,
U_3+\dots\right)\xi
\eqno(2.47)
$$
with the coefficients satisfying
$$
U_k^T =(-1)^{k+1} \eta\,  U_k \, \eta^{-1}, ~k=1,\, 2, \dots.
\eqno(2.48)
$$
\smallskip
{\bf Lemma 2.5.} {\it By a gauge transformation of the form (2.45) the system
(2.47) can be reduced to the canonical form
$$
\partial_z \tilde \xi =\left( {1\over z} \mu + R_1 + z\, R_2 + z^2\, R_3
+\dots \right)\tilde\xi
\eqno(2.49)
$$
where the matrices $R_1$, $R_2$, \dots satisfy
$$
\eqalignno{
R_k^T &=(-1)^{k+1} \eta\, R_k \, \eta^{-1}
& {(2.50)}\cr
\left(R_k\right)_\beta^\alpha &\neq 0 ~{\rm only
~if}~\mu_\alpha-\mu_\beta=k , ~k=1,\, 2, \dots.
& {(2.51)}\cr}
$$
}

Observe that there is only a finite number of nonzero matrices $R_k$.

{\bf Proof.} Gauge equivalence of the system (2.47) to a system (2.49) 
with the
matrices $R_k$ satisfying (2.51) is a wellknown fact (see, e.g. [Ga]).
Namely,
from the recursion relations
$$
R_n =U_n + n\, G_n +\left[ G_n, \mu\right] +\sum_{k=1}^{n-1} \left(G_{n-k}
U_k -R_k G_{n-k}\right)
\eqno(2.52)$$
we determine uniquely the matrix entries
$$
\left( R_n\right)_\beta^\alpha ~{\rm for}~ \mu_\alpha -\mu_\beta =n
$$
and
$$
\left(G_n\right)_\beta^\alpha ~{\rm for}~ \mu_\alpha -\mu_\beta \neq n
$$
and we put
$$
\left(G_n\right)_\beta^\alpha=0  ~{\rm for}~ \mu_\alpha -\mu_\beta =n.
$$
Using induction it can be easily seen that the matrices $R_n$ satisfy
the $\eta$-symmetry/antisym\-metry conditions (2.50) and the matrices $G_n$
satisfy the orthogonality conditions (2.53)
$$
G_n^T = (-1)^{n+1} \eta\, G_n\, \eta^{-1} + \sum_{k=1}^{n-1} (-1)^{n+k+1}
G_k^T \eta\, G_{n-k}\eta^{-1}.
\eqno(2.53)$$
Lemma is proved.

\medskip
We will call (2.49) {\it the normal form} of the system (2.47).
The ambiguity in the choice of the normal form  will be described below.
\smallskip
{\bf Lemma 2.6.} {\it The matrix solution of the system (2.49) is
$$
\xi = z^\mu z^R
\eqno(2.54)$$
where
$$
R:= R_1 +R_2 +\dots.
\eqno(2.55)
$$
}
{\bf Proof.} From (2.51) we obtain the identity
$$
z^\mu R_k z^{-\mu} =z^k R_k, ~k=1, \, 2, \dots.
\eqno(2.56)
$$
So, differentiating the matrix-valued function (2.54) one obtains
$$
\eqalign{\partial_z\xi &= {\mu\over z} z^\mu z^R + {1\over z} z^\mu
R z^R
\cr
 &= \left( {\mu\over z} + R_1 + R_2 z^2 + \dots \right) z^\mu z^R.\cr}
$$
Lemma is proved.
\medskip
{\bf Exercise 2.9.} Prove that the monodromy around $z=0$ of the solution
$$
\Xi_0 :=z^\mu z^R
$$
has the form
$$
\Xi_0 \left( z\, e^{2\pi i}\right) = \Xi_0 (z) M_0
$$
$$
M_0 = \exp 2\pi i (\mu +R).
\eqno(2.57)
$$
\medskip
We will now represent the parameters of the normal form (2.49) in a
geometric way. Let $\V$ be a linear space equipped with a symmetric
nondegenerate bilinear form $<~,~>$ and an antisymmetric operator
$$
\mu: \V\to \V, ~~<\mu\, a,b>+<a, \mu\,b>=0.
$$
Let us assume, for simplicity, the operator $\mu$ to be diagonalizable.
Let
$$
{\rm spec}\,:= \{ \mu_1, \dots, \mu_n\}
$$
be the spectrum of $\mu$. Denote $e_1$, \dots, $e_n$ the correspondent
eigenvectors. We define a filtration on $V$ 
$$
0=F_0 \subset F_1 \subset F_2 \subset \dots \subset \V
$$
$$
F_k :={\rm span}\, \{ e_\alpha | \mu_\alpha +k \not\in \, {\rm spec}\} .
\eqno(2.58)
$$
Obviously, for a non-resonant $\mu$ the filtration consists of two terms
$0=F_0 \subset F_1 =\V$. 

The asociated graded space
$$
\eqalign{
\V_* &= \oplus_{k\geq 1} \V_k \cr
\V_k &=F_k /F_{k-1}\cr}
$$
is isomorphic to $\V$ due to the natural isomorphism
$$
\V_k \simeq {\rm Ker}\, (\mu + k-1) \cap F_k \subset \V.
$$

A linear operator
$$
R: \V \to \V
$$
is called {\it $\mu$-nilpotent} if it commutes with $\exp 2\pi i \mu$ 
$$
R\, e^{2\pi i\mu}=e^{2\pi i\mu}R
\eqno(2.59)
$$
and
$$
R(F_k)\subset F_{k-1}, ~k=1,\, 2, \dots.
\eqno(2.60)
$$ 
The associated operator
$$
R_*:\V_*\to\V_*
$$
has a natural grading
$$
R_* =\oplus_{k\geq 1} R_k
\eqno(2.61)$$
where the operator $R_k$ shifts the grading by $-k$
$$
R_k \left( \V_m\right) \subset \V_{m-k} ~{\rm for~any}~m>k.
\eqno(2.62)
$$
Writing all the operators by matrices in the basis of eigenvectors
of $\mu$ one obtains
$$
\left(R_k\right)_\beta^\alpha =\cases{R_\beta^\alpha ~{\rm
if}~\mu_\alpha-\mu_\beta =k\cr
0 ~~{\rm otherwise.}\cr}
\eqno(2.63)
$$
We say that the $\mu$-nilpotent operator is {\it $\mu$-skew-symmetric} if
$$
\{ R\, x,y\} + \{ x, R\, y\} =0 ~{\rm for~any}~x,\, y \in\V
\eqno(2.64)
$$
where
$$
\{ x, y\} := \left< e^{\pi i \mu}x, y\right>.
\eqno(2.65)
$$
The corresponding graded components $R_k$ satisfy the following conditions
$$
\left< R_k x, y\right> = (-1)^{k+1} \left< x, R_k y\right> ~{\rm
for~any}~x,\,
y\in\V .
\eqno(2.66)
$$

We conclude that the normal form of the system (2.47) is a quadruple
$$
\left( \V , <~,~>, \mu, R\right)
\eqno(2.67)
$$
where
$$
\V = {\rm span}\, (e_1, \dots, e_n)
$$
is $n$-dimensional space with a bilinear symmetric form
$$
<e_\alpha,e_\beta>=\eta_{\alpha\beta}
$$
and an antisymmetric operator
$$
\mu ={\rm diag}\, (\mu_1, \dots, \mu_n)
$$
and a $\mu$-nilpotent $\mu$-skew-symmetric operator $R:\V\to\V$.

Let us describe now the ambiguity in the choice of the normal form data. 
We say that
$$
G:\V\to\V
\eqno(2.68)
$$
is a $\mu$-{\it parabolic orthogonal operator} if 
$$
G=1+\Delta
\eqno(2.69)
$$
where $\Delta$ is a $\mu$-nilpotent operator and $G$ satisfies the
following orthogonality condition w.r.t. the bilinear form (2.65)
$$
\{ G\,x, G\,y\} =\{ x,y\}.
\eqno(2.70)
$$
Representing $\Delta$ as a sum of the graded components
$$
\Delta \simeq \Delta_* =\oplus_{k\geq 1}\Delta_k
\eqno(2.71)
$$
one rewrites the orthogonality condition in the form
$$
\left( 1-\Delta_1^T + \Delta_2^T -\dots \right) \, \eta\, \left(
1+\Delta_1 + \Delta_2 +\dots \right) =\eta.
\eqno(2.72)
$$
\smallskip
{\bf Exercise 2.10.} Prove that the monodromy operator (2.57) is orthogonal
w.r.t.
the bilinear form (2.65).
\medskip
Clearly all $\mu$-parabolic orthogonal operators $G$ form a group
denoted  ${\cal G}\left(\mu, <~,~>\right)$. The space of all
$\mu$-nilpotent $\mu$-skew-symmetric operators $R$ coincides with the Lie
algebra of the nilpotent group. The group  ${\cal G}\left(\mu,
<~,~>\right)$ acts on this space by conjugations
$$
R\mapsto G^{-1}R\, G.
\eqno(2.73)
$$
In the grading components one has
$$
\eqalign{
R_1 &\mapsto R_1\cr
R_2 &\mapsto R_2 +[R_1, \Delta_1]\cr
R_3 &\mapsto R_3 + [R_2, \Delta_1] - \Delta_1 R_1 \Delta_1 + [R_1,
\Delta_2] + \Delta_1^2 R_1\cr}
\eqno(2.74)
$$
etc. Two $\mu$-nilpotent $\mu$-skew-symmetric operators related
by a conjugation (2.73) will be called {\it equivalent}.
\smallskip
{\bf Lemma 2.7.} {\it The set of all normal forms at $z=0$ of a given
system
(2.47) is in one-to-one correspondence with the orbit of one normal form
w.r.t. the action (2.73) of the  group  ${\cal G}\left(\mu, <~,~>\right)$.}

{\bf Proof.} Let us consider two normal forms of the same system (2.47)
$$
\eqalign{
\partial_z \xi &= \left( {\mu\over z} + R_1 + R_2 z+\dots \right) \xi\cr
\partial_z \tilde\xi &= \left( {\mu\over z} + \tilde R_1 + \tilde R_2
z+\dots \right) \tilde\xi.\cr}
$$
They must be related by a gauge transformation
$$
\tilde\xi = G(z)\xi
$$
where
$$
G(z) = 1 + \Delta_1 z + \Delta_2 z^2 + \dots
$$
satisfies (2.45c). Explicitly we obtain a system of relations identical 
to (2.53) $$
\tilde R_n =R_n + n\, \Delta_n + [\Delta_n, \mu] + \sum_{k=1}^{n-1}
\left( \Delta_{n-k}R_k -\tilde R_k\Delta_{n-k}\right), ~n=1, \, 2, \dots
$$
From this system we recursively prove that
$$
\left( \Delta_n\right)_\beta^\alpha =0 ~{\rm unless}~ \mu_\alpha-\mu_\beta
=n.
$$
So
$$
G:= G(1) = 1+\Delta_1 +\Delta_2 +\dots
$$
is a $\mu$-parabolic orthogonal operator $G:\V\to\V$. From (2.45c) we derive
the orthogonality condition (2.70) in the form (2.72).

Let us derive the relation (2.73) for the operators
$$
R=R_1 +R_2 +\dots, ~~\tilde R = \tilde R_1 + \tilde R_2 +\dots.
$$
Since $\xi=z^\mu z^R$ is a solution to (2.49) (see Lemma 4),
$$
\tilde \xi =G(z) z^\mu z^R
$$
must be a solution to the system with tilde. So we must have
$$
G(z) z^\mu z^R = z^\mu z^{\tilde R} C
$$
for an invertible matrix $C$. Let us rewrite the last equation  in the form
$$
z^{-\mu} G(z) z^\mu z^R = z^{\tilde R} C.
$$
Using the identities 
$$
z^{-\mu} \Delta_k z^\mu = z^{-k}\Delta_k
$$
we finally obtain
$$
\left( 1+\Delta_1 +\Delta_2 +\dots \right) z^R =G\,z^R =z^{\tilde R}.
$$
Expanding
$$
\eqalign{
G\,z^R &= G\left( 1+R\log z +{R^2\over 2!}\log^2 z+\dots\right) \cr
= z^{\tilde R} C &= \left( 1+\tilde R\log z +{\tilde R^2\over 2!}\log^2
z+\dots\right) C\cr}
$$
and equating the coefficients in front of various powers of $\log z$ we
obtain
$$
\eqalign{
C &= G\cr
G\,R &= \tilde R\, G.\cr}
$$
Lemma is proved.
\medskip
{\bf Definition 2.4.} A quadruple
$$
\left( \V, <~,~>, \mu, [R]\right)
\eqno(2.75)
$$
where $\V$ is $n$-dimensional linear space with a bilinear symmetric form
$<~,~>$, an antisymmetric diagonalizable operator $\mu$ and a class of
equivalence $[R]$ of normal forms (2.49) of the system (2.47) will be called
{\it monodromy data at $z=0$} of this system.
\smallskip
{\bf Lemma 2.8.} {\it Two systems of the form
$$
\partial_z \xi^{(i)} = \left( {\mu\over z} +\sum_{k\geq 1} U_k^{(i)}
z^{k-1}\right) \xi^{(i)}, ~~i=1, \, 2
\eqno(2.76)$$
satisfying
$$
{U_k^{(i)}}^T =(-1)^{k+1} \eta\, U_k^{(i)} \eta^{-1}
$$
are equivalent w.r.t. a gauge transform of the form (2.45) {\rm iff} they
have the same monodromy data (2.75).}

{\bf Proof.} Let the gauge transformations
$$
\tilde \xi^{(i)} =G^{(i)}(z)\xi^{(i)}, ~~i=1,\,2
$$
reduce the systems (2.76) to the normal forms
$$
\partial_z \tilde \xi^{(i)} =\left({\mu\over z} +\sum_{k\geq 1} R_k^{(i)}
z^{k-1} \right) \, \tilde\xi^{(i)}, ~~i=1,\,2.
$$
If
$$
R^{(1)}=G\,R^{(2)}G^{-1}
$$
with
$$
G=1+\Delta_1 +\dots \in {\cal G} \left( \mu, <~,~>\right)
$$
then the gauge transformation
$$
\tilde\xi^{(2)} =(1+z\, \Delta_1 +z^2\Delta_2 +\dots )\, \tilde\xi^{(1)}
$$
establishes a gauge equivalence of the systems (2.76) for $i=1$ and $i=2$.
Thus the systems (2.76) are gauge equivalent with
$$
\xi^{(2)} ={G^{(2)}}^{-1}(z) \left( 1+z\, \Delta_1 +z^2\Delta_2
+\dots\right) G^{(1)}(z) \, \xi^{(1)}.
$$
Conversely, from Lemma 2.7 it follows that gauge equivalent systems have the
same monodromy data. Lemma is proved.
\medskip
We will now return to Frobenius manifolds. The last component (2.41b) of the
system determining horizontal sections of the connection $\tilde\nabla$
is a linear system of ODEs with rational coefficients of the form (2.47).
The coefficients $\U_\beta^\alpha(t)
=E^\epsilon(t)c_{\epsilon\beta}^\alpha(t)$ depend on the point $t$ of the
Frobenius manifold as on the parameter. The solutions $\xi$ take values in
the space $\V = T_tM$. We may identify the tangent planes at different
points $t$ using the Levi-Civita connection $\nabla$ on $M$. Actually,
the space $\V$ is equipped with an additional structure, namely, a
marked vector $e\in
\V$. This is an eigenvector of the linear operator $\mu$ with the
eigenvalue $-d/2$.
 
We now show
\smallskip
{\bf Isomonodromicity Theorem (first part).} {\it The monodromy data at
$z=0$ of the system (2.41b) do not depend on $t\in M$.}

{\bf Proof.} The matrix $\mu$ and the bilinear form $<~,~>$ are
$t$-independent by construction. We will now construct a $t$-independent
representative $R$ of the normal form (2.49) of the equation (2.41b)
$$
\partial_z\xi = \left( \U+{\mu\over z}\right)\,\xi
$$
$$
\U_\beta^\alpha(t)=E^\epsilon(t) c_{\epsilon\beta}^\alpha(t).
$$
Let us choose a basis $h_1(t;z)$, \dots, $h_n(t;z)$ of solutions of the
system
$$
\partial_\alpha\partial_\beta h_\gamma(t;z) =z\,
c_{\alpha\beta}^\epsilon(t) \partial_\epsilon h_\gamma(t;z)
$$
of the form (2.33), (2.34). Multiplying, if necessary, the series
$$
h_\alpha(t;z) =t_\alpha +\sum_{p\geq 1} h_{\alpha,p}(t) z^p
$$
by a matrix-valued series
$$
M(z)\equiv \left( M_\beta^\alpha (z)\right) = 1+\sum z^k M_k 
$$
$$
h_\alpha(t;z)\mapsto \sum_\epsilon h_\epsilon(t;z) M_\alpha^\epsilon(z)
$$
with $t$-independent coefficients $M_1$, $M_2$, \dots we obtain the
identity
$$
\left< \nabla\,h_\alpha(t;-z),\nabla\,h_\beta(t;z)\right>
=\eta_{\alpha\beta}
$$
(see Exercise 2.8). Let us construct a $n\times n$-matrix series $G(t;z)
=\left( G_\beta^\alpha(t;z)\right)$
$$
\eqalign{G_\beta^\alpha (t;z) &= \eta^{\alpha\epsilon} \partial_\epsilon
h_\beta(t;z)\cr
G(t;z)&= 1+\sum_{k\geq 1} G_k(t) z^k.\cr}
$$
The identity (2.45c) reads
$$
G^T(t;-z)\,\eta\,G(t;z) \equiv \eta.
$$

We do now a gauge transform in the system (2.41)
$$
\xi = G(t;z)\, \tilde\xi.
$$
Since $G$ is a matrix solution of the equations
$$
\tilde\nabla_\alpha G=0, ~~\alpha=1, \dots, n,
$$
we obtain
$$
\partial_\alpha \tilde\xi =0, ~~\alpha=1, \dots, n.
$$
The system (2.41b) after the gauge transform will read in a form (2.47) 
with the matrices $U_1$, $U_2$, \dots satisfying the symmetry/antisymmetry
conditions (2.48)
$$
\partial_z \tilde\xi =\left({\mu\over z} + U_1 + z\, U_2 + \dots \right)
\, \tilde\xi.
$$
The full system of the last two equations remains to be compatible after the 
gauge transform. This implies $t$-independency of the coefficients $U_1$,
$U_2$, \dots. Hence the normal form of the system (2.41b) does not depend on
$t$. Theorem is proved.
\medskip
{\bf Definition 2.5.} The monodromy data (2.75) of the system (41b) are 
called {\it monodromy data of the Frobenius manifold at $z=0$.}
\medskip
Explicitly,
$$
\eqalignno{ 
{R_1}_\beta^\alpha &= \U_\beta^\alpha ~{\rm
for}~\mu_\alpha-\mu_\beta =1
& {(2.77a)} \cr
{R_2}_\beta^\alpha &= \sum_{\mu_\alpha-\mu_\gamma \neq 1}
{\U_\gamma^\alpha \U_\beta^\gamma \over
\mu_\alpha-\mu_\gamma-1} ~{\rm for} ~\mu_\alpha-\mu_\beta =2
& {(2.77b)} \cr}
$$
etc.
\smallskip
{\bf Exercise 2.11.} Prove, using the formula (2.36), that the normalized
coefficients $A_{\alpha\beta}$, $B_\alpha$, $C$ in the quasihomogeneity
equation (WDVV3) for the free energy $F(t)$ have the form
$$
A_{\alpha\beta} ={R_1}_\alpha^\epsilon \eta_{\epsilon\beta},
\eqno(2.78)
$$
particularly,
$$
r_\alpha ={R_1}_1^\alpha,
\eqno(2.79)
$$
$$
B_\alpha ={R_2}_\alpha^\epsilon \eta_{\epsilon 1}
\eqno(2.80)
$$
$$
C=-{1\over 2} {R_3}_1^\epsilon \eta_{\epsilon 1}.
\eqno(2.81)
$$
\medskip
{\bf Remark 2.3.} We defined our monodromy data as {\it formal
invariants},
i.e., all the gauge transforms (2.45) were defined by formal power series
$G(z) = 1+ G_1 z + \dots$. It is wellknown, however, that at a regular
singularity of the system (2.41b) formal invariants coincide with analytic
ones [CL]. In other words, all the normalizing transformations are
convergent series for sufficiently small $|z|$.

We obtain
\smallskip
{\bf Theorem 2.1.} {\it For any Frobenius manifold with the monodromy data
$\left( \V, <~~>, \mu, [R]\right)$ there exists a fundamental matrix
$$
\eqalign{
\Xi^0(t;z)& = H(t;z) z^\mu z^R\cr
H(t;z) &= 1+H_1(t)z +H_2(t)z^2 =\dots\cr}
\eqno(2.82)$$
for the system (2.41) defining horizontal sections of $\tilde\nabla$. The
power series converges for sufficiently small $|z|$. Here $R$ is a
representative of the class $[R]$. 

A change of the representative 
$$
R\mapsto \tilde R = G^{-1}R\,G
$$
where
$$
G=1+\Delta_1 +\Delta_2 +\dots
$$
is a $\mu$-parabolic orthogonal operator, transforms
$$
\Xi^0 \mapsto \tilde{\Xi}^0=\Xi^0 G.
$$
The power series transforms as follows
$$
\tilde{H}(t;z) \left( 1+z\, \Delta_1 +z^2\Delta_2 +\dots \right) = H(t;z).
$$

A choice of the fundamental matrix $\Xi^0(t;z)$ determines a system of
deformed flat coordinates $\tilde t_1(t;z)$, \dots, $\tilde t_n(t;z)$
such that the gradients $\nabla\,\tilde t_1$, \dots, $\nabla\, \tilde t_n$
are the columns of the matrix $\Xi^0$
$$
\left( \tilde t_1(t;z), \dots, \tilde t_n(t;z)\right) =
\left( h_1(t;z), \dots, h_n(t;z)\right) z^\mu z^R.
\eqno(2.83)
$$ 
The functions $\tilde t_\alpha
(t;z)$ are determined uniquely up to $t$-independent shifts and
$\mu$-parabolic transformations. }
\medskip
We will call the functions $\tilde t_\alpha(t;z)$ {\it the normalized
deformed flat coordinates} on the Frobenius manifold. Recall that the
columns of the matrix $H(t;z)$ are gradients of the functions
$h_1(t;z)$, \dots, $h_n(t;z)$. The Taylor expansions of these functions
for small $|z|$ have the form
$$
h_\alpha (t;z) = t_\alpha + \sum_{p=1}^\infty h_{\alpha, p}(t) z^p
$$
with the coefficients $h_{\alpha, p}(t)$ satisfying the system of
recursion relations (2.34). But now the coefficients are determined uniquely
within the ambiguity given by the action of the group of
$\mu$-parabolic orthogonal transformations and up to a $t$-independent
shift.
\smallskip
{\bf Exercise 2.12.} For a $\mu$-nilpotent operator
$$
R=R_1 + R_2 +\dots
$$
define the operators $R_{k,l}$ putting
$$
\eqalign{R_{0,0} &=1\cr
R_{k,0} &=0, ~~k>0\cr
R_{k,l} &=\sum_{i_1+\dots i_l=k} R_{i_1} \dots R_{i_l}\cr}
\eqno(2.84)
$$
Prove that the normalized deformed flat coordinates have the following
expansion near $z=0$:
$$
\tilde t_\alpha(t;z) =\sum_{k, \,l\geq 0} \sum_{p=0}^k \sum_\epsilon
h_{\epsilon,p}(t) \left( R_{k-p,l}\right)_\alpha^\epsilon z^{k+\mu_\alpha}
{\log^l z\over l!}.
\eqno(2.85)
$$
\medskip
{\bf Exercise 2.13.} Derive the following quasihomogeneity conditions for
the
gradients of the functions $h_{\alpha,k}(t)$
$$
{\cal L}_E \nabla\, h_{\alpha,k} = \left( k+{1\over 2} (d-2) +
\mu_\alpha\right)\,
\nabla\,h_{\alpha,k} +\sum_{\epsilon,\, p} \nabla\,h_{\epsilon,k-p}
\left( R_{p}\right)_\alpha^\epsilon.
\eqno(2.86)
$$
\medskip
These quasihomogeneity conditions together with the relations (2.34) can
serve as the recursive definition of the functions $h_{\alpha,k}$ starting
from $h_{\alpha,0}=t_\alpha$.

We conclude this Lecture with descrption of the monodromy data of
Frobenius manifold with good analytic properties (partuclarly, of quantum
cohomologies).
\smallskip
{\bf Proposition 2.2.} {\it For the Frobenius manifold of the form (2.9) all
the matrices $R_2$, $R_3$, \dots vanish and
$$
\left( R_1\right)_\beta^\alpha = \sum_\epsilon r_\epsilon
c_{\epsilon\beta}^\alpha.
\eqno(2.87)
$$
Here the numbers $r_\epsilon$ enter into the Euler vector field (2.4),
$c_{\alpha\beta}^\gamma$ are the structure constants of the cubic
part of $F(t)$ (2.9).}

{\bf Proof.} Due to Isomonodromicity Theorem it is sufficient to compute
the monodromy data of the operator (2.41b) at the point $t=t_0$ of 
``classical limit''. We have
$$
\U_\beta^\alpha (t_0) = \sum_\epsilon r_\epsilon c_{\epsilon\beta}^\alpha.
$$
The algebra $A_{t_0}$ is graded by the degree
$$
\deg e_\alpha =q_\alpha.
$$
On the other hand, the vector
$$
\sum_\epsilon r_\epsilon e_\epsilon
$$
has degree one. Thus the operator (2.38) of multiplication in $A_{t_0}$
by this
vector increases degrees by one. Hence
$$
\U_\beta^\alpha \neq 0 ~{\rm only~if}~q_\alpha-q_\beta=1,
$$
and the system
$$
\partial_z \xi =\left(\U(t_0) +{\mu\over z}\right)\,\xi
$$
is already in the normal form (2.49) with $R=R_1 =\U(t_0)$. Proposition is
proved.
\medskip
{\bf Corollary.} {\it In the quantum cohomology of a manifold $X$
the monodromy data at $z=0$ is the operator $R=R_1$ of multiplication
by the first Chern class $c_1(X)$  acting in the classical cohomologies
$H^*(X)$.} 
\smallskip
{\bf Example 2.3.} Let us explain who is the deformed connection (2.28) in
the
case of quantum cohomology of a ``sufficiently good'' (for 
example, smooth projective) $2d$-dimensional manifold $X$ (the
assumptions and notations are as in Lecture 1 above). Let $\phi\in H^*(X;
\C )$ be an arbitrary element. We will construct function $\tilde t_\phi
(t; z)$, $t=(t', t'')$ as in (2.10), such that for any $\phi$ it satisfies

$$
\tilde\nabla d\tilde t_\phi =0.
\eqno(2.88)
$$ 
Taking $\phi=\phi_1$, \dots,
$\phi=\phi_n$ for a basis in $H^*(X,\C
)$ we will obtain a system of flat coordinates of the deformed connection.

Denote $Q$ the grading operator (1.4). We introduce a line bundle ${\cal L}$
on the moduli space $X_{[\beta],l}$. The fibre of this bundle in the point
$\left(\beta, p_1, \dots, p_l\right)\in X_{[\beta],l}$ is the cotangent
line
to the Riemann sphere at the first marked point $p_1$. Let
$$
\sigma_1 := c_1\left( {\cal L}\right) \in H^2\left( X_{[\beta],l}\right).
$$
Put
$$
\tilde t_\phi (t;z) =z^{-{d\over 2}}\sum_{[\beta],l} \left< { z^Q
z^{c_1(X)} \phi \over 1-z \sigma_1} \otimes 1 \otimes
e^{t''}\right>_{[\beta],l}
e^{\int_{S^2}\beta^*(t')}.
\eqno(2.89)
$$
Here we define the symbols
$$
\left< {a_1\over 1-z \sigma_1} \otimes a_2 \otimes \dots \otimes a_k
\right>_{[\beta],l}
$$
as the formal series in $z$ using the expansion
$$
{1\over 1-z \sigma_1} = 1 + z \sigma_1 + z^2 \sigma_1^2 + \dots ,
$$
and
$$
\left< a_1 \sigma_1^m \otimes a_2 \otimes \dots \otimes
a_k\right>_{[\beta], l} =
\cases{0, ~k\neq l\cr
\int_{X_{[\beta],l}} \sigma_1^m \wedge p_1^*(a_1)\wedge p_2^*(a_2) \wedge
\dots
\wedge p_l^*(a_l), ~k=l.\cr}
\eqno(2.90)
$$
These are particular {\it gravitational descendants} arising in the
description of coupling of the topological sigma-model (=quantum
cohomology) to topological gravity [Wi2, Dij1, Dij2, Du3].
\smallskip
{\bf Theorem 2.2.} {\it The function $\tilde t_\phi (t;z)$ for any
$\phi\in
H^*(X)$ satisfies the equation (2.88).
}

Proof. Let us choose some basis $\phi_1$, \dots, $\phi_n$ in $H^*(X)$. The
formula (2.89) for the functions $\left( \tilde t_{\phi_1} (t;z), \dots,
\tilde t_{\phi_n}(t;z)\right)$ can be rewritten in the form
$$
\left( \tilde t_{\phi_1} (t;z), \dots,
\tilde t_{\phi_n}(t;z)\right) =
\left( h_1 (t;z), \dots, h_n(t;z) \right) z^\mu z^R
$$
where the formal series $h_\alpha (t;z)$ have the form
$$
h_\alpha (t;z) = t_\alpha + \sum_{p=1}^\infty h_{\alpha,p}(t) z^p,
~\alpha=1, \dots, n
$$
$$
h_{\alpha, p} (t) = \sum_{[\beta], l} \left< \sigma_1^p \phi_\alpha
\otimes 1 \otimes e^{t''} \right>_{[\beta],l} e^{\int_{S^2} \beta^* (t')},
\eqno(2.91)
$$
$$
\mu=\diag (\mu_1, \dots, \mu_n), ~~\mu_\alpha = q_\alpha -{d\over 2}, ~
\phi_\alpha \in H^{2 q_\alpha}(X),
$$
$R$ is the matrix of the operator of multiplication by the first Chern
class $c_1(X)$ (cf. (2.83) above). To demonstrate (2.88) it is sufficient to
prove that the coefficients $H_{\alpha,p}(t)$ satisfy the recursion
relations
$$
\partial_\lambda \partial_\mu h_{\alpha,p}(t) = c_{\lambda\mu}^\nu 
\partial_\nu h_{\alpha,p-1}(t), ~p\geq 1,
\eqno(2.92)
$$
(see (2.34)) and
$$
\L \nabla h_{\alpha,p} = \left( p +{d-2 \over 2} + \mu_\alpha \right)
h_{\alpha, p} 
+ \nabla h_{\epsilon, p-1} (R)_\alpha^\epsilon.
\eqno(2.93)
$$
The last equation is the particular case of (2.86) for the case of quantum
cohomology where $R=R_1$ is the matrix of multiplication by the first
Chern class and $R_2 =R_3 =\dots =0$ (see above Corollary from 
Proposition 2.2). The first relation (2.92) follows from
the genus 0 topological recursion relations of Dijkgraaf and Witten [DW]
(see the derivation in [Du3, Du7]). The second follows from the recursion
relations of Hori [Ho]. Theorem is proved.
\medskip
In general we do not know anything about analytic properties of the series
(2.89) for big $|z|$.
But if the quantum cohomology algebra is
semisimple (conjecturally, this is the case for Fano varieties $X$, see
below Lecture 3) then the asymptotic behaviour of the
series (2.89) for big $|z|$ is under control. The Stokes parameters of this
asymptotic behaviour will give us in Lecture 4 additional parameters
of the Frobenius manifold to determine it uniquely.
\medskip
{\bf Example 2.4.} We will now compute the flat coordinates of the
deformed
connection $\tilde \nabla$ for the Frobenius manifolds arising in the
singularity theory. We consider here only the case of simple singularities
$f(x)$
(see the definition in [AGV, Ar2]). (The formulation of K.Saito's
theory of primitive forms in the setting of Frobenius manifolds for more
general singularities can be found in [Sab], [Man3], [Ta].)
Simple singularities are labelled
by the simply-laced Dynkin
diagrams $A_n$, $D_n$, $E_6$, $E_7$, $E_8$. Denote $f_t(x)$ the
correspondent versal deformation. The variable $x$ is one-dimensional for
$A_n$
or $x=(x_1, x_2)$ for other simple singularities. The parameters $t=(t^1,
\dots, t_n)$ for $A_n$, $D_n$, $E_n$. Explicitly
$$
\eqalign{
A_n : ~f_t(x) &= x^{n+1} + a_n x^{n-1} + \dots +a_1
\cr
D_n :~ f_t (x) & = x_1^{n-1} + x_1 x_2^2 +a_{n-1} x_1^{n-2} + \dots + a_1
+
b\, x_2
\cr
E_6 : ~f_t(x) &= x_1^4 + x_2^3 + a_6 x_1^2 x_2 + a_5 x_1 x_2 + a_4 x_1^2 +
a_3 x_2 + a_2 x_1 + a_1
\cr
E_7 : ~f_t(x) &= x_1^3 x_2+ x_2^3  + a_1 + a_2 x_2 + a_3 x_1 x_2 + a_4 x_1
x_2^2 + a_5 x_1 + a_6 x_1^2 + a_7 x_2^2
\cr
E_8 : ~f_t(x) &= x_1^5 + x_2^3 + a_8 x_1^3 x_2 + a_7 x_1^2 x_2 + a_6 x_1^3
+ a_5 x_1 x_2 + a_4 x_1^2 + a_3 x_2 + a_2 x_1 + a_1.
\cr}
$$
The coefficients $a_i$ are certain polynomials of the flat coordinates
$t^1$, \dots, $t^n$. The dependence on the flat coordinates satisfies the
following two remarkable identities [EYY2]
$$
\phi_\alpha(x;t) \phi_\beta(x;t) =c_{\alpha\beta}^\gamma(t)
\phi_\gamma(x;t) +K_{\alpha\beta}^a (x;t) {\partial f_t(x)\over \partial
x^a}
\eqno(2.94)$$
$$
\partial_\alpha\phi_\beta(x;t) = {\partial K_{\alpha\beta}^a(x;t)\over
\partial x^a}.
\eqno(2.95)$$
Here
$$
\phi_\alpha(x;t) ={\partial f_t(x)\over \partial t^\alpha}, ~\alpha = 1,
\dots, n,
$$
$K_{\alpha\beta}^a(x;t)$ are certain polynomials, the index $a$ takes
only one value for $A_n$ or two values for $D_n$, $E_n$. The coefficients
$c_{\alpha\beta}^\gamma(t)$ coincide with the structure constants of the
correspondent Frobenius manifold.
\smallskip
{\bf Theorem 2.3.} {\it The oscillatory integrals
$$
\tilde t_C(t;z) = z^{N-2\over 2} \int_C e^{z\, f_t(x)} d^N x
\eqno(2.96)
$$
($N=1$ for $A_n$ and $N=2$ for $D_n$, $E_n$) are flat coordinates of the
deformed connection $\tilde \nabla$. Here $C$ is any $N$-dimensional cycle
in $\C^N$ that goes to infinity along the directions where ${\rm Re}\,
z\,f_t(x)\to -\infty$.
}

Proof. We are to prove that the functions
$$
\xi_\alpha =\partial_\alpha \tilde t(t;z) =  z^{N\over 2} \int_C
\phi_\alpha (x;t) e^{z\,
f_t(x)} d^N x
$$
satisfy the system (2.88). Using (2.94), (2.95) we obtain
$$
\eqalign{
\partial_\alpha \xi_\beta &=
 z^{N\over 2} \int_C {\partial K_{\alpha\beta}^a(x;t)\over \partial x^a}
e^{z\, f_t(x)} d^N x
\cr
&+
z^{N+2\over 2} \int_C \left[ c_{\alpha\beta}^\gamma
(t) \phi_\gamma (x;t) + K_{\alpha\beta}^a(x;t) {\partial f_t(x)\over
\partial x^a} \right] e^{z\, f_t(x)} d^N x
\cr
&= z\, c_{\alpha\beta}^\gamma (t) \xi_\gamma + z^{N\over 2} \int_C
{\partial\over \partial x^a} \left[ K_{\alpha\beta}^a (x;t) e^{z\,
f_t(x)}\right]\, d^N x
\cr
&= z\, c_{\alpha\beta}^\gamma (t) \xi_\gamma
\cr}
$$
(we used Stokes formula
$$
\int_C {\partial\over \partial x^a} v^a d^N x =0
$$
for any vector field $v^a$ vanishing on the boundary of $C$ at infinity).

To demonstrate the second equation (2.41b) it is sufficient to prove that
$$
z\partial_z \tilde t_C = \L \tilde t_C + {d-2\over 2}\tilde t_C.
$$
Here the Euler vector field and $d$ have the form
$$
E=\sum_{\alpha=1}^n {d_\alpha\over h} t^\alpha \partial_\alpha, ~d_\alpha
= m_\alpha +1
$$
$$
d=1-{2\over h}
$$
where $h$ is the Coxeter number and $m_\alpha$ are the exponents of the
correspondent Weyl group $W(A_n)$, $W(D_n)$, $W(E_n)$ (see Lecture 5
below). One can assign certain degrees $r_1$, $r_2$ to the variables
$x_1$, $x_2$ ($r_1$ only for $A_n$) in such a way that the whole
deformation $f_t(x)$ be a quasihomogeneous function of $t^1$, \dots,
$t^n$, $x_1$, $x_2$ of the degree 1. Explicitly,
$$
\eqalign{
A_n: ~ r_1 &= {1\over n+1}
\cr
D_n: ~ r_1 &= {1\over n-1}, ~r_2 ={n-2\over 2n-2}
\cr
E_6: ~r_1 &= {1\over 4},~ r_2 ={1\over 3}
\cr
E_7 : ~ r_1 &= {2\over 9}, ~ r_2 ={1\over 3}
\cr
E_8 : ~ r_1 &= {1\over 5}, ~ r_2 ={1\over 3}
\cr}
$$
Note that the coefficients $a_\alpha$ of the versal deformations 
are quasihomogeneous polynomials of $t^1$, \dots, $t^n$ of the degrees
$d_\alpha /h = \deg t^\alpha$. The quasihomogeneity can be recasted in the
form of the following Euler identity
$$
\sum_a r_a x^a {\partial f_t(x)\over \partial x^a} + \sum_\alpha
{d_\alpha\over h} t^\alpha {\partial f_t(x) \over \partial t^\alpha}
=f_t(x).
$$
Using this identity we obtain
$$
\eqalign{
\partial_z \tilde t_C &= {N-2\over 2 z} \tilde t_C 
+ z^{N-2\over 2} \int_C f_t(x) e^{z\, f_t(x)} d^N x 
\cr
&=  {N-2\over 2 z} \tilde t_C
+ z^{N-2\over 2} \int_C
\left[
\sum_a r_a x^a {\partial f_t(x)\over \partial x^a} + \sum_\alpha
{d_\alpha\over h} t^\alpha {\partial f_t(x) \over \partial t^\alpha}
\right]\, e^{z\, f_t(x)} d^N x
\cr
&=  {N-2\over 2 z} \tilde t_C + {1\over z} \L \tilde t_C 
+z^{N-4\over 2} \int_C \sum_a {\partial\over \partial x^a}
\left( r_a x^a e^{z\, f_t(x)} \right) \, d^N x
-{r_1+r_2\over z} \tilde t_C
\cr
&=  {N-2\over 2 z} \tilde t_C -{1\over z} \left( r_1+r_2-{N-2\over
2}\right) \tilde t_C.
\cr}
$$
It remains to check that in all these cases
$$
r_1+r_2 -{N-2\over 2} = {2-d\over 2} ={h+2\over 2 h}.
$$
Theorem is proved.
\medskip 
Let us show that, for some cycles $C_1$, \dots, $C_n$ the oscillatory
integrals $\tilde t_{C_1}(t;z)$, \dots,  $\tilde t_{C_n}(t;z)$ give
independent flat coordinates of the deformed connection $\tilde\nabla$.
First we will rewrite, following [AGV], the integral (2.96) as a Laplace-type
transform of an appropriate function $p_C(\lambda; t)$:
$$
\tilde t_C(t; z) =z^{N-2\over 2} \int_0^\infty e^{z\lambda} p_C(\lambda;
t)
\, d\lambda
\eqno(2.97)
$$
where the integration is to be taken along any ray in the half-plane ${\rm
Re}\, z\lambda <0$. Put
$$
p_C(\lambda; t) := \oint_{C(\lambda)} {d^Nx\over df_t(x)}.
\eqno(2.98)
$$
Here the Gelfand - Leray form $ {d^Nx / df_t(x)}$ is defined by the
equation
$$
d^Nx = df_t(x) \wedge  {d^Nx\over df_t(x)},
$$
the $(N-1)$-cycle $C(\lambda)$ is the intersection of $C$ with the level
surface
$$
V_\lambda(t) =\{ x ~|~ f_t(x)=\lambda \}.
\eqno(2.99)
$$

Fixing a noncritical for $f_t(x)$ value $\lambda_0$ we obtain the {\it
period map} 
$$
t\mapsto \left[  {d^Nx\over df_t(x)}\right] \in H^{N-1}\left(
V_{\lambda_0}(t)\right)
\eqno(2.100)
$$
where the square brackets denote the cohomology class of the form. The map
is defined for those $t$ when $\lambda_0$ is not a critical value of
$f_t(x)$. In the
coordinates the map reads
$$
t\mapsto \left( p_{\sigma_1}(\lambda_0;t) , \dots,
p_{\sigma_n}(\lambda_0;t)\right)
\eqno(2.101)
$$
for a basis 
$$
\sigma_1, \dots, \sigma_n \in H_{N-1} \left( V_{\lambda_0}(t); {\bf
Z}\right).
$$
The period map is known to be a local diffeomorphism [Loo]. Now, choosing
a basis of $N$-cycles $C_1$, \dots, $C_n$ such that the $(N-1)$-cycles
$C_1(\lambda)$, \dots, $C_n(\lambda)$ are linearly independent we obtain
independent flat coordinates   
 $\tilde t_{C_1}(t;z)$, \dots,  $\tilde t_{C_n}(t;z)$.

Finally, we note that nodegeneracy of the period map was not essential to
prove independency of the oscillatory integrals. One could use, instead,
the analysis of the asymptotic behaviour of the integrals when
the $\lambda$ goes to one of the critical values of $f_t(x)$ and $C$ is 
the correspondent vanishing cycle. We hope, however, that this digression 
into the singularity theory would help to understand the constructions of 
Lecture 5. 
\medskip
{\bf Remark 2.4.} The two main classes of examples of Frobenius manifolds
look so different. There are, however, some unexpected relationships
between these two classes of two-dimensional topological field theories.
That is, main playing characters of a two-dimensional topological field
theory constructed from quantum cohomology turn out to coincide with those
coming from singularity theory. This phenomenon was first discovered
in quantum cohomologies of Calabi - Yau varieties [COGP]. It was called
{\it
mirror conjecture} (now partially proved [Gi2 - Gi4]). In the last lecture
we will
present our version of mirror construction for semisimple Frobenius
manifolds. Particularly, we will express
the deformed flat coordinates of $\tilde\nabla$ on any semisimple
Frobenius manifold satisfying certain nondegeneracy condition by
oscillatory integrals, and we also obtain an analogue of the residue
formulae (1.20), (1.21). Some general approaches to mirror conjecture were 
recently proposed in [Wi4, Gi3, LLY, BK].

\vfill\eject
\def\nbh{{neighborhood }}
\def\bn{{{\cal B}_n}}
\def\rh{{Riemann - Hilbert boundary value problem }}
\def\rr{{\rm right}}
\def\ll{{\rm left}}
\def\fo{{\rm formal}}
\def\arg{{\rm arg}\,}
\def\cc{{{\cal C}(\mu,R)}}
\def\c0{{{\cal C}_0 (\mu, R)}}

\def\L{{{\cal L}_E}}
\def\nab{{\tilde\nabla}}

\def\res{\mathop{\rm res}}
\def\deg{\mathop{\rm deg}}
\def\wdvv{WDVV equations of asociativity }
\def\V{{\cal V}}
\def\U{{\cal U}}
\def\C{{\bf C}}
\def\diag{{{\rm diag}\,}}
\centerline{Lecture 3}
\medskip
\centerline{ \bf Semisimplicity and canonical coordinates.}
\medskip
In this Lecture we introduce the class of semisimple Frobenius manifolds 
and obtain the main geometrical tool of dealing with them: the canonical 
coordinates [Du3].

We recall that a commutative associative algebra $A$ is 
called 
{\it semisimple} if it contains no nilpotents, i.e., nonzero vectors
$a\in A$ such that
$$
a^m=0
$$
for some positive integer $m$. If all vectors $a\in A$ are nilpotents
then the algebra $A$ is called {\it nilpotent}.
\smallskip
{\bf Lemma 3.1.} {\it Any finite-dimensional Frobenius algebra over $\C$
is isomorphic to orthogonal direct sum
$$
A\simeq A_0 \oplus \C \oplus \dots \oplus \C
\eqno(3.1)
$$
where $A_0$ is a nilpotent algebra.}

{\bf Proof.} Let $\lambda_0,\, \lambda_1, \dots, \lambda_k\in A^*$ be
the roots of the commutative algebra $A$, i.e., such linear functions
$\lambda_i :A\to \C$ that for any $a\in A$ the eigenvalues of the operator
of multiplication by $a$ are $\lambda_0(a), \dots, \lambda_k(a)$. If
$\lambda_0\equiv 0$ is not a root then the algebra is semisimple. Let
$$
A_j :=\cap_{a}{\rm Ker}\, \left( a\cdot -\lambda_j(a)\right)^n, 
~j=0,\, 1, \dots, k
$$
be the correspondent root subspaces. Particularly, $A_0$ consists of all
nilpotent vectors of $A$.
It is easy to see
that the algebra is
decomposed into the orthogonal direct sum of the root subspaces
$$
A=\oplus_{j=0}^k A_j
$$
$$
A_j\cdot A_i =0 ~~{\rm for}~i\neq j
\eqno(3.2)
$$
and, thus,
$$
\lambda_i\left(A_j\right) =0, ~~{\rm for}~ i\neq j.
$$
For $i>0$ let $0\neq v_i\in A_i$ be an eigenvector, i.e., such a vector
that
$$
v\, v_i =\lambda_i(v)\, v_i 
$$
for any $v\in A$. If $\lambda_i (v_i) =0$ then
$$
v_i^2 =0.
$$
So the vector $v_i$ is nilpotent. Thus it must belong to $A_0$. But 
$A_0 \cap A_i =0$. Hence $\lambda_i (v_i)\neq 0$. Put
$$
\pi_i = {v_i\over \lambda_i(v_i)}.
$$
We obtain that 
$$
\pi_i^2 =\pi_i.
$$

Let us prove that each $A_i$ is one-dimensional for $i\neq 0$. If
$w_i\neq0$ is another eigenvector in $A_i$ then
$$
w_i\, \pi_i =\lambda_i(w_i)\pi_i =\lambda_i(\pi_i)w_i =w_i.
$$
So $w_i$ is proportional to $\pi_i$. Similarly, one can see that in $A_i$
there are no vectors adjoint to $\pi_i$. Indeed, if $\pi_i'$ is an adjoint
vector of the height one, i.e.,
$$
\left( \pi_i -\lambda_i(\pi_i)\right) \, \pi_i' =\pi_i
$$
then
$$
\eqalign{\pi_i\,\pi_i'&=\pi_i +\pi_i'\cr
 &=\lambda_i(\pi_i')\, \pi_i.\cr}
$$
This contradicts to linear independence of $\pi_i$ and $\pi_i'$. So, any
of $A_i$ for $i=1, \dots, k$ is a one-dimensional subalgebra in $A$
generated by the vector $\pi_i$ such that
$$
\pi_i^2=\pi_i.
\eqno(3.3)$$
From (3.2) it follows that
$$
\pi_i \pi_j =0 ~~{\rm for }~i\neq j.
\eqno(3.4)
$$
So $\pi_1, \dots, \pi_k$ are the basic idempotents of the semisimple part
$A_1\oplus\dots \oplus A_k$ of $A$. Lemma is proved.
\medskip
{\bf Corollary 3.1.} {\it Any $n$-dimensional semisimple Frobenius algebra
over $\C$ is isomorphic to an orthogonal direct sum of $n$ copies of $\C$.
The basic idempotents $\pi_1$ \dots, $\pi_n$ of the algebra are determined
uniquely up to reordering.}
\medskip
{\bf Definition 3.1.} A Frobenius manifold $M$ is called {\it semisimple}
if
the algebras $T_tM$ are semisimple for generic $t\in M$.
\medskip
Semisimplicity of an algebra is an open property. So, if at some point
$t=t_0\in M$ the algebra $T_tM$ is semisimple then it remains semisimple
in some neighborghood of $t_0$.
\smallskip
{\bf Exercise 3.1.} Prove that the function
$$
F(t_1, t_2, t_3, t_4) ={1\over 2} t_1^2 t_4 + t_1 t_2 t_3 + f(t_2)
$$
and the Euler vector field
$$
E=t_1 \partial_1 -t_3 \partial_3 - 2 t_4 \partial_4
$$
give a solution of WDVV (with $d=3$) for an arbitrary function $f(t_2)$.
\medskip
So, nonsemisimple Frobenius manifolds may depend on functional parameters.
This fact is wellknown to experts in mirror symmetry: \wdvv 
provide
no information about GW invariants of Calabi - Yau three-folds. For a
Calabi - Yau three-fold the Frobenius structure is identicaly nilpotent.

In the opposite case  semisimple Frobenius manifolds depend on finite 
number of parameters. Part of these parameters were described in Lecture
2:
they are monodromy data at $z=0$ of the connection $\tilde\nabla$. In
Lecture 4 for semi-simple Frobenius manifolds we will define also
monodromy
data at $z=\infty$. We will show that the full list of the monodromy data
is a complete local invariant of a semisimple Frobenius manifold. We will
also
describe the global structure of these manifolds (i.e., the analytic
continuation of the local structure) in terms of the monodromy data.

Conjecturally, semisimplicity holds true for quantum cohomology of Fano
varieties (see below the example for ${\bf CP}^2$). This conjecture is
partially supported by the results of [TX].
 \smallskip
{\bf Theorem 3.1.} {\it Let $u_1(t)$, \dots, $u_n(t)$ be the eigenvalues
of
the operator of multiplication by the Euler vector field
$$
\det \left( \U (t) -\lambda \cdot 1\right) =(-1)^n \prod_{i=1}^n
\left(\lambda-u_i(t)\right),
\eqno(3.5)
$$
$$
\U_\beta^\alpha(t) =E^\epsilon(t) c_{\epsilon\beta}^\alpha(t).
\eqno(3.6)
$$
Near a semisimple point $t_0\in M$ they can serve as local coordinates. In
these coordinates
$$
{\partial\over\partial u_i} \cdot {\partial\over\partial u_j}
=\delta_{ij} {\partial\over\partial u_i}
\eqno(3.7)
$$
$$
e=\sum_{i=1}^n {\partial\over\partial u_i}
\eqno(3.8)
$$
$$
E=\sum_{i=1}^n u_i {\partial\over\partial u_i}
\eqno(3.9)
$$
$$
<~,~> = \sum_{i=1}^n \eta_{ii}(u) du_i^2 ~~{\rm where}~\eta_{ii}(u)
={\partial t_1\over \partial u_i}, ~~t_1 := \eta_{i\epsilon} t^\epsilon.
\eqno(3.10)
$$
}

{\bf Main Lemma.} {\it Let $M$ be a complex-analytic manifold with a
structure of Frobenius algebras on the tangent planes $T_tM$ depending
analytically on $t$ and satisfying FM1 and FM2 (the quasihomogeneity FM3
not included). Then local coordinates $u_1$, \dots, $u_n$ exist near a
semisimple point $t_0\in M$ such that
$$
{\partial\over\partial u_i}\cdot {\partial\over\partial u_j}
=\delta_{ij}{\partial\over\partial u_i}.
\eqno(3.11)
$$
}

{\bf Proof.} Near a semisimple point $t_0$ one can choose a frame 
of basic idempotents $\pi_1$, \dots, $\pi_n$
$$
\pi_i \cdot \pi_j = \delta_{ij}\pi_i
$$
depending analytically on the point. It is sufficient to show that the Lie
brackets $[\pi_i,\pi_j]$ of these vector fields vanish. Let us use the
deformed flat connection
$$
\tilde\nabla_u v=\nabla_uv+z\, u\cdot v
$$
(no $\tilde\nabla_{d/dz}$ component since we do not assume quasihomogeneity).
Let us introduce the coefficients $\Gamma_{ij}^k$ and $f_{ij}^k$ from the
expansions
$$
\eqalign{
\nabla_{\pi_j} \pi_i &=\sum_k \Gamma_{ij}^k \pi_k\cr
[\pi_i,\pi_j] &=\sum_k f_{ij}^k \pi_k.\cr}
$$
Computing linear in $z$ terms of the curvature
$$
\nab_{\pi_i}\nab_{\pi_j} -\nab_{\pi_j}\nab_{\pi_i} -\nab_{[\pi_i,\pi_j]}=0
$$
we obtain the equation
$$
\Gamma_{kj}^l \delta_i^l+\Gamma_{ki}^l\delta_{kj}-\Gamma_{ki}^l\delta_j^l
-\Gamma_{kj}^l \delta_{ki} = f_{ij}^l \delta_k^l
$$
valid for arbitrary four indices $i$, $j$, $k$, $l$ (no summation w.r.t.
repeated indices in these formulas!). For $l=k$ we obtain
$$
f_{ij}^k=0.
$$
Lemma is proved.
\smallskip
{\bf Proof of Theorem.} Let us prove that
$$
{\cal L}_E \left( {\partial\over\partial u_i}\right) =-
{\partial\over\partial u_i}
\eqno(3.12)
$$
where $u_1$, \dots, $u_n$ are the local coordinates constructed in Main
Lemma. We use the equation (2.23) of the axiom FM3. This reads
$$
\L (a\cdot b) -\L(a)\cdot b -a\cdot \L b =a\cdot b
\eqno(3.13)
$$
for any vector fields $a$ and $b$. Applying this to the basic idempotents
$\pi_i =\partial/\partial u_i$ for
$a=\pi_i$, $b=\pi_j$, $i\neq j$ we obtain
$$
\L(\pi_i) \cdot \pi_j +\pi_i \cdot \L(\pi_j) =0.
$$
Hence
$$
\L (\pi_i)\cdot \pi_j =0 ~{\rm for }~i\neq j.
$$
So $\L (\pi_i) =\lambda_i\, \pi_i $ with some factor $\lambda_i$. Applying
now (3.13) to the case $a=b=\pi_i$ we obtain
$$
\lambda_i \pi_i -2 \lambda_i \pi_i =\pi_i.
$$
So $\lambda_i=-1$. This proves (3.12). Writing the vector-field $E$ in the
coordinates $(u_1, \dots, u_n)$
$$
E=\sum_{i=1}^n E^i(u) {\partial\over\partial u_i}
$$
we obtain from (3.12)
$$
{\partial E^i\over\partial u_j} =0, ~~i\neq j, ~~{\partial
E^i\over\partial
u_i}=1.
$$
Doing, if necessary, a shift of the coordinates $(u_1, \dots, u_n)$ we
arrive at the formula (3.9). The eigenvectors of the operator of
multiplication by this vector-field are $\partial /\partial u_1$, \dots,
$\partial /\partial u_n$. The correspondent eigenvalues are $u_1$, \dots,
$u_n$. To complete the proof of Theorem we observe that the basic
idempotents of a Frobenius algebra are pairwise orthogonal
$$
<\pi_i, \pi_j> =<e, \pi_i\cdot\pi_j> =0 ~{\rm for}~i\neq j.
$$
Consider the 1-form $dt_1$. By definition for any vector $v$
$$
\partial_v t_1 =dt_1 (v) =<e,v>
$$
where $e$ is the unity of the algebra. Hence
$$
<\pi_i, \pi_i> =<e,\pi_i\cdot\pi_i> =\left<e, {\partial\over\partial
u_i}\right> ={\partial t_1 \over \partial u_i}.
$$
Theorem is proved.
\medskip
{\bf Corollary 3.2.} {\it All the points $t\in M$ where the eigenvalues of
$\left( E(t)\cdot \right)$ are pairwise distinct are semisimple.}
\medskip
{\bf Definition 3.2.} The coordinates $(u_1, \dots, u_n)$ constructed in
Theorem 1 are called {\it canonical coordinates} of the Frobenius
manifold.
\medskip
The canonical coordinates near any point are defined uniquely up to
permutations. We will use Latin indices for canonical coordinates and we
put
$$
\partial_i :={\partial\over\partial u_i}.
$$
We will also show explicitly all the sums w.r.t. Latin indices not
distinguishing between upper and lower indices. Recall that Greek indices
are used for flat coordinates and
$$
\partial_\alpha ={\partial\over\partial t^\alpha}.
$$
The rules of tensor algebra (raising and lowering indices using
$\eta^{\alpha\beta}$ and $\eta_{\alpha\beta}$, the Einstein summation rule
etc.) will be applied only to Greek indices.
\smallskip
We make now an algebraic digression about semisimple Frobenius algebras
over $\C$. Let $(A, <~,~>)$ be such an algebra with a basis $e_1=e$,
$e_2$, \dots, $e_n$ and the multiplication table
$$
e_\alpha\cdot e_\beta =c_{\alpha\beta}^\gamma e_\gamma.
$$
Let $\pi_1$, \dots, $\pi_n$ be the idempotents of $A$. Introduce the basis
of {\it normalized} idempotents
$$
f_i ={\pi_i\over \sqrt{<\pi_i,\pi_i>}}, ~~i=1, \dots, n
$$
choosing arbitrary    signs of the square roots. Let us introduce the
matrix
$\Psi=\left(\psi_{i\alpha}\right)$ putting
$$
e_\alpha =\sum_{i=1}^n \psi_{i\alpha} f_i, ~~\alpha=1, \dots, n.
$$
\smallskip
{\bf Exercise 3.2.} Prove the following formulae
$$
\eqalignno{
\Psi^T\Psi &=\eta
& {(3.14)}
\cr
\psi_{i1} &=\sqrt{<\pi_i,\pi_i>}
&{(3.15)}
\cr
f_i &=\sum_{\alpha,\,\beta=1}^n
\psi_{i1}\psi_{i\beta}\eta^{\beta\alpha}e_\alpha
& {(3.16)}
\cr
c_{\alpha\beta\gamma} &=\sum_{i=1}^n
{\psi_{i\alpha}\psi_{i\beta}\psi_{i\gamma}\over \psi_{i1}}.
& {(3.17)}
\cr}
$$
\medskip
On a semisimple Frobenius manifold the matrix $\Psi$ depends on the point.
The above formula give
$$
\eqalignno{
<~,~> &= \sum_{i=1}^n \psi_{i1}^2(u) du_i^2
& {(3.18)}
\cr
\partial_\alpha &= \sum_{i=1}^n {\psi_{i\alpha}(u)\over \psi_{i1}(u)}
\partial_i
& {(3.19)}
\cr
\partial_i &=\sum_{\alpha, \, \epsilon}\eta^{\alpha\epsilon}
\psi_{i\epsilon}(u)\psi_{i1}(u) \partial_\alpha
& {(3.20)}
\cr}
$$
or, equivalently,
$$
dt^\alpha =\sum_{i=1}^n \psi_i^\alpha (u) \psi_{i1}(u) du_i, ~~
{\rm where} ~ \psi_i^\alpha := \eta^{\alpha\epsilon}\psi_{i\epsilon}.
\eqno(3.21)
$$

We will now rewrite the connection $\nab$ in the frame of normalized
idempotents
$$
f_i ={\partial_i\over \sqrt{<\partial_i,\partial_i>}}.
\eqno(3.22)
$$
We recall that the horizontal sections $\xi$ satisfy the compatible system
$$
\eqalignno{
\partial_\alpha \xi &= z\, C_\alpha \xi,
~~\left(C_\alpha\right)_\gamma^\beta := c_{\alpha\gamma}^\beta
& {(3.23)}
\cr
\partial_z\xi &= \left( \U +{\mu\over z}\right)\, \xi, ~~ \U_\gamma^\beta
:=E^\epsilon c_{\epsilon\gamma}^\beta, ~\mu=\diag(\mu_1, \dots,
\mu_n).
& {(3.24)}
\cr}
$$
The operator $\U$ of multiplication by the Euler vector field becomes
diagonal in the basis $f_1, \dots, f_n$
$$
\Psi\, \U\, \Psi^{-1} =:U =\diag (u_1, \dots, u_n).
\eqno(3.25)
$$
We introduce also the matrix 
$$
V:= \Psi\, \mu\, \Psi^{-1}
\eqno(3.26)
$$
of the operator $\mu$ in the same basis. From antisymmetry
$$
<\mu\, a, b>+<a, \mu\,b>=0
$$
it follows antisymmetry of the matrix $V$
$$
V^T+V=0.
\eqno(3.27)
$$
\smallskip
{\bf Lemma 3.2.} {\it After the gauge transformation
$$
y=\Psi\xi
\eqno(3.28)
$$
the system (2.41) reads
$$
\eqalignno{
\partial_iy &= \left( z\, E_i +V_i\right)\, y, ~~i=1, \dots,
n
& {(3.29)}
\cr
\partial_zy &= \left( U+{V\over z}\right)\,y.
& {(3.30)}
\cr}
$$
Here $E_i$ are the matrix unities
$$
\left( E_i\right)_{ab} =\delta_{ia}\delta_{ib},
\eqno(3.31)
$$
$V_i$ are skew-symmetric matrices uniquely determined by the equations
$$
[U,V_i] =[E_i,V].
\eqno(3.32)
$$
The matrices $V$ and $\Psi$ satisfy the differential equations
$$
\eqalignno{
\partial_i \Psi &= V_i \Psi
& {(3.33)}
\cr
\partial_i V &= [V_i, V].
& {(3.34)}
\cr}
$$
}

Observe that the matrices $V_i$ are defined in those points where the
canonical coordinates are pairwise distinct. Symbolically, (3.32) can be
recasted in the form
$$
V_i = ad_{E_i} {ad_U}^{-1} \left( V\right).
\eqno(3.35)
$$

{\bf Proof.} Using (3.25) one obtains
$$
\partial_i \xi = z\, \Pi_i \xi
$$
where $\Pi_i$ is the operator of multiplication by $\pi_i$. By definition
of $\Psi$
$$
\Psi\, \Pi_i \Psi^{-1} =E_i.
$$
So
$$
\partial_i y =z\, E_i y + \tilde V_i y
$$
where
$$
\tilde V_i := \partial_i \Psi \cdot \Psi^{-1}.
\eqno(3.36)
$$
Using the orthogonality (3.14) we obtain antisymmetry of $\tilde V_i$. From
compatibility
$$
\partial_i\partial_j y =\partial_j\partial_i y
$$
it follows
$$
\left[ E_i, \tilde V_j\right] = \left[ E_j, \tilde V_i\right]
$$
for any $i, \, j$. This implies existence of a symmetric matrix $\Gamma$
such that
$$
\tilde V_i = \left[ E_i, \Gamma \right], ~~i=1, \dots, n.
$$
The off-diagonal entries of $\Gamma$ are determined uniquely. In the points
of $M$ where $u_i\neq u_j$ for any $i\neq j$ we thus obtain a uniquely
determined skew-symmetric matrix $\tilde V$ such that
$$
\left[ U,\tilde V_i\right]=\left[ E_i, \tilde V\right], ~~i=1, \dots, n.
$$
Doing the gauge transformation (3.28) in the system (2.41b) we obtain
$$
\partial_z y= \left( U+{V\over z}\right)\, y.
$$
The compatibility $\partial_i\partial_z = \partial_z\partial_i$ implies
$$
V=\tilde V,
$$
$$
\partial_i V =[V_i,V], ~~i=1, \dots, n.
$$
The definition (3.36) of the matrix $\tilde V_i = V_i$ reads
$$
\partial_i \Psi = V_i \Psi.
$$
Lemma is proved.
\medskip
{\bf Exercise 3.3.} Let us consider $V=\left(V_{ij}(u)\right)$ as a
function
of $u=(u_1,\dots. u_n)$ with the values in the Lie algebra $so(n)$. Prove
that the equations (3.34) can be considered as time-dependent Hamiltonian
systems
$$
{\partial V\over \partial u_i} =\left\{ V, H_i(V;u)\right\} , ~~i=1, \dots,
n
\eqno(3.37)
$$
with the quadratic Hamiltonians
$$
H_i(V;u) ={1\over 2} \sum_{j\neq i} {V_{ij}^2\over u_i-u_j}, ~~i=1, \dots,
n
\eqno(3.38)
$$
w.r.t. the standard linear Poisson bracket on $so(n)$
$$
\left\{ V_{ij}, V_{kl}\right\}
= V_{il}\delta_{jk}-V_{jl}\delta_{ik}+V_{jk}\delta_{il}-V_{ik}\delta_{jl}.
\eqno(3.39)
$$
The canonical coordinates $u_1, \dots, u_n$ play the role of the time
variables of these Hamiltonian systems.
\smallskip
{\bf Exercise 3.4.} Prove that $\{ H_i, H_j\} =0$ for any $i$, $j$. From
this and from commutativity of the flows (3.37) derive that the form
$$
\sum_{i=1}^n H_i(V;u) du_i
$$
is closed for any solution $V(u)$ of the system (3.34). That means that
(locally) there exists a function $\tau (u)$ such that
$$
{\partial \log\tau (u)\over \partial u_i} =H_i(V(u); u), ~i=1, \dots, n.
\eqno(3.40)
$$

This is called {\it tau-function} of the solution of the system $V(u)$. 
In the next Lecture we will show that the system (3.33), (3.34) can be 
solved by
reducing to certain linear Riemann - Hilbert boundary value problem.
The tau-function will coincide with the Fredholm determinant of the
correspondent system of integral equations (see [JM, Mi]).

Importance of the tau-function in topological field theory is clear from
the following 
\smallskip
{\bf Theorem 3.2} [DZ2]. {\it Let $X$ be a smooth projective manifold such
that
the quantum cohomology of $X$ is semisimple. Then the generating function
$F^{(1)}(t)$ of elliptic Gromov - Witten invariants of $X$ is given by the
formula
$$
F^{(1)}(t) =\log {\tau (u)\over J^{1/24}}|_{u=u(t)}
$$
where $\tau(u)$ is the above tau-function and 
$$
J=\det \left( {\partial t^\alpha \over \partial u_i}\right) =
\psi_{11} \dots \psi_{n1}.
$$
}
\medskip
Particularly, from this theorem it follows validity of Conjectures 0.1 
and 0.2 of recent paper of Givental [Gi5].
\bigskip
We prove now the converse to 
Lemma 2 statement.

Let $V(u)$, $\Psi(u)$ be a solution of the system (3.33), (3.34) with a
diagonalizable matrix $V(u)$. We observe first that the product
$$
\Psi^{-1}(u) V(u)\, \Psi(u)
$$
does not depend on $u$. We can therefore find a constant matrix $C$ in
such a way that
$$
\Psi^{-1}(u) V(u)\, \Psi(u)=C\, \mu\, C^{-1}
$$
where 
$$
\mu=\diag (\mu_1, \dots , \mu_n)
$$
is a constant diagonal matrix. Doing a change
$$
\Psi(u)\mapsto \Psi(u)C
$$
we obtain another solution of the linear system (3.33) such that 
$$
\Psi^{-1}(u) V(u)\, \Psi(u)=\diag (\mu_1, \dots, \mu_n).
\eqno(3.41)
$$
After these preliminaries we formulate
\smallskip
{\bf Lemma 3.3.} {\it Let $V(u)=\left(V_{ij}(u)\right)$,
$\Psi(u)=\left(\psi_{i\alpha}(u)\right)$ be a solution of the system 
(3.33), (3.34)
satisfying (3.41). Then the formulae (3.18), (3.21), (3.17) define a 
Frobenius structure on the domain
$$
u_i\neq u_j ~{\rm for}~i\neq j, ~~\psi_{11}(u)\dots \psi_{n1}(u)\neq 0.
\eqno(3.42)
$$
}
{\bf Proof.} From the antisymmetry of the matrices $V_i$ it follows 
that
$$
\partial_i \left( \Psi^T \Psi\right) =0, ~i=1 \dots, n.
$$
Put
$$
\eta=\left(\eta_{\alpha\beta}\right) =\Psi^T\Psi,
~~\left(\eta^{\alpha\beta}\right) = \eta^{-1}.
\eqno(3.43)
$$
Next step is to prove that the 1-forms 
$$
\sum_{i=1}^n \psi_i^\alpha \psi_{i1} du_i, ~~{\rm where}~
\psi_i^\alpha=\eta^{\alpha\epsilon}\psi_{i\epsilon}
$$
are closed. From (3.33) we obtain
$$
\partial_j \psi_{i\alpha} ={V_{ij}\over u_j-u_i}\psi_{j\alpha}~~
{\rm for~any}~i\neq j, ~{\rm any}~\alpha.
$$
From this the identity
$$
\partial_j\left( \psi_i^\alpha\psi_{i1}\right) =
\partial_i\left( \psi_j^\alpha\psi_{j1}\right)
$$
follows. This proves local existence of the functions $t^\alpha$ such that
$$
dt^\alpha = \sum_{i=1}^n \psi_i^\alpha\psi_{i1}du_i.
$$
The differentials $dt^1$, \dots, $dt^n$ are independent on the domain 
(3.42). 
So $t^1, \dots, t^n$ serve as local coordinates on the domain. From the
orthogonality (3.14) we obtain that
$$
\partial_\alpha =\sum_{i=1}^n {\psi_{i\alpha}\over \psi_{i1}} \partial_i.
$$
The last step is to prove the symmetry
$$
\partial_\delta \left( \sum_{i=1}^n
{\psi_{i\alpha}\psi_{i\beta}\psi_{i\gamma}\over\psi_{i1}}\right) 
=
\partial_\gamma \left( \sum_{i=1}^n
{\psi_{i\alpha}\psi_{i\beta}\psi_{i\delta}\over\psi_{i1}}\right).
$$
To prove this we are to use another consequence of (3.33)
$$
\partial_i\psi_{i\alpha}=-\sum_{k\neq i} \partial_k \psi_{i\alpha}
$$
valid for any $i$, any $\alpha$.
We leave this computation as an exercise for the reader. Lemma is proved.
\medskip
{\bf Corollary 3.3.} {\it Classes of local equivalence of semisimple
Frobenius
manifolds such that 1 is an eigenvalue of $\nabla\, E$ of the multiplicity
$k$ depend on 
$$
k-1+{n(n-1)\over 2}
$$
parameters.
}

{\bf Proof.} Take the initial data 
$$
V(u^0) =\left(V_{ij}(u^0)\right)
\eqno(3.44)
$$
of the antisymmetric matrix $V$ in a point $u^0 = (u_1^0, \dots, u_n^0)$
with
$$
u_i^0\neq u_j^0 ~{\rm for}~i\neq j.
$$
Solving the system (3.34) of commuting ODEs we obtain locally uniquely the
matrix-valued function $V(u)$ and, therefore, the matrices $V_i(u)$. The
solution $\Psi(u)$ of the linear system (3.33) such that
$$
\Psi^{-1}(u) V(u)\,\Psi(u)=\diag (\mu_1, \dots, \mu_n)=\mu
$$
is determined uniquely up to multiplication by a matrix
$$
\eqalign{
\Psi(u) &\mapsto \Psi(u)C \cr
C^{-1}\mu\, C &= \mu.\cr}
$$
The matrices $C$ preserving the direction of the eigenvector $e$ of $\mu$
with the eigenvalue $\mu_1=-d/2$ produce equivalences of the Frobenius
manifolds. The assumption about the multiplicity of the eigenvalue 1 of
$\nabla\,E$ means that the eigenvalue $\mu_1=-d/2$ of $\mu$ has also the
multiplicity $k$. The vectors $C\,e$ considered up to rescalings must be
eigenvectors of $\mu$ with the same eigenvalue $\mu_1$. The directions of
these vectors give $k-1$ parameters additional to the initial data (3.44).
Corollary is proved.

\vfill\eject
\def\nbh{{neighborhood }}
\def\bn{{{\cal B}_n}}
\def\rh{{Riemann - Hilbert boundary value problem }}
\def\rr{{\rm right}}
\def\ll{{\rm left}}
\def\fo{{\rm formal}}
\def\arg{{\rm arg}\,}
\def\cc{{{\cal C}(\mu,R)}}
\def\c0{{{\cal C}_0 (\mu, R)}}

\def\L{{{\cal L}_E}}
\def\nab{{\tilde\nabla}}

\def\res{\mathop{\rm res}}
\def\deg{\mathop{\rm deg}}
\def\wdvv{WDVV equations of asociativity }
\def\V{{\cal V}}
\def\U{{\cal U}}
\def\C{{\bf C}}
\def\diag{{{\rm diag}\,}}
\centerline{Lecture 4}
\medskip
\centerline{ \bf Stokes matrices and classification of semisimple
Frobenius manifolds.}
\medskip
In the previous Lecture we parametrized semisimple Frobenius manifolds $M$  
by
initial data of the system (3.33), (3.34) of differential equations in a
point $t\in
M$ such that $u_i(t)\neq u_j(t)$ for $i\neq j$.
Typically, however, one has no ``natural'' point in the Frobenius
manifold to specify the initial data (3.44). E.g., for Frobenius manifolds
with
good analytic properties the ``natural'' point would be $t_0 =(t''=0,
t'=-\infty)$. But in this point the Frobenius algebra $T_{t_0}M$ 
typically is nilpotent (one can keep in mind the example of quantum
cohomology where $t_0$ is the point of classical limit). So in this
point $u_1(t_0)=u_2(t_0)=\dots =u_n(t_0)$. This is a singular point for 
the system (3.33), (3.34).

Instead, we will use [Du3, Du7] the monodromy data of the system (3.30) as
the
parameters. Recall that the system is gauge equivalent to the equations
(2.41) determinining the horizontal sections of the connection $\nab$.
Part
of the monodromy data has already been defined in Lecture 2. Namely, this
part is the monodromy at $z=0$ of the system (3.30) gauge equivalent to
(2.41b).
Recall that for a system (3.30) with the monodromy data at $z=0$ of the
form
$\left( \V, <~,~>, \mu, [R]\right)$ a fundamental matrix solution 
$Y_0(z)$ exists such that
$$
Y_0=\Phi(z) z^\mu z^R
\eqno(4.1)
$$
where
$$
\Phi(z) = \Psi + z\, \Phi_1 + z^2 \Phi_2 + \dots
\eqno(4.2)
$$
is an invertible matrix holomorphic for small $|z|$ satisfying
$$
\Phi^T(-z)\Phi(z) =\eta.
\eqno(4.3)
$$
Let us describe the ambiguity in the choice of the normalized solution
(4.1).

Let $\c0$ be the group of all invertible matrices $C$ such that
$$
z^\mu z^R C\, z^{-R} z^{-\mu}
=C_0 + z \, C_1 + \dots
\eqno(4.4)
$$
is a matrix-valued polynomial in $z$.
\smallskip
{\bf Lemma 4.1.} {\it Two solutions $Y(z)$, $\tilde Y(z)$ of the system
(3.30)
have
the same form (4.1) {\rm iff} they are related by a right multiplication
by
a
matrix $C\in \c0$.}

{\bf Proof.} If
$$
Y=\Phi(z) z^\mu z^R, ~~\tilde Y=\tilde\Phi(z) z^\mu z^R
$$
satisfy (3.30) then
$$
\tilde Y(z)=Y(z)C
$$
for a constant matrix $C$. We have
$$
\Phi^{-1}(z)\tilde\Phi(z) = z^\mu z^R C\,z^{-R}z^{-\mu}.
$$
Hence the r.h.s. must be a polynomial. The converse statement is obvious.
Lemma is proved.
\medskip
{\bf Exercise 4.1.} Show that the matrices in $\c0$ commute with $\exp
2\pi
i\mu$ and that they must have the form
$$
C=C_0 + C_1 + C_2 +\dots
$$
with
$$
\left( C_k\right)_\beta^\alpha \neq 0 ~~{\rm
only~if}~\mu_\alpha-\mu_\beta=k, ~~k=0,\, 1, \dots.
$$
Particularly, the matrix $C_0$ commutes with $\mu$.
\medskip
{\bf Remark 4.1.} In the case of a Frobenius manifold we have an
additional
structure of the monodromy data of (3.30) at $z=0$. Namely, an eigenvector
$e$ of the matrix $V$ with the eigenvalue $\mu_1 = -d/2$ must be marked.
It corresponds to the unity of $M$. We must therefore to impose an
additional constraint on the matrix $C$ in (4.4): the component $C_0$ must
preserve the marked vector. Observe that the marked vector corresponds to
the first column $\psi_{i1}$ of the matrix $\Psi$.
\medskip
The second part is the monodromy data at $z=\infty$ that we are going to
define now.

We first describe the monodromy data at $z=\infty$ of the system
$$
{dy\over dz} =\left( U+{1\over z}V\right)\,y
\eqno(4.5)
$$
with arbitrary $n\times n$ matrices of the form
$$
U=\diag (u_1, \dots, u_n), ~~u_i\neq u_j
\eqno(4.6)
$$
$$
V^T=-V
\eqno(4.7)
$$
being a diagonalizable matrix
$$
\Psi^{-1} V\,\Psi =\mu =\diag (\mu_1,\dots ,\mu_n).
\eqno(4.8)
$$
The point $z=\infty$ is an irregular singularity of the system (3.30). So
in
the problem of normal form of the system (3.30) we are to distinguish
formal
gauge equivalences
$$
\eqalign{
y &\mapsto G\left({1\over z}\right) y\cr
 G\left({1\over z}\right) &= 1+{G_1\over z} +{G_2\over z^2} + \dots
\cr}
$$
and analytic ones, where the series converges for sufficiently large
$|z|$.
\smallskip
{\bf Definition 4.1.} Two systems 
$$
{dy^{(i)}\over dz} =\left( U+{1\over z}V^{(i)}\right)\,y^{(i)}, ~~i=1,\,2
\eqno(4.9)
$$
are called {\it analytically equivalent at $z=\infty$} if there exists
a gauge transform
$$
y^{(2)} =  G\left( z\right) \, y^{(1)}
\eqno(4.10)
$$
with the matrix-valued function $ G\left( z\right)$ analytic at
$z=\infty$ satisfying $G(\infty)=1$ and the orthogonality condition
$$
 G^T\left(- z\right) G\left(z\right)=1.
$$
The {\it monodromy at $z=\infty$ of the system (3.30)} is the class of
analytic
equivalence of this system.
\medskip
Below we will explain how one can parametrize the monodromy at infinity by
Stokes matrices of the system (3.30). But first we will show that the
system
is, to some extent, uniquely determined by the monodromy at $z=0$ and
$z=\infty$.
\smallskip
{\bf Lemma 4.2.} {\it Let (4.9) be two systems analytically equivalent at
$z=\infty$. Then the matrix $G$ establishing the gauge equivalence is a
rational function of $z$ of the form
$$
G= 1 +{G_1\over z} +{G_2\over z^2}+\dots +{G_m\over z^m}.
\eqno(4.11)
$$
}
{\bf Proof.} Let the given gauge transform (4.10) be analytic for $|z|>M$
for
some constant $M$. Choose  a point $z_0$ with $|z_0|>M$ and the
fundamental matrix solutions $Y^{(i)}(z)$ of the systems (4.9) with the
initial data
$$
Y^{(i)}(z_0) =1, ~i=1,\,2.
$$
For any $z$ with $|z|>M$ we must have
$$
G(z) Y^{(1)}(z) = Y^{(2)}(z)C
$$
for some constant nondegenerate matrix $C$. The solutions $Y^{(1,\,2)}(z)$
can be continued analytically along any path in $\C\setminus 0$. The
formula
$$
G(z) =Y^{(2)}(z)C\,{Y^{(1)}(z)}^{-1}
$$
gives analytic continuation of $G(z)$ (recall that $\det Y^{(i)}(z) =\exp
(z-z_0) \sum_i u_i  \neq 0$). We obtain a single-valued analytic
function in ${\bf \bar C}\setminus 0$ such that $G(\infty)=1$. Near the
point of regular singularity $z=0$ the entries of the matrices
$Y^{(1,\,2)}(z)$ grow not faster than some power of $|z|$. Hence also
$G(z)$ has at most power growth at $z=0$. So it must be a rational
function having a pole only at $z=0$. Lemma is proved.
\medskip
{\bf Exercise 4.2.} Prove that the determinant of the matrix (4.11) is 
identically equal to 1.
\medskip
{\bf Remark 4.2.} Gauge transformations with rational $G(z)$ are called
{\it
Schlesinger transformations} [JM]. For the case of Frobenius manifolds
they induce certain symmetries of WDVV, i.e., changes of variables
$$
\eqalign{
t &\mapsto \hat t
\cr
F &\mapsto \hat F
\cr}
$$
mapping solutions to solutions. We give here the explicit form [Du7] of
such
symmetries for the case $m\leq 1$ in (4.11).

Type 1. $G=\,$const, $G\mu =\mu G$, $G$ permutes the two eigenvectors of
$\mu$ with the numbers 1 and $\kappa$. Then
$$
\eqalign{
\hat t_\alpha &= \partial_\alpha \partial_\kappa F(t)
\cr
{\partial^2 \hat F\over \partial \hat t^\alpha \partial \hat t^\beta} &=
{\partial^2  F\over \partial  t^\alpha \partial t^\beta}
\cr
\hat\eta_{\alpha\beta} &= \eta_{\alpha\beta}.
\cr}
\eqno(4.12)
$$

Type 2. 
$$G= 1+{A\over z}
$$
where
$$
A_{ij}={\psi_{i1}\psi_{j1}\over t_1}.
$$
Then
$$
\eqalign{
\hat t^1 &= {1\over 2} {t_\sigma t^\sigma \over t_1}
\cr
\hat t^\alpha & = {t^\alpha\over t_1}, ~\alpha \neq 1,\, n
\cr
\hat t^n &= -{1\over t_1}
\cr
\hat F &= t_1^{-2} \left[ F - {1\over 2} t^1 t_\sigma t^\sigma\right]
\cr
\hat\eta_{\alpha\beta} &= \eta_{\alpha\beta}.
\cr}
\eqno(4.13)
$$
Also one may take superposition of (4.13) with any transformation of the
form (4.12)
\medskip
We classify now the systems of the form (3.30) having the same monodromies 
at $z=0$ and $z=\infty$. We will show that, generically, these systems 
must coincide. There remain, however, some subtleties in the nongeneric
situation. The ambiguity of the reconstruction of the system (3.30)
starting 
from the monodromy data at $z=0$ and $z=\infty$ will be completely 
described in terms of the monodromy at $z=0$.

Let us choose a representative $R$ in the class of equivalence $[R]$ of 
the monodromy data at $z=0$ of the system
(3.30). Let us consider the centralizer of the monodromy matrix
$$
M_0 = \exp 2 \pi i (\mu +R)
\eqno(4.14)
$$
in the group of invertible matrices, i.e., the matrices $C$ commuting 
with $M_0$
$$
C^{-1}M_0 C = M_0.
\eqno(4.15)
$$
For any such a matrix $C$ the product
$$
z^\mu z^R C z^{-R} z^{-\mu} =\sum_k A_k z^k
\eqno(4.16)
$$
is a matrix-valued Laurent polynomial in $z$. Particularly, for the matrix
$C\in \c0$ the r.h.s. of (4.16) contains only nonnegative powers of $z$.
Denote
$\cc$ the quotient group of the centralizer (4.15) over the subgroup
$\c0$.
\smallskip
{\bf Example 4.1.} For a nonresonant $\mu$ the group $\cc$ consists of one 
element.
\smallskip
{\bf Example 4.2.} The group $\cc$ with a resonant $\mu$ and $R=0$ is not 
trivial. It is isomorphic to the subgroup of ``upper triangular'' parabolic
matrices in the centralizer of $\exp 2 \pi i \mu$
$$
C= \dots +C_{-2} +C_{-1} +1
\eqno(4.17)
$$
where
$$
\left( C_k\right) _\beta^\alpha \neq 0 ~~{\rm only 
~if}~\mu_\alpha-\mu_\beta=k, ~k = -1, \, -2, \dots .
\eqno(4.18)
$$
\medskip
Let the two systems of the form (4.9) have the same monodromy data at
$z=0$ 
and $z=\infty$. We will asociate with such a pair a matrix $C\in \cc$ where
$\mu$, $R$ are the monodromy data of the systems (4.9) at $z=0$. Let 
$Y^{(1)}(z)$ be the matrix solution of the system
$$
\partial_z Y^{(1)} =\left( U+{1\over z}V^{(1)}\right)\, Y^{(1)}
\eqno(4.19)
$$
of the form
$$
Y^{(1)}(z) =\Phi(z) z^\mu z^R.
$$
Let $G_0(z)=1 +O(z)$ and $G_\infty(z)=1+O\left(1/z\right)$ be the gauge 
transformations of the system (4.19) to another system of the same form
$$   
\partial_z Y^{(2)} =\left( U+{1\over z}V^{(2)}\right)\, Y^{(2)}.
\eqno(4.20)
$$
The matrix-valued functions $G_0(z)$ and $G_\infty(z)$ are assumed to be 
analytic near $z=0$ and $z=\infty$ resp. Near $z=0$ we obtain a solution
$$
 Y_0^{(2)}(z) =G_0(z) Y^{(1)}(z)
$$
of (4.20). Continuing $Y^{(1)}(z)$ analytically along a ray $\rho$ in the 
neighborhood of infinity we produce another solution of the system (4.20)
$$
 Y_\infty^{(2)}(z) =G_\infty(z) Y^{(1)}(z).
$$
Continuing  $Y_\infty^{(2)}(z)$ back along the same ray $\rho$, we obtain 
two matrix solutions of (4.20) defined in a neighborhood of $z=0$. They
must 
be related by a multiplication by an invertible matrix $C_{12}$
$$
 Y_\infty^{(2)}(z)= Y_0^{(2)}(z) C_{12}.
$$
We rewrite the last equation in the form
$$
G_0^{-1}(z) G_\infty(z)=\Phi(z) z^\mu z^R C_{12} z^{-R} z^{-\mu}\Phi^{-1}(z).
\eqno(4.21)
$$
The r.h.s. must be a meromorphic function near $z=0$. That 
means, particularly, 
that the matrix $C_{12}$ commutes with the monodromy matrix $M_0$.
We arrive at
\smallskip
{\bf Theorem 4.1.} {\it The set of all systems 
$$
\partial_z \tilde Y =\left( U+{1\over z}\tilde V\right)\,\tilde Y
$$
of the form (3.30) having the monodromy data at $z=0$ and $z=\infty$
coinciding
with those of the given system
$$
\partial_z Y =\left( U+{1\over z}V\right)\, Y
$$
is in one-to-one correspondence with the elements of the group $\cc$.
}

{\bf Proof.} The above construction associates with the pair of this systems
an element $C=C_{12}$ of the centralizer of $M_0$. It remains to
show that
the two systems coincide {\it iff} 
$$
z^\mu z^R C\, z^{-R} z^{-\mu}
$$
is a polynomial in $z$. Indeed, if this is the case then the r.h.s. of 
(4.21) is analytic at $z=0$. Hence $G_\infty(z)$ is analytic at $z=0$.
Using 
the normalization $G_\infty (\infty)=1$
we conclude that $G_\infty (z)\equiv 1$. The converse statement is obvious.
Theorem is proved.
\medskip
We proceed now to a ``quantative'' description of the monodromy at infinity
of systems of the form (3.30). We first show that all the systems (3.30)
with given pairwise distinct values of $u_1$, \dots, $u_n$
are
gauge
equivalent at $z=\infty$ w.r.t. {\it formal} gauge transformations. It is 
sufficient to construct a gauge transformation
$$
\tilde Y = G(z) \,Y
\eqno(4.22)
$$
of the system (3.30) to the system with constant coefficients
$$
\partial _z \tilde Y=U\,\tilde Y.
\eqno(4.23)
$$
\smallskip
{\bf Lemma 4.3.} {\it For any system (3.30) there exists a unique formal
series
$$
G(z) = 1 +{A_1\over z} +{A_2\over z^2} +\dots
\eqno(4.24)
$$
satisfying
$$
G^T(-z)G(z) =1
\eqno(4.25)
$$
such that (4.22) transforms (3.33) to the system (4.23) with constant
coefficients.}

{\bf Proof.} For the coefficients of the formal series (4.24) one obtains 
the recursion relations
$$
\eqalign{
[U,A_1] &=V\cr
[U,A_{k+1}] &=A_k V-k\,A_k, ~~k=1,\,2,\dots.\cr}
$$
Representing
$$
A_k=B_k+D_k
$$
with an off-diagonal matrix $B_k$ and a diagonal one $D_k$ we obtain
$$
\eqalign{
B_1 &={ad_U}^{-1}(V)\cr
D_k &={1\over k} \diag (B_kV)\cr
B_{k+1} &={ad_U}^{-1}\left(A_kV-k\,A_k\right)\cr}
$$
where `$\diag$' stands for the diagonal part of the matrix. This proves 
existence and uniqueness of the series $G(z)$.

Let us choose a fundamental matrix $Y(z)$ for the system (3.30) such that
$$
Y^T(-z)Y(z)\equiv 1.
$$
Then
$$
G(z)Y(z)
$$
is a formal solution of the system (4.23). Hence for an appropriate
constant 
invertible matrix $C$
$$
G(z) Y(z) = e^{z\,U}C.
$$
Computing the product
$$
\left( G^{-1}(z)\right)^T G^{-1}(-z) = e^{z\,U} \left( 
C\,C^T\right)^{-1}e^{-z\,U}
$$
we conclude that $C\,C^T=1$ since the l.h.s. is a formal series in inverse
powers of $z$ of the form $1+O(1/z)$. This proves the orthogonality 
relation (4.25). Lemma is proved.
\medskip
The series $G(z)$ typically diverges. However, in certain sectors of the 
complex $z$-plane near $z=\infty$ it serves as the asymptotic development
of an actual solution of the original system.

We recall that a series
$$
a_0 +{a_1\over z} +{a_2\over z^2}+\dots
$$
is an asymptotic expansion of the function $f(z)$ for $|z|\to\infty$ in 
the sector
$$
\alpha<\arg z<\beta
$$
if for any $n$
$$
z^n \left[ f(z) -\sum_{k=0}^n {a_k\over z^k}\right] \to 0
$$
as $|z|\to\infty$ uniformly in the sector
$$
\alpha+\varepsilon <\arg z <\beta -\varepsilon
$$
for any sufficiently small positive $\varepsilon$. This fact will be 
denoted briefly
$$
f(z)\sim a_0 +{a_1\over z} +{a_2\over z^2}+\dots, ~~|z|\to\infty, 
~~\alpha<\arg z<\beta.
$$

Let us denote 
$$
Y_\fo (z) =\left( 1+{A_1\over z} +{A_2\over z^2} +\dots \right)\, e^{z\,U}
\eqno(4.26)
$$
where the coefficients of the formal series are defined in Lemma 4.3. We 
say that a matrix solution $Y(z)$ of the system (3.30) has asymptotic 
development
$$
Y(z)\sim Y_\fo(z), ~~|z|\to\infty,
~~\alpha<\arg z<\beta
$$
if in the same sector
$$
Y(z) \, e^{-z\,U} \sim Y_\fo (z)\, e^{-z\,U} =  1+{A_1\over z} +{A_2\over 
z^2} +\dots .
$$
\smallskip
{\bf Definition 4.2.} A line $\ell$ through the origin in the complex 
$z$-plane is called {\it admissible} for the system (3.30) if
$$
{\rm Re}\, z\,(u_i-u_j)|_{z\in \ell\setminus 0}\neq 0 ~~{\rm for 
~any}~i\neq j.
\eqno(4.27)
$$
\medskip
Let us fix an admissible line $\ell$ and an orientation on it. According 
to the orientation the line splits into the negative and positive parts
$\ell_-$ and $\ell_+$ resp. Let the parts have the equations
$$
\eqalign{
\ell_+ &=\{ z \, | \arg z=\phi\}\cr
\ell_- &= \{ z\, | \arg z= \phi-\pi\} .\cr}
\eqno(4.28)
$$
We construct two sectors
$$
\eqalign{
\Pi_\rr :~~\phi-\pi-\varepsilon < &\arg z<\phi+\varepsilon\cr
\Pi_\ll :~~\phi-\varepsilon < &\arg z<\phi+\pi+\varepsilon\cr}
\eqno(4.29)
$$
for sufficiently small positive $\varepsilon$.
\smallskip
{\bf Theorem 4.2.} {\it There exist unique solutions $Y_{\rr/\ll}(z)$
of (3.30) analytic in the sectors $\Pi_{\rr/\ll}$ resp. having the 
asymptotic development
$$
Y_{\rr/\ll}(z)\sim Y_\fo (z)
\eqno(4.30)
$$
as $|z|\to\infty$ in these sectors.}

Proof see in [BJL1].
\medskip
We are now ready to define Stokes matrices of the system (3.30). In the 
narrow sector
$$
\Pi_+ : ~~\phi-\varepsilon< \arg z<\phi+\varepsilon
\eqno(4.31)
$$
we have two solutions. They must be related by multiplication by a matrix
$$
Y_\ll (z) =Y_\rr (z) S, ~~z\in \Pi_+.
\eqno(4.32)
$$
Similarly, in the opposite narrow sector $\Pi_-$
$$
Y_\ll (z) = Y_\rr (z) S_-, ~~z\in \Pi_-.
\eqno(4.33)
$$
\smallskip
{\bf Definition 4.3.} The matrices $S$, $S_-$ are called {\it Stokes
matrices
of the system (3.30)}.
\smallskip
{\bf Lemma 4.4.} {\it Two systems with equal Stokes matrices w.r.t. the
same
admissible oriented line $\ell$ are analytically equivalent near $z=\infty$.}

{\bf Proof.} Let $Y_{\ll/\rr}^{(1)}(z)$,  $Y_{\ll/\rr}^{(2)}(z)$ be the 
solutions of the correspondent systems with the needed asymptotic 
developments in the sectors $\Pi_{\ll/\rr}$ resp. Let us consider the 
following piecewise analytic matrix-valued function $G(z)$ defined for 
sufficiently large $|z|$ such that 
$$
G(z)=\cases{Y_\rr^{(2)}(z){Y_\rr^{(1)}}^{-1}(z), ~~z\in \Pi_\rr\cr
{Y_\ll^{(2)}(z)Y_\ll^{(1)}}^{-1}(z), ~~z\in \Pi_\ll.\cr}
$$
In the sectors $\Pi_+$, $\Pi_-$ we have
$$
\eqalign{
Y_\ll^{(1,\,2)}(z) &=Y_\rr^{(1,\,2)}(z) S, ~~z\in\Pi_+\cr
Y_\ll^{(1,\,2)}(z) &=Y_\rr^{(1,\,2)}(z) S_-, ~z\in\Pi_-.\cr}
$$
So $G(z)$ is a single-valued analytic function for $|z|>M$ for some big 
constant $M$. In the sectors $\Pi_{\rr/\ll}$
$$
G(z)\sim 1+O\left({1\over z}\right).
$$
Hence $z=\infty$ is a removable singularity for this function, and
$$
G(\infty)=1.
$$
This function $G(z)$ establishes the needed gauge transformation between 
the systems. Lemma is proved.
\medskip
We will now describe the algebraic properties of the Stokes matrices. We 
first describe explicitly all non-admissible lines. Each of them 
consists of two {\it Stokes rays} 
$$
R_{ij}:= \left\{ z~ | z=-ir(\bar u_i-\bar u_j), ~r\geq 0\right\} ,~~i\neq j
\eqno(4.34)
$$
(some of them may coincide). We explain: for $z\in R_{ij}$
$$
\big | e^{z\,u_i}\big | = \big | e^{z\, u_j}\big |;
$$
on the right of $R_{ij}$
$$
\big | e^{z\,u_i}\big | < \big | e^{z\, u_j}\big |
$$
and on the left of $R_{ij}$
$$
\big | e^{z\,u_i}\big | > \big | e^{z\, u_j}\big |.
$$
The ray $R_{ji}$ is the opposite one to $R_{ij}$. An admissible line $\ell$
must contain no Stokes rays. The sectors $\Pi_{\rr/\ll}$ can be extended 
up to the first nearest Stokes ray (see [BJL1]).
\smallskip
{\bf Theorem 4.3.} {\it The Stokes matrices $S=(s_{ij})$, $S_-$ of the 
system (3.30) satisfy the following properties
$$
S_- =S^T.
\eqno(4.35)$$
$$
\eqalignno{
s_{ii} &=1, ~i=1, \dots, n
& {(4.36a)}
\cr
s_{ij} &\neq 0 ~~{\rm only ~if}~R_{ij}\subset \Pi_\ll.
& {(4.36b)}
\cr}
$$}

{\bf Proof.} We know that, for any two matrix solutions $Y_1(z)$, $Y_2(z)$
of the system (3.30), the product
$$
Y_1^T(-z)Y_2(z)
$$
does not depend on $z$. Let us choose for $z\in \Pi_\rr$ 
$Y_2(z) =Y_\rr(z)$, $Y_1(z)=Y_\ll(z)$. Using the asymptotic developments
$$
\eqalign{
Y_\rr(z) &\sim G(z) \, e^{z\,U}\cr
Y_\ll(-z) &\sim G(-z) \, e^{-z\,U}\cr}
$$
valid for $z\in \Pi_+$, with $G(z)$ being defined in Lemma 4.3, and the 
orthogonality condition (4.25) we obtain
$$
Y_\ll^T(-z)Y_\rr(z)\equiv 1, ~~z\in \Pi_+.
$$
Let us continue analytically this formula in the counter-clockwise direction
through the ray $\ell_+$. We obtain after the analytic continuation
$$
\eqalign{
Y_\rr(z) &\mapsto Y_\ll(z) S^{-1}\cr
Y_\ll(-z) &\mapsto Y_\rr (-z) S_-.\cr}
$$
So
$$
S_-^T Y_\rr^T(-z)Y_\ll(z)S^{-1}\equiv 1, ~~z\in \Pi_-.
$$
As above we show that 
$$  
 Y_\rr^T(-z)Y_\ll(z)\equiv 1, ~~z\in \Pi_-.
$$
Hence
$$
S_-^T =S.
$$

Let us now prove (4.36). Comparing the asymptotic developments of the both 
sides of (4.32) for $z\in \Pi_+$ we conclude that
$$
e^{z\,U}S\,e^{-z\,U}\sim 1, ~~|z|\to\infty, ~z\in \Pi_+.
$$
This means that
$$
e^{z(u_i-u_j)} s_{ij}\sim \delta_{ij},  ~~|z|\to\infty, ~z\in \Pi_+.
$$
For the diagonal terms this implies $s_{ii}=1$. For the off-diagonal 
terms we have
$$
\big | e^{z(u_i-u_j)}\big | \to \infty ~~{\rm for} ~~|z|\to\infty, ~z\in 
\Pi_+
$$
if $R_{ij}\subset \Pi_\rr$. So, for those pairs $i\neq j$ for which 
$R_{ij}\subset \Pi_\rr$ we must have $s_{ij}=0$. The opposite ray 
$R_{ji}\subset \Pi_\ll$. And
$$
\big | e^{z(u_j-u_i)}\big | \to 0 ~~{\rm for} ~~|z|\to\infty, ~z\in
\Pi_+.
$$
So $s_{ji}$ need not to be zero. Lemma is proved.
\medskip
We see that the Stokes matrix $S$ contains  $n(n-1)/2$
independent parameters. 

To complete the list of the monodromy data we define the central 
connection matrix
$$
Y_0(z) =Y_\rr(z)C, ~~z\in \Pi_+
\eqno(4.37)
$$
(observe: the branchcut in the definition of $Y_0(z)$ is to be chosen 
along $\ell_-$).

The monodromy $(\mu, R)$ at $z=0$, the monodromy $S$ at $z=\infty$, and 
the central connection matrix $C$ are not independent. First of all, we have
the following {\it cyclic relation}
$$
C^{-1}S^T S^{-1}C=M_0 =\exp 2\pi i (\mu+R).
\eqno(4.38)
$$
This expresses a simple topological fact: on the punctured plane 
$\C\setminus 0$ a loop around infinity is homotopic to a loop around the 
origin. Another property comes from the orthogonality relations 
$$
\eqalign{
S &= C \, e^{-\pi i R} e^{-\pi i \mu} \eta^{-1} C^T\cr
S^T &= C\, e^{\pi i R}e^{\pi i \mu} \eta^{-1} C^T.\cr}
\eqno(4.39)
$$
We leave the proof of these identities as an exercise for the reader.

The matrix $C$ is defined up to transformations of the form
$$
C\mapsto B\,C, ~~B\,S\,B^T=B
\eqno(4.40a)
$$
preserving the relations (4.38), (4.39), and
$$
C\mapsto C\, C_0, ~~C_0 \in \c0
\eqno(4.40b)
$$
corresponding to a change of the solution $Y_0(z)$.
\smallskip
{\bf Exercise 4.3.} Prove that classes of equivalence (4.40) of
central 
connection matrices of systems (3.30) with a given monodromy $(\mu, R)$ at 
the origin and a given monodromy $S$ at infinity are in one-to-one 
correspondence with the group $\cc$.
\medskip
The properties (4.38) and (4.39) typically specify the central connection
matrix
$C$ of the system with given $\mu$, $R$, $S$ essentially uniquely with an 
ambiguity (4.40) that does not affect the Frobenius structure. This
reflects the 
claim of Theorem 4.1 (here ``typically'' means triviality of the group 
$\cc$). Anyhow, the following uniqueness theorem holds.
\smallskip
{\bf Lemma 4.5.} {\it If two systems 
$$
\partial_z Y^{(1,\,2)} =\left( U+{1\over z}V^{(1,\, 2)}\right)\, Y^{(1,\, 2)}
$$
have the same matrices $\mu$, $R$, $S$ (w.r.t. the same admissible 
oriented line $\ell$), $C$ then $V^{(2)}=V^{(1)}$.}

The proof is similar to that of Lemma 4.4. We leave it as an exercise.
\medskip
Let us return to semisimple Frobenius manifolds. Starting from a point 
$t_0\in M$ such that the eigenvalues $u_1(t_0)$, \dots, $u_n(t_0)$ of the
operator $\U(t_0) =\left(E(t_0)\cdot\right)$ are pairwise distinct, 
ordering these eigenvalues, and choosing signs of the square roots of
$<\partial_i, \partial_i>$, and fixing an oriented line $\ell$
on the complex $z$-plane admissible for the points  $u_1(t_0)$, \dots, 
$u_n(t_0)$ we define the Stokes matrix $S=S(t_0)$ and the central connection
matrix $C=C(t_0)$. We will now prove that these matrices do not change
under small variations of $t_0$. Observe that the property of admissibility
of the line $\ell$ is stable under small perturbations of $t_0$.
\smallskip
{\bf Isomonodromicity Theorem (second part).} {\it The Stokes matrix $S$
and the central connection matrix $C$ do not depend on the point of a 
semisimple Frobenius manifold.}

{\bf Proof.} Due to Lemma 4.3 the coefficients $A_1$, $A_2$, \dots
 of the solution $Y_\fo (z;u)$ are analytic functions on $u$. From
the uniqueness of $Y_\fo (z;u)$ it easily follows that
$$
\partial_i Y_\fo (z;u) =(z\, E_i +V_i)\,Y_\fo (z;u), ~i=1, \dots, n.
$$
The same statements are true for the solutions $Y_{\rr/\ll}(z;u)$ and, as 
we already know from Lecture 2, for the solution $Y_0(z;u)$. Using the 
definitions
$$
S=Y_\rr^{-1}(z;u)Y_\ll(z;u), ~~z\in \Pi_+
$$
$$
C=Y_\rr^{-1}(z;u)Y_0(z;u), ~~z\in \Pi_\rr
$$
we obtain
$$
\partial_iS=0, ~~\partial_iC=0.
$$
Theorem is proved.
\medskip
Together with the results of Lecture 2 we conclude that the monodromy data 
$\mu$, $R$, $S$, $C$ do not depend on the point of the Frobenius manifold.

We will now show how to reconstruct the semisimple Frobenius manifiold
starting from the monodromy data.

To reconstruct the operator (3.30) and the solutions $Y_{\rr/\ll}$, $Y_0$ 
for given $u_1$, \dots, $u_n$, $( \mu, R, S, C)$ one is to solve certain 
\rh. Let $D$ be the disk
$$
|z|<\rho
$$
for some $\rho>0$, $P_\rr$ and $P_\ll$ the two components of $\C\setminus 
\ell$ intersected with the external parts of the disk. We are to construct
a piecewise-analytic function 
$$
\Phi(z)=\cases{\Phi_\rr(z) , ~~z\in P_\rr\cr
\Phi_\ll(z), ~~z\in P_\ll\cr
\Phi_0(z), ~~z\in D\cr}
$$
continues in the closures of $P_\rr$, $P_\ll$, $D$ resp. 
such that:

1). on the positive (i.e., that belonging to $\ell_+$) part of the 
common boundary of
$P_\rr$ and $P_\ll$ the boundary values of the functions are related by
$$
\Phi_\ll(z) =\Phi_\rr (z) e^{z\,U}S\,e^{-z\,U}.
\eqno(4.41)
$$

2). on the negative part of the common boundary of $P_\rr$ and $P_\ll$ 
the boundary values of the functions are related by
$$
\Phi_\ll(z) =\Phi_\rr (z) e^{z\,U}S^Te^{-z\,U}.
\eqno(4.42)
$$

3). on the common  boundary of $D$ and $P_\rr$ the boundary values of the 
functions are related by
$$
\Phi_0(z) =\Phi_\rr(z) e^{z\,U} C z^{-R}z^{-\mu}.
\eqno(4.43)
$$

4).  on the common  boundary of $D$ and $P_\ll$ the boundary values of 
the   
functions are related by
$$
\Phi_0(z) =\Phi_\ll(z) e^{z\,U} S^{-1}C z^{-R}z^{-\mu}.
\eqno(4.44)
$$ 

5). for $|z|\to \infty$ within $P_{\rr / \ll}$ 
$$
\Phi_{\rr / \ll}(z)\to 1.
\eqno(4.45)
$$
\smallskip
{\bf Theorem 4.4.} {\it If the \rh  1 - 5 has a unique solution at a point
$u^0=\left(u^0_1, \dots, u^0_n\right)$, $u_i^0\neq u_j^0$ for $i\neq j$, 
then the unique solution $\Phi =\Phi(z; u_1, \dots, u_n)$ exists for $u$
sufficiently close to $u^0$ and it is an analytic function of $u$. It can 
be continued analytically to a meromorphic function on the universal covering
of the space
$$
\C^n\setminus \diag :=\{ (u_1, \dots, u_n) \, | \, u_i\neq u_j ~{\rm 
for}~i\neq j\} .
\eqno(4.46)
$$}

Proof follows from the general theory of Riemann - Hilbert boundary
value problems (see in [Mi, Ma]).
\medskip
Having a solution $\Phi =\left( \Phi_\rr(z;u), \Phi_\ll (z;u), 
\Phi_0(z;u)\right)$ of the \rh we can reconstruct the solutions
$$
\eqalign{
Y_{\rr/\ll}(z;u) &= \Phi_{\rr/\ll}(z;u)e^{z\,U}\cr
Y_0(z;u) &= \Phi_0 (z;u) z^\mu z^R.\cr}
\eqno(4.47)
$$
Let us introduce notations for the coefficients of the expansion
of the matrix $\Phi_0(z;u) =\left( {\Phi_0}_{i\alpha}(z;u)\right)$
near $z=0$
$$
 {\Phi_0}_{i\alpha}(z;u) =\sum_{p=0}^\infty \phi_{i\alpha,\,p}(u) z^p.
\eqno(4.48)
$$
Observe
$$
\phi_{i\alpha, \, 0}(u) =\psi_{i\alpha}(u).
\eqno(4.49)
$$
\smallskip
{\bf Isomonodromicity Theorem (third part).} {\it Let the \rh
(4.41) - (4.45) for given $\mu$, $R$, $S$, $C$ satisfying (4.36), (4.38),
(4.39) 
have a unique solution $\Phi = \Phi(z; u^0)$ at a point $u^0
=(u^0_1, \dots , u^0_n)$, $u_i^0\neq u_j^0$ for $i\neq j$ such that
$$
\prod_{i=1}^n \phi_{i1, \, 0}(u^0)\neq 0.
\eqno(4.50)
$$
Then the formulae
$$
\eqalignno{
\eta_{\alpha\beta} &= \sum_{i=1}^n \phi_{i\alpha, \, 0}(u) 
\phi_{i\beta,\,0}(u)
& {(4.51)}
\cr
e &= \sum_{i=1}^n \partial_i
& {(4.52)}
\cr
E &= \sum_{i=1}^n u_i\partial_i
& {(4.53)}
\cr
t_\alpha &= \sum_{i=1}^n \phi_{i1,\, 1}(u) \phi_{i\alpha, \, 0}(u)
& {(4.54)}
\cr
c_{\alpha\beta\gamma} &= \sum_{i=1}^n
{\psi_{i\alpha}\psi_{i\beta}\psi_{i\gamma}\over \psi_{i1}}
& {(4.55)}
\cr
F &= {1\over 2} \sum_{i=1}^n \left[\eta^{\alpha\beta}
\phi_{i\alpha,\, 1} \phi_{i\beta,\, 0} \phi_{i 1,\,1}^2
-\phi_{i1,\, 
2}\phi_{i1,\, 1} -\phi_{i1,\, 0}\phi_{i1,\, 3}\right]
& {(4.56)}
\cr}
$$
define a semisimple Frobenius structure on a small neighborhood of
$u^0$.}

{\bf Proof.} Let us define the matrix-valued functions $Y_{\rr/\ll}(z;u)$,
$Y_0(z;u)$ by the formulae (4.47) and prove that they satisfy the linear 
system (3.29), (3.30) with
$$
\eqalignno{
V(u) &= [U,A_1(u)]
& {(4.57)}
\cr
V_i(u) &= [E_i, A_1(u)]
&{(4.58)}
\cr}
$$
where the matrix $A_1(u)$ is defined from the asymptotic development 
$$
A_1(u) := \lim_{|z|\to\infty,\, z\in \Pi_+} z\, \left( \Phi_\rr(z; u) 
-1\right).
\eqno(4.59)
$$
Let us consider the piecewise-analytic function
$$
Y(z;u) =\cases{
Y_\rr(z;u), ~~z\in \Pi_\rr\cr
Y_\ll(z;u), ~~z\in \Pi_\ll\cr
Y_0(z;u), ~~~~z\in D.\cr}
$$
We prove first that the matrix $Y(z;u)$ is invertible for any $z$, $u$. 
Indeed, $\det Y(z;u) \, e^{-z \sum u_i}$ is a piecewise-analytic
function of $z$ having no jumps on the intersections of the domains 
$\Pi_\rr$, $\Pi_\ll$, $D$ and going to 1 when $|z|\to\infty$. Thus
$$
\det Y(z;u)\equiv  e^{z(u_1+\dots +u_n)}.
$$
We introduce now piecewise-analytic functions
$$
G_i(z;u) := \partial_i Y(z;u) \cdot Y ^{-1}(z;u).
$$
From construction of $S$, $C$ it follows that $G_i(z;u)$ has no jumps on 
the intersections of the domains $\Pi_\rr$, $\Pi_\ll$, $D$. So it is an 
analytic matrix-valued function on $\C\setminus 0$. At $|z|\to\infty$
it has the asymptotic development
$$
\eqalign{
G_i(z;u) &=\partial_i
\left[ \left( 1+{A_1\over z} +\dots \right) \, e^{z\, U}\right] \, 
e^{-z\, U} \left( 1-{A_1\over z}+\dots \right)\cr
 &\sim z\, E_i +V_i +O\left({1\over z}\right).\cr}
$$
At $z=0$ the function $G_i(z;u)$ goes to a finite limit
$$
\eqalign{
G_i(z;u) &=\partial_i\left[ \left( \Psi(u) +O(z)\right)z^\mu z^R\right]
z^{-R} z^{-\mu} \left[ \Psi^{-1}(u) + O(z)\right]\cr
 &=\partial_i \Psi(u) \cdot \Psi^{-1}(u) + O(z)\cr}
$$
due to constancy of $\mu$, $R$. Hence
$$
G_i(z;u) = z\, E_i +V_i
$$
and 
$$
\partial_i Y = \left( z\, E_i +V_i\right)\, Y, ~~i=1, \dots, n.
$$
Particularly, 
$$
\partial_i\Psi = V_i \Psi.
$$

Similarly, considering the piecewise-analytic function
$$
G_z : = \partial_z Y(z; u) \cdot Y^{-1}(z;u)
$$
we obtain that
$$
G_z = U+{V\over z}
$$
where the matrix $V=V(u)$ is defined in (4.57).

To prove the orthogonality conditions
$$
\eqalign{
\Phi_{\rr/\ll}^T(-z;u) \Phi_{\rr/\ll}(z;u) &\equiv 1\cr
\Phi_0^T(-z;u) \Phi_0(z;u) &\equiv \eta\cr}
$$
we will consider the piecewise-analytic matrix-valued function
$$
G(z) := \cases{
Y_\rr(z;u) Y_\ll^T(-z;u), ~~z\in \Pi_\rr\cr
Y_\ll(z;u)Y_\rr^T(-z;u), ~~z\in \Pi_\ll.\cr}
$$
For $z\in \Pi_+\cap \Pi_\rr$
$$
\eqalign{
G(z) &= Y_\rr (z;u) Y_\ll^T(-z;u)\cr
 &=Y_\rr(z;u)SY_\rr^T(-z;u).\cr}
$$
For $z\in \Pi_+\cap \Pi_\ll$
$$
\eqalign{
G(z) &= Y_\ll(z;u) Y_\rr^T(-z;u) \cr
 &=Y_\rr(z;u) S Y_\rr^T(-z;u).\cr}
$$
So, $G(z)$ has no jumps on $\ell_+$. Similarly, it has no jumps on 
$\ell_-$. Using (4.39) one obtains that for $z\in \Pi_\rr$ near $z=0$
$$
\eqalign{
G(z) &= Y_0(z;u) e^{\pi i R} e^{\pi i \mu} \eta^{-1} Y_0^T(-z;u) 
=\Phi_0(z;u) \eta^{-1} \Phi_0^T(-z;u)\cr
 &=1 +O(z).\cr}
$$
A similar computation gives the same behaviour of $G(z)$ at $z\to 0$, 
$z\in \Pi_\ll$. So $G(z)\equiv 1$. This proves the orthogonality conditions.

The equations (4.54), (4.56) is the spelling of (2.35), (2.36). Note that
the functions
$t_1(u)$, \dots, $t_n(u)$ are independent coordinates in the points $u$
where the product 
$$
\prod_{i=1}^n \psi_{i\alpha}(u) \neq 0.
$$
Theorem is proved.
\medskip
{\bf Exercise 4.4.} Show that the product (4.50) does not vanish
identically
unless the matrix
$$
e^{z\, U} S\, e^{-z\, U}
$$
is independent on one of the variables $(u_1, \dots, u_n)$.
\medskip
The Isomonodromicity Theorem gives a structure of semisimple Frobenius 
manifold on small domains in the space of isomonodromy deformations of 
the operator
$$
L={d\over dz} -\left( U+{V\over z}\right)
$$
with rational coefficients. The parameters of these Frobenius manifold
are the monodromy data
$$
(\mu, e, R, S, C)
\eqno(4.60)
$$
of the operator satisfying the above properties (4.6), (4.7). Here $e$ is
a marked 
eigenvector of the matrix $V$ with the eigenvalue $\mu_1$ ($\mu_1$
being the marked 
diagonal entry of the matrix $\mu$). The choice of $e$ corresponds to the 
choice
of the first column of the matrix $\Psi$ in the formulae (4.50) - (4.56).
(We need not 
to fix the bilinear form $<~,~>$. It is given by (4.51).) It also 
demonstrates that, locally, any semisimple Frobenius manifold can be 
realized in such a way.
\smallskip
{\bf Exercise 4.5.} We say that the Stokes matrix $S$ is {\it reducible}
if it has the form $S=S'\oplus S''$ w.r.t. some decomposition of the set of
indices $\{1, \dots, n\} =I'\cup I''$ into a union of two non-empty 
non-intersecting subsets. Prove that a reducible matrix $S$ can make a part
of the monodromy data only if $\exp 2\pi i \mu_1$ is the eigenvalue
of both the matrices ${S'}^T {S'}^{-1}$ and ${S''}^T {S''}^{-1}$. Prove 
that the
 Stokes matrix of a reducible Frobenius manifold is reducible
(see Exercise 2.5).
\medskip
We will now describe the structure of analytic continuation of semisimple
Frobenius manifolds. According to Theorem 4.4 and due to the formulae
(4.54), (4.56) 
the functions $t_\alpha$ and $F$ can be continued analytically to 
meromorphic functions on the universal covering of $\C^n\setminus \diag$. 
Since the canonical coordinates are defined up to reordering, the structure
of analytic continuation of the Frobenius manifold with given monodromy 
data (4.60) is described by an action of the fundamental group
$$
\pi_1 \left( \left( \C^n\setminus\diag\right)/S_n, (u_1^0, \dots, 
u_n^0)\right) =\bn
$$
(the braid group) on the monodromy data computed at a given point $u^0$.
The global structure of the Frobenius manifold is described by the 
stationary subgroup $\bn^0\subset\bn$ of the given monodromy data (4.60).

To compute the action of the braid group $\bn$ on the monodromy data, and 
also to describe the dependence of the monodromy data on the admissible 
oriented line $\ell$, we will briefly present here the theory of Stokes 
factors (see [BJL1]).

Let us label all the Stokes rays (4.34) of the system (3.30) in the 
counter-clockwise order starting from the first one in $\Pi_\rr$. We obtain 
the rays
$$
\eqalign{
R^{(1)}, \dots, R^{(m)} ~&{\rm in} ~\Pi_\rr
\cr
R^{(m+1)}, \dots, R^{(2m)} ~&{\rm in}~\Pi_\ll.\cr}
\eqno(4.61)$$
We will use the cyclic labelling $R^{(k\pm 2m)}=R^{(k)}$. Observe that 
the narrow sectors $\Pi_+$ and $\Pi_-$ contain no Stokes rays. For generic
$(u_1, \dots, u_n)$ one has
$$
m={n(n-1)\over 2}
$$
but some coincidences of the Stokes rays may happen in the nongeneric 
situation when there are three $u_i$, $u_j$, $u_k$ on a line or two pairs
$u_i, ~u_j$ and $u_k, ~u_l$ on two parallel lines. Let us consider the sector
of $z$-plane from $R^{(k)} e^{- i \varepsilon\over 2}$ to $R^{(m+k)} 
e^{-i \varepsilon}$. According to Theorem 4.2 there exists a unique
solution
$Y^{(k)}(z)$ of (3.30) such that
$$
Y^{(k)}(z) \sim Y_\fo (z), ~~|z|\to\infty
\eqno(4.62)
$$
within the above sector. This solution can be extended preserving the 
asymptotics into the open sector
$$
\Pi_k: ~{\rm from}~R^{(k-1)} ~{\rm to}~R^{(m+k)}.
\eqno(4.63)
$$
On the intersection of two subsequent sectors one has a constant matrix
$K_j$ defined by
$$
Y^{(j+1)}(z) =Y^{(j)}(z)K_j, ~~z\in \Pi_j\cap \Pi_{j+1}.
\eqno(4.64)
$$
\smallskip
{\bf Lemma 4.6.} 
$$
Y_\rr =Y^{(1)}, ~~Y_\ll = Y^{(m+1)}
\eqno(4.65)
$$
$$
S=K_1 \dots K_m.
\eqno(4.66)
$$

Proof is obvious.
\smallskip
{\bf Definition 4.4.} The matrices $K_j$ are called {\it Stokes factors}
of the matrix $S$.
\medskip
{\bf Exercise 4.6.} Prove that
$$
K_{m+j}K_j^T =1.
\eqno(4.67)
$$
\medskip
How to find the Stokes factors knowing the Stokes matrix $S$ and the 
configuration of pairwise distinct complex numbers $u_1$, \dots, $u_n$?
The clue is in the following property of Stokes factors (see [BJL1]).
\smallskip
{\bf Lemma 4.7.} {\it All the diagonal entries of $K_j$ equal 1. Of the 
off-diagonal entries $\left(K_j\right)_{ab}$ all equal zero but those for 
which the Stokes ray $R_{ba}$ coincides with $R^{(j)}$.}

{\bf Proof.} On $z\in \Pi_j\cap\Pi_{j+1}$ one must have
$$
e^{z\,U}K_je^{-z\,U}\to 1 ~~{\rm as}~|z|\to\infty.
$$
Hence $\left(K_j\right)_{aa}=1$ (as in the proof of Theorem 4.3). On the 
intersection the absolute values
$$
\big | e^{z(u_a-u_b)}\big |
$$
can go to either $+\infty$ or  $0$ for any pair $a\neq b$ but those for 
which $R_{ab}$ or $R_{ba}$ coincides with $R^{(j)}$. Indeed, the whole 
intersection $\Pi_j\cap\Pi_{j+1}$ lies on the right from the oriented line
$$
R^{(m+j)} \cup \left(-R^{(j)}\right).
$$
If
$$
R_{ab} =R^{(m+j)}, ~~R_{ba}=R^{(j)}
$$
then on the right from the oriented line one has
$$
\big | e^{z(u_a-u_b)}\big | \to 0 ~~{\rm as}~|z|\to\infty.
$$
Lemma is proved.
\medskip
{\bf Theorem 4.5.} {\it Any Stokes matrix $S$ with the above properties 
can be uniquely factorized into the product $S=K_1\dots K_m$ of Stokes 
factors of the above form.}

Proof see in [BJL1].
\medskip
From the factorization (4.66) it follows that the Stokes matrix does not
change
if one deforms the admissible line $\ell$ not intersecting any of the 
Stokes rays. We describe now what happens if the oriented admissible line
$\ell =\ell_+ \cup (-\ell_-)$ passes through the Stokes ray $R$ moving
counter-clockwise. 
Instead, 
one may consider a deformation of one of the Stokes rays $R$ passing through
$\ell_+$ moving clockwise.
\smallskip
{\bf Lemma 4.8.} {\it After the above deformation the new solutions 
$Y'_{\rr/\ll}$, the new Stokes matrix $S'$, and the new connection matrix 
$C'$ have the form
$$
\eqalignno{
Y_\rr &=Y'_\rr K_R^T
&{(4.68a)}
\cr
Y'_\ll &= Y_\ll K_R
&{(4.68b)}
\cr
S' &= K_R^T S\, K_R
&{(4.68c)}
\cr
C' &=K_R^T C
&{(4.68d)}
\cr}
$$
(the last formula holds true modulo the ambiguity (4.40).
Here $K_R$ is the Stokes factor corresponding to the Stokes ray $R$.}

Proof follows from Lemma 4.6 and from Exercise 4.6.
\medskip
We are now ready to compute the action of the braid group $\bn$
on the monodromy data describing the analytic continuation of the 
Frobenius manifold. First, the action of $\bn$ on the monodromy at $z=0$ 
is trivial. We now compute the action on the Stokes matrix $S$. Let us 
assume that the canonical coordinates $(u_1,\dots, u_n)$ are ordered in 
such a way that $S$ is an upper triangular matrix. We choose the standard
generators $\sigma_1$, \dots, $\sigma_{n-1}$ of the braid group $\bn$. The
generator $\sigma_i$ is given by a deformation of $(u_1, \dots, u_n)$ 
such that:

1). $u_k$ remains fixed for $k\neq i,\,i+1$.

2). $u_i$ and $u_{i+1}$ are permuted moving counterclockwise.

Let us deform $(u_1, \dots, u_n)$ in the coefficients of the operator
$$
L={d\over dz} -\left( U+{V(u)\over z}\right).
$$
Due to isomonodromicity the matrices $S$ and $C$ remain unchanged until 
some of the Stokes rays passes through $\ell$. After this we are to 
reorder the canonical coordinates to preserve upper triangularity of the 
Stokes  matrix and, then, to compute the new matrices $S'$ and $C'$ using
Lemma 4.8. We are to recall here that the operator $L$ for a given
ordering 
of the canonical coordinates $(u_1, \dots, u_n)$ is determined up to a 
transformation
$$
L\mapsto J\,L\,J
\eqno(4.69)
$$
where $J$ is an arbitrary diagonal matrix of the form
$$
J=\diag (\pm 1, \dots, \pm 1).
\eqno(4.70)
$$
Thus the matrices
$$
S ~{\rm and }~J\,S\,J, ~~C~{\rm and}~J\, C
\eqno(4.71)
$$
must be identified. So what we need is actually an action of ${\cal B}_n$
on the 
classes of equivalence of the matrices $S$ and $C$ w.r.t. the identifications
(4.71).

The result is given by
\smallskip
{\bf Theorem 4.6.} {\it The analytic continuation of a semisimple
Frobenius 
manifold is described by the following action
$$
\eqalign{
S &\mapsto \beta (S)\cr
C &\mapsto \beta (C)\cr}
\eqno(4.72)
$$
of the braid group $\bn\ni\beta$ on the Stokes matrix $S=(s_{ij})$ and 
the central connection matrix $C$. For the standard generator 
$\beta=\sigma_i$ the action has the form
$$
\eqalign{
\sigma_i(S) &=K^{(i)}(S)\, S\, K^{(i)}(S)\cr
\sigma_i(C) &= K^{(i)}(S)\, C\cr}
\eqno(4.73)
$$
where
$$
\eqalign{
\left( K^{(i)}(S)\right) _{kk} &=1, ~k=1, \dots, n, ~~k\neq i, \, i+1\cr
\left( K^{(i)}(S)\right)_{i+1,\, i+1} &=-s_{i,\,i+1}\cr
\left( K^{(i)}(S)\right)_{i,\,i+1} &= \left( 
K^{(i)}(S)\right)_{i+1,\,i}=1\cr}
\eqno(4.74)
$$
all other entries of the matrix $ K^{(i)}(S)$ are equal to 
zero.}

{\bf Proof. } Let us assume that, during the deformation $\sigma_i$, the 
coordinates $u_i$ and $u_{i+1}$ remain sufficiently close to each other. 
Then all the Stokes rays but $R_{i,\,i+1}$ and $R_{i+1,\, i}$ will be 
only slightly deformed and they will return to their original positions
(with renumbering $i\leftrightarrow i+1$) after the end of the deformation.
But the rays  $R_{i,\,i+1}$ and $R_{i+1,\, i}$ interchange their positions
rotating clockwise. Particularly, it is the ray $R = R_{i+1,\, i}$
who passes through the positive half-line $\ell_+$ rotating clockwise.
At the very last moment before the collision the configuration of the 
Stokes rays is such that $R^{(1)} = R_{i+1,\, i}$ and $R^{(m+1)}= 
R_{i,\,i+1}$, and we may assume that $R^{(1)}$ and $R^{(m+1)}$ contain no 
other Stokes rays. From Theorem 4.5 we obtain a factorization of $S$ into
the 
product of upper triangular Stokes factors
$$
S= K_1 K_2 \dots K_m
$$
where the only nonzero off-diagonal entry of the matrix $K_1$ sits in the 
$(i,i+1)$ box, and all the factors $K_2$, \dots, $K_m$ have zero on the 
$(i,i+1)$ place. From this we obtain that
$$
\left( K_1\right)_{i,\,i+1} =s_{i,\,i+1}.
$$
We are now to apply the formulae (4.68) to compute the new matrices $S'$,
$C'$
with
$$
K_R =K_{m+1} =\left(K_1^T\right)^{-1}.
$$
After this applying the permutation $i\leftrightarrow i+1$ we arrive at 
the formulae (4.74). Theorem is proved.
\medskip
{\bf Example 4.3.} For $n=3$ the generators $\sigma_1$, $\sigma_2$ of
${\cal 
B}_3$ act as follows in the space of Stokes matrices
$$
S=\left(\matrix{1 & x & y\cr
0 & 1 & z\cr
0 & 0 & 1\cr}\right)
$$
$$
\sigma_1(x,y,z) =(-x,z,y-x\,z), ~~\sigma_2(x,y,z) =(y, x-y\,z,-z).
\eqno(4.75)
$$
\smallskip
{\bf Exercise 4.7.} Prove that the braid
$$
\zeta = (\sigma_1 \dots \sigma_{n-1})^n
\eqno(4.76)
$$
acts trivially on Stokes matrices.
\medskip
The braid $\zeta$ generates the center of $\bn$ (see [Bi]). So
the quotient $\bn/{\rm center}$ acts on the space of Stokes matrices.
For $n=3$ the quotient is isomorphic to the modular group $PSL_2({\bf Z})$
[{\it ibid}].

Let $\bn(S,C)\subset \bn$ be the stationary subgroup of the class of 
equivalence (4.71) of the pair $S$, $C$. We realize it as a subgroup
in the fundamental group
$$
\pi_1\left( \left[\C^n\setminus\diag\right]/S_n, (u_1^0, \dots, u_n^0)\right)
$$
and construct the corresponding covering
$$
M(S,C)\to  \left[\C^n\setminus\diag\right]/S_n
$$
i.e., such a covering that the group of deck transformations of the fiber
is isomorphic to $\bn(S,C)$. From Theorem 4.6 it follows 
\smallskip
{\bf Theorem 4.7.} {\it 1). For a given monodromy data $(\mu, e, R, S, C)$
the Frobenius structure extends from a small \nbh of $u^0$ to a dense 
open subset in the manifold $M(S,C)$. This Frobenius structure on $M(S,C)$
we denote $Fr(\mu, e, R, S, C)$.

2). Let $(\mu, e, R, S, C)$ be the monodromy data of a semisimple Frobenius
manifold $M$ computed at the point $u^0=(u^0_1, \dots, u^0_n)$ w.r.t. an
admissible oriented line $\ell$. Let $M^0$ be the open part of the 
Frobenius manifold $M$ consisting of all points $t\in M$ such that all the 
eigenvalues $u_1(t)$, \dots, $u_n(t)$ of the operator of multiplication
by the Euler vector field are pairwise distinct. Then the map
$$
M^0 \to Fr(\mu, e, R, S, C)
$$
is well-defined and it is an equivalence of Frobenius manifolds.}
\medskip
{\bf Example 4.4.} Let us compute the monodromy data of quantum cohomology
of ${\bf CP}^2$, i.e., of the solution (1.15) of \wdvv. The monodromy
at $z=0$ is completely determined by the classical cohomology 
$H^*\left({\bf CP}^2\right)$ together with the first Chern class 
$c_1\left({\bf CP}^2\right)$ (see Lecture 2). We obtain
$$
\mu =\diag (-1, 0, 1), ~~R=\left(\matrix{0 & 0 & 0\cr 3 & 0 & 0\cr 0 & 3 &
0\cr}\right).
$$

Let us compute the Stokes matrix in the semisimple point
$$
t_1=t_3=0, ~~{\rm arbitrary}~t_2 ~{\rm with}~{\rm Re}\, t_2<R.
\eqno(4.77)
$$
Here $R$ is the radius of convergence (1.16). Let us denote $q=\exp t_2$.
The system (2.41) for horizontal sections $\left(\xi_1, \xi_2,
\xi_3\right)
=\left ( \partial_1\tilde t, \partial_2 \tilde t, \partial_3 \tilde t\right)$
of the connection $\nab$ can be reduced to two third-order equations
$$
\eqalign{
\partial_2^3 \phi &= z^3 q\,\phi\cr
(z\partial_z)^3 \phi &= 27 z^3 q\, \phi\cr}
\eqno(4.78)
$$
for the function
$$
\phi =\phi(t_2 z) = {\xi_1\over z},
$$
$$
\left(\xi_1, \xi_2, \xi_3\right)=\left( z\, \phi, {1\over 3} z\partial_z 
\phi, {1\over 9} \partial_z (z\partial_z\phi)\right).
$$
The system (4.78) is equivalent to one equation
$$
(z\partial_z)^3 \Phi = 27 z^3 \Phi
\eqno(4.79)
$$
using the quasihomogeneity
$$
\phi(t_2, z) =\Phi\left( z \, q^{1/3}\right).
\eqno(4.80)
$$
The problem is reduced to computation of the Stokes matrix of the 
generalized hypergeometric equation (see [DM]). We are to carefully select 
the basis of formal solutions of (4.79) at $z\to\infty$ correspondent to
the 
basis of columns of $Y_\fo(z)$ of the solution (4.26) of the 
gauge-equivalent system (3.30).

The operator $\U$ of multiplication by the Euler vector field in the basis
$e_1 =\partial_1$, $e_2 =\partial_2$, $e_3 =\partial_3$ has the matrix
$$
\U(t) =\left( \matrix{ 0 & 0 & 3q\cr 3 & 0 & 0\cr 0 & 3 & 0 \cr}\right) ~~
t=(0, t_2, 0), ~q=e^{t_2}.
\eqno(4.81)
$$
The canonical coordinates (i.e., the eigenvalues of $\U$) in the point 
(4.77) take the values
$$
u_1 = 3 q^{1/3}, ~u_2 =3 \bar\epsilon^2 q^{1/3}, ~ u_3 =3 \epsilon^2 
q^{1/3}
\eqno(4.82)
$$
where
$$
\epsilon =\exp{\pi i\over 3}.
$$
The correspondent idempotents of the quantum cohomology algebra are
$$
\eqalign{
\pi_1 &= {1\over 3}\left( e_1 +q^{-1/3} e_2 +q^{-2/3}e_3\right)\cr
\pi_2 &= {1\over 3}\left( e_1 +\epsilon^2 q^{-1/3} e_2 +
\bar\epsilon^2 q^{-2/3}e_3\right)\cr
\pi_3 &= {1\over 3}\left( e_1 +\bar\epsilon^2 q^{-1/3} e_2 
+ \epsilon^2 q^{-2/3}e_3\right)\cr}.
$$
The invariant metric 
$$
<\pi_1, \pi_1> ={1\over 3} q^{-2/3}, ~~<\pi_2, \pi_2> ={1\over 3} 
\bar\epsilon^2 q^{-2/3}, ~~<\pi_3, \pi_3>={1\over 3} \epsilon^2 q^{-2/3}.
$$
Evaluating the square root we obtain the normalized idempotents
$$
\eqalign{
f_1 &= {1\over \sqrt{3}}\left( q^{1/3} e_1 + e_2 
+q^{-1/3}e_3\right)\cr   
f_2 &= {1\over \bar\epsilon \sqrt{3}}\left( q^{1/3} e_1 +\epsilon^2  
e_2 + \bar\epsilon^2 q^{-1/3}e_3\right)\cr
f_3 &= {1\over \epsilon \sqrt{3}}\left( q^{1/3} e_1 +\bar\epsilon^2 
e_2 + \epsilon^2 q^{-1/3}e_3\right)\cr}.
$$
This gives the matrix $\Psi=\left(\psi_{i\alpha}\right)$
$$
\Psi ={1\over \sqrt{3}} \left(\matrix{
q^{-1/3} & 1 & q^{1/3}\cr
\bar\epsilon q^{-1/3} & -1 & \epsilon q^{1/3}\cr
\epsilon q^{-1/3} & -1 & \bar\epsilon q^{1/3}\cr}
\right).
\eqno(4.83)
$$
We can easily compute the matrix $V$ in the point of interest (cf.
[MM]). 
But what we need is to determine the asymptotic structure of the 
solutions of (4.79) at $z\to\infty$. We must choose the basis $\tilde 
t_1^\infty$,  $\tilde
t_2^\infty$,  $\tilde
t_3^\infty$ of the coordinates $\tilde t$ such that the matrix
$$
Y_{ij}:={\partial_i \tilde t_j^\infty\over \psi_{i1}}
$$
has the development (4.26), i.e.,
$$
Y_{ij
} \sim \left(\delta_{ij} +O\left({1\over z}\right)\right) \, e^{z\, u_j}, 
~~i, \, j=1, \, 2, \, 3.
\eqno(4.84)
$$
This gives the three solutions $\phi_1$,  $\phi_2$, $\phi_3$ of the system 
(4.78) such that
$$
\phi_j ={1\over z} {\partial\over \partial t^1} \tilde t_j^\infty =
{1\over z} \sum_{i=1}^3 \partial_i \tilde t_j^\infty ={1\over z} \sum_{i=1}^3
\psi_{i1} Y_{ij}.
$$
For the correspondent basic solutions of (4.79) we obtain the needed 
developments
$$
\eqalign{
\Phi_1 
&\sim 
{1\over \sqrt{3}} 
{e^{3z}\over z} 
\left( 1+O\left({1\over z}\right)\right)\cr
\Phi_2 
&\sim 
{\bar\epsilon\over \sqrt{3}} 
{e^{3\bar\epsilon^2 z}\over z} 
\left(1+O\left({1\over z}\right)\right)\cr
\Phi_3 
&\sim 
{\epsilon\over \sqrt{3}} 
{e^{3\epsilon^2 z}\over z}
\left( 1+O\left({1\over z}\right)\right).\cr}
\eqno(4.85)
$$
We are now to compute the Stokes matrix of the equation (4.79) with
respect to
the bases of solutions  having the asymptotic developments (4.85) in the 
$\rr/\ll$ half-planes $\Pi_{\rr/\ll}$ with some admissible oriented line
$\ell$. 

The Stokes rays of equation (4.79) have the form
$$
\eqalign{
R_{12} &= \{ -\rho\,\epsilon \, | \, \rho\geq 0\}\cr
R_{13} &= \{ \rho\,\bar\epsilon \, | \, \rho\geq 0\}\cr
R_{23} &= \{ \rho ~~| \, \rho \geq 0\}\cr}
\eqno(4.86)
$$
the rays $R_{21}$, $R_{31}$, $R_{32}$ are the opposite to the above. We 
choose the admissible line
$$
\ell =\{ r\, e^{i\alpha} \, | \, -\infty < r < \infty\}
\eqno(4.87)
$$
for a fixed small $\alpha>0$ oriented according to the positive direction 
of $r$. We will use now a suitable Meijer function [Lu] to compute
the Stokes matrix.
\smallskip
{\bf Lemma 4.9.} {\it The function
$$
g(z) ={1\over (2\pi )^2 i} \int_{-c-i\infty}^{-c+i\infty} \Gamma^3(-s) 
e^{\pi i s} z^{3s} ds
\eqno(4.88)
$$
defined for $z\neq 0$,
$$
-{5 \pi\over 6} <\arg z < {\pi \over 6}
\eqno(4.89)
$$
where $c$ is any positive number, satisfies (4.79). The analytic 
continuation of this function has the asymptotic development
$$
g(z)\sim {1\over \sqrt{3}} \bar\epsilon {e^{3\bar\epsilon^2 z}\over z} 
=\Phi_2(z), ~~ |z|\to \infty
\eqno(4.90)
$$
in the sector
$$
-{5\pi\over 3} < \arg z < \pi.
\eqno(4.91)
$$
It satisfies the identity
$$
g\left(z\, e^{2\pi i}\right) -3\, g\left( z\, e^{4\pi i\over 3}\right)
+ 3\, g\left( z\, e^{2\pi i\over 3}\right) -g \left( z \right) =0.
\eqno(4.92)
$$}

{\bf Proof} (cf [Lu]). Using the Stirling formula
$$
\log \Gamma(z) = \left( z-{1\over 2} \right) \log z - z + {1\over 2} \log
(2\pi) + O\left( {1\over z} \right) 
$$
and
$$
\lim_{|y|\to\infty} |\Gamma(x+iy)| e^{{\pi\over 2} {y}} |y|^{{1\over 2}-x}
=\sqrt{2\pi}, ~x, \, y ~{\rm real}
$$
we prove uniform convergence of the integral in the domain (4.89) and 
independence it on $c$. Differentiation gives
$$
(z\partial_z)^3 g = 
{27\over (2\pi )^2 i} \int_{-c-i\infty}^{-c+i\infty} s^3 \Gamma^3(-s)
e^{\pi i s} z^{3s} ds.
$$
Using the property of gamma-function
$$
s\, \Gamma(-s) =-\Gamma(1-s)
$$
and doing a shift $s\to s+1$ we obtain for the r.h.s. the integral
$$
{27\over (2\pi )^2 i} \int_{-c-i\infty}^{-c+i\infty} \Gamma^3(-s)
e^{\pi i s} z^{3(s+1)} ds
=27 z^3 g.
$$
To derive the asymptotic development we use Laplace method. Representing
the integrand in the form $\exp phase$ and using Stirling formula one obtains
the following asymptotic development for
$$
\eqalign{
phase &= 3\log \Gamma(-s) +\pi i s + 3 s \log z\cr
 &\sim -3\left( s+{1\over 2}\right) \log s + 3 s \log z + (3-2\pi i)s\cr}
$$
valid for
$$
-{3\pi\over 2}<\arg s<-{\pi\over 2}.
\eqno(4.93)
$$
For big $|z|$ the phase has critical point at
$$
s\sim z\, e^{-{2\pi i\over 3}} -{1\over 2}.
$$
This critical point is in the domain (4.89) if
$$
-{5\pi\over 6} < \arg z < {\pi \over 6}.
$$
In the critical point the value of the phase is
$$
phase_0 \sim -{3\over 2} \log z + 3 z e^{-{2\pi i\over 3}} + {3\over 2} 
\log 2 \pi + \pi i
$$
and the second $s$-derivative at this point
$$
phase''_0 \sim -{3 e^{2\pi i\over 3}\over z}.
$$
Applying Laplace formula for the integral
$$
g(z) \sim {1\over (2\pi)^2 i} {1\over \sqrt{2\pi}} {e^{phase_0}\over
\sqrt{phase''_0}}
$$
we obtain the asymptotics (4.90). The asymptotics remains valid in a wider 
sector
$$
-{5\pi\over 3} < \arg z < \pi.
$$
Indeed, during this analytic continuation, i.e., counterclockwise until
$R_{32}$ and clockwise until $R_{21}$ the exponential $e^{z\, u_2}$
remains dominant.

To derive the identity (4.92) we observe that the equation is invariant
w.r.t. the rotation
$$
z\mapsto z\, e^{2\pi i\over 3}.
$$
This generates a linear operator, $A$, in the 3-dimensional space of 
solutions of (4.79). Let us prove that all the eigenvalues of $A$ are
equal 
to 1. Indeed, near $z=0$ all the solutions have the form
$$
\Phi(z) = \sum_{m=0}^\infty
{z^{3m}\over (m!)^3} \left[
a_m + b_m \log z + c_m \log^2 z\right]
\eqno(4.94)
$$
where $a_0$, $b_0$, $c_0$ are arbitrary parameters and the coefficients
$a_m$, $b_m$, $c_m$ for $m>0$ are uniquely determined from the recursion 
relations
$$
\eqalign{
c_m &= c_{m-1}\cr
b_m +{2\over m} c_m &= b_{m-1}\cr
a_m +{1\over m} b_m +{2\over 3 m^2} c_m &= a_{m-1}.\cr}
$$
The operator
$$
\left( A\phi\right)(z) = \Phi \left( z\, e^{2\pi i\over 3}\right)
$$
in the basis of solutions of the form (4.94) with only one nonzero of
$a_0$,
$b_0$, $c_0$  is given by a triangular matrix with all 1 on the 
diagonals. Writing Cayley - Hamilton theorem
$$
\left( A-1\right)^3 =0
$$
we obtain
$$
A^3 g - 3 A^2 g + 3 A\ g -g =0.
$$
This gives the identity (4.92). Lemma is proved.
\medskip
Let us construct the three solutions
$\Phi^\rr(z) =\left( \Phi_1^\rr(z), \Phi_2^\rr(z), \Phi_3^\rr(z)\right)$ 
having the asymptotic behaviour of the form (4.85)
$$
\Phi_j^\rr(z) \sim \Phi_j(z), ~~|z|\to\infty, ~-\pi <\arg z< {\pi\over 
3}, ~j=1, \, 2, \, 3.
$$
We can take
$$
\Phi^\rr(z) =\left( -g\left( e^{2\pi i\over 3}z\right), g\left(z\right), 
g\left( e^{-{2\pi i\over 3}}z\right)\right).
\eqno(4.95)
$$
Similarly, the components of the vector-function $\Phi^\ll(z)$ must have the
asymptotics 
$$
\Phi_j^\ll(z) \sim \Phi_j(z), ~~|z|\to\infty, ~0 <\arg z< {4\pi\over
3}, ~j=1, \, 2, \, 3.
$$
We take
$$
\Phi^\ll(z) = \left( -g\left( e^{-{4\pi i\over 3}}z\right), g\left( 
e^{-{2\pi i}} z\right)-3 g\left( e^{-{4\pi i\over 3}} z\right),
g\left( e^{-{2\pi i\over 3}}z\right)\right).
\eqno(4.96)
$$
The only novelty to be proved is the formula for $\Phi_2^\ll$. Indeed, 
from Lemma 4.9 it follows that
$$
\Phi_2^\ll (z) =  g\left(
e^{-{2\pi i}} z\right)-3 g\left( e^{-{4\pi i\over 3}} z\right)\sim \Phi_2(z),
~|z|\to\infty, ~{\pi\over 3}<\arg z<{4\pi\over 3}.
$$
Using the identity (4.92) we may rewrite this function as
$$
\Phi_2^\ll(z) = g\left( z\right) - 3 g\left( e^{-{2\pi i\over 3}}z\right) 
\sim \Phi_2(z), ~|z|\to\infty, ~0<\arg z<{\pi\over 3}.
$$

Applying again the identity (4.92) we obtain that in the sector
$$
0< \arg z< {\pi \over 3}
$$
$$
\left( \Phi_1^\ll(z), \Phi_2^\ll(z), \Phi_3^\ll(z)\right) =
\left( \Phi_1^\rr(z), \Phi_2^\rr(z), \Phi_3^\rr(z)\right) \, S
$$
with 
$$
S=\left( \matrix{ 1 & 0 & 0\cr
3 & 1 & 0\cr
-3 & -3 & 1\cr}\right).
\eqno(4.97)
$$
This is the Stokes matrix of the quantum cohomology of ${\bf CP}^2$. 
Changing the sign of the normalized idempotent $f_3$ we can reduce $S$ to 
the form
$$
S=\left( \matrix{ 1 & 0 & 0\cr
3 & 1 & 0\cr
3 & 3 & 1\cr}\right).
\eqno(4.97')$$

The matrix (4.97) was obtained from physical considerations in [CV2].
The main argument was that, in Landau - Ginzburg models of 2D topological
field theory the entries of the Stokes matrix must be integers. Then,
since the eigenvalues of $S^T S^{-1}$ must all be 1, one arrives at the
following Diophantine equation for the entries
$$
x^2 + y^2 + z^2 - x y z = 0
$$
where
$$
S=\left(\matrix{1 & x & y\cr 0 & 1 & z \cr 0 & 0 & 1\cr}\right).
$$
All the integer solutions to the equation have the form
$$
x = 3 x_1, ~y=3 y_1, ~ z=3 z_1
$$
where $x_1, y_1, z_1$ are integer solutions to Markoff equations
$$
x_1^2 + y_1^2 + z_1^2 - 3 x_1 y_1 z_1=0.
$$
The solutions of Markoff equation are known to be all equivalent to
$(1,1,1)$ modulo the action (4.75) of the braid group. This solution of
Markoff equation just corresponds to the Stokes matrix (4.97')
 \bigskip
In the next Lecture we will construct polynomial Frobenius manifolds 
starting from an arbitrary finite Coxeter group. Particularly, for the
Coxeter groups with simply-laced Dynkin diagrams these coincides with
the Frobenius manifolds of the singularity theory.
It can be shown, using this construction, that, in this examples, $S$ is 
the variation operator of the singularity computed in 
the so-called marked basis of vanishing cycles [AGV].

Here we define a remarkable operation of {\it tensor product} of 
Frobenius manifolds. We are motivated by the results of
Kaufmann, Kontsevich, 
and Manin [KM, Ka] describing quantum cohomology of the direct 
product of two varieties.

Let $M'$, $M''$ be two Frobenius manifolds of the dimensions $n'$ and $n''$
resp. We say that a Frobenius manifold $M$ of the dimenion $n'\, n''$ is 
the tensor product $M=M'\otimes M''$ if it has the following structure.

1). The tangent planes $TM$ with the bilinear form $<~,~>$ and the unity 
vector field $e$ are represented as
$$
\left( TM, <~,~>, e\right) = \left( TM'\otimes TM'', 
<~,~>'\otimes  <~,~>'', 
e' \otimes e''\right) 
$$
(as usually, we identify the tangent planes in different points using the 
Levi-Civita flat connection). Thus, the flat coordinates on $M$ have 
double labels
$$
t=\left( t^{\alpha'\,\alpha''}\right), ~~1\leq \alpha'\leq n', ~1\leq 
\alpha''\leq n''.
$$
The unity vector field is
$$
e={\partial\over \partial t^{1'\, 1''}}.
$$
The matrix of $<~,~>$ has the form
$$
\eta_{\alpha'\alpha''\, 
\beta'\beta''}=\eta_{\alpha'\beta'}\eta_{\alpha''\beta''}.
$$

2). In the points
$$
t\in M, ~~t^{\alpha'\alpha''}=0 ~{\rm for}~ \alpha'>1, ~\alpha''>1
\eqno(4.98)
$$
the algebra $T_tM$ is the tensor product
$$
T_tM= T_{t'}M'\otimes T_{t''}M'',
$$
$$
\eqalign{
t' &= \left( t^{2'1''}, \dots, t^{n'1''}\right)\cr
t'' &= \left( t^{1' 2''}, \dots, t^{1'n''}\right)\cr}
$$
i.e., 
$$
c_{\alpha'\alpha''\, \beta'\beta''}^{\gamma'\gamma''}(t) =
c_{\alpha'\beta'}^{\gamma'}(t') c_{\alpha''\beta''}^{\gamma''}(t'').
$$
In these formulae $\eta_{\alpha'\beta'}$, $c_{\alpha'\beta'}^{\gamma'}$
and $\eta_{\alpha''\beta''}$, $ c_{\alpha''\beta''}^{\gamma''}$ are the 
invariant bilinear form and the structure constants of the Frobenius 
manifolds $M'$ and $M''$ resp.

3). The charge
$$
d_M = d_{M'} + d_{M''}
$$
and the Euler vector field on $M$ has the form
$$
E=\sum_{\alpha',\, \alpha''}t^{\alpha'\alpha''} \left( 
1-q_{\alpha'}-q_{\alpha''}\right) {\partial\over \partial 
t^{\alpha'\alpha''}}
+\sum r_{\alpha'} {\partial\over\partial t^{\alpha'1''}} + \sum 
r_{\alpha''} {\partial\over\partial t^{1'\alpha''}}.
\eqno(4.99)
$$
Here
$$
\eqalign{
E' &= \sum_{\alpha'=1}^{n'} \left[ (1-q_{\alpha'}) t^{\alpha'} + 
r_{\alpha'}\right] \partial_{\alpha'}\cr
E'' &= \sum_{\alpha''=1}^{n''} \left[ (1-q_{\alpha''}) t^{\alpha''} +
r_{\alpha''}\right] \partial_{\alpha''}\cr}
$$
are the Euler vector fields on $M'$ and $M''$ resp.

For any two semisimple Frobenius manifolds $M'$, $M''$
we will now describe
their tensor product $M=M'\otimes M''$ in terms of the monodromy data of
the factors.
\smallskip
{\bf Lemma 4.10.} {\it 1). If $M=M'\otimes M''$ with semisimple $M'$ and
$M''$
then $M$ is semisimple.

2). Let $t_0'\in M'$, $t_0''\in M''$ be two points such that {\rm a)} 
${t_0^1}'={t_0^1}''$, and {\rm b)} the values of the canonical coordinates
$u_{i'}=u_{i'}(t_0')$, $i'=1, \dots, n'$, $u_{i''}=u_{i''}(t_0'')$, $i''=1, 
\dots, n''$ satisfy the properties
$$
\eqalign{u_{i'} &\neq u_{j'}, ~i'\neq j'\cr
u_{i''} &\neq u_{j''}, ~ i''\neq j''\cr
u_{i'} +u_{i''} &\neq u_{j'} +u_{j''}, ~(i',\, i'') \neq (j', \, j'').\cr}
$$
Let $\ell$ be a line on the $z$-plane such that for any $z\in 
\ell\setminus 0$
$$
\eqalign{
{\rm Re}\, \left[ z\, (u_{i'} -u_{j'})\right] &\neq 0, ~i'\neq j'\cr
{\rm Re}\, \left[ z\, (u_{i''} -u_{j''})\right] &\neq 0, ~i''\neq j''\cr
{\rm Re}\, \left[ z\, (u_{i'} -u_{j'})\right] +
{\rm Re}\, \left[ z\, (u_{i''} -u_{j''})\right] &\neq 0, ~(i',\, i'')\neq 
(j',\, j'').\cr}
$$
Then the Stokes matrix $S$ of $M$ in the point $t_0$ with the coordinates
$$
\eqalign{ t^{\alpha'1''} &= t_0^{\alpha'}, ~\alpha'=1, \dots, n'\cr
 t^{1'\alpha''} &= t_0^{\alpha''}, ~\alpha''=1, \dots, n''\cr
t^{\alpha'\alpha''} &=0, ~\alpha'>1, ~\alpha''>1\cr}
\eqno(4.100)
$$
is the tensor product of the Stokes matrices $S'$ of $M'$ in the point 
$t_0'$ and $S''$ of $M''$ in the point $t_0''$
$$
S=S'\otimes S''.
$$}

{\bf Proof.} If $t_0'\in M'$, $t_0''\in M''$ are semisimple points of the 
Frobenius manifolds then the point (4.100) will be a semisimple point of
$M$. 
The idempotents of the algebra
$$
T_{t_0'}M'\otimes T_{t_0''}M''
$$
are tensor products $\pi_{i'}\otimes\pi_{i''}$, $i'=1, \dots, n'$, 
$i''=1, \dots, n''$. The operator of multiplication by the Euler vector 
field (4.99) in the point (4.100) has the form
$$
\U = 1'\otimes \U'' + \U'\otimes 1'' - t^1 \, 1'\otimes 1''
$$
where $t^1 = t_0^{1'} =t_0^{1''}$. The eigenvalues of this operator are
$$
u_{i'} +u_{i''} - t^1, ~~1\leq i'\leq n', ~1\leq i''\leq n''.
$$
These are the values of the canonical coordinates on $M$ in the points of 
the $(n'+n''-1)$-dimensional locus (4.98).

Let $Y_{\rr/\ll}'(z; t_0')$, $Y_{\rr/\ll}''(z; t_0'')$ be the solutions 
of the system (3.30) for $M'$ and $M''$ resp. with the asymptotic
behaviour 
(4.26) in the right/left half-planes w.r.t. the admissible line $\ell$.
Then 
the solutions of the system (3.30) for $M$ with the needed asymptotic 
development (4.26) are
$$
\eqalign{
Y_\rr (z; t_0) &= e^{-z\, t^1} Y_\rr'(z; t_0') \otimes Y_\rr''(z; t_0''), 
~z\in \Pi_\rr\cr
Y_\ll (z; t_0) &= e^{-z\, t^1} Y_\ll'(z; t_0') \otimes Y_\ll''(z; t_0''),
~z\in \Pi_\ll\cr}
$$
This proves Lemma.
\medskip
{\bf Theorem-Definition 4.8.} {\it Let
$$
M=Fr \left( \mu'\otimes 1+1\otimes \mu'', e'\otimes e'', R'\otimes 
1+1\otimes R'', S'\otimes S'', C'\otimes C''\right)
$$
$$
M'=Fr\left( \mu', e', S', C'\right)
$$
$$
M''=Fr\left( \mu'', e'', S'', C''\right).
$$
Then
$$
M=M'\otimes M''.
$$}

{\bf Proof.} Let $u'=(u_{1'}. \dots, u_{n'})\in M'$,  $u''=(u_{1''}. \dots, 
u_{n''})\in M''$ be two regular points of these Frobenius manifolds, 
i.e., such points that the \rh of the form (4.41) - (4.45) for each of the
manifolds
has unique solution $\left( Y_0', Y_\rr', Y_\ll'\right)$ and  $\left( 
Y_0'', Y_\rr'', Y_\ll''\right)$ resp. Doing, if necessary, a diagonal shift
$$
u_{i'}\mapsto u_{i'}+c, ~i'=1, \dots, n'
$$
we may also assume that
$$
t^{1'} (u') = t^{1''}(u'') =: t^1.
$$
Then the functions
$$
\eqalign{
Y_0 &= e^{-z\, t^1} Y_0'\otimes Y_0''\cr
Y_\rr &= e^{-z\, t^1} Y_\rr'\otimes Y_\rr''\cr
Y_\ll &= e^{-z\, t^1} Y_\ll'\otimes Y_\ll''\cr}
$$
will give the solution of the \rh for the manifold $M$. It follows that 
the matrix $\Psi$ in the points (4.100) is also a tensor product
$$
\Psi =\left( \psi_{i'\alpha'}(u') \psi_{i''\alpha''}(u'')\right).
$$
Using the formulae of Isomonodromicity Theorem we conclude that 
$M=M'\otimes M''$. Theorem is proved.
\medskip
{\bf Example 4.5.} Let $M$ be the Frobenius manifold corresponding to 
quantum cohomology of ${\bf CP}^1$, i.e.,
$$
F={1\over 2} t_1^2 t_2 + e^{t_2}
$$
$$
E = t_1 \partial_1 + 2 \partial_2.
$$
The monodromy data are
$$
\mu = \diag \left( -1/2, 1/2\right), ~R=\left(\matrix{0 & 0\cr 2 & 
0\cr}\right), S=\left( \matrix{1 & 2\cr 0 & 1\cr}\right)
$$
(the computation of $S$ is similar to the above computation of the Stokes 
matrix of quantum cohomology of ${\bf CP}^2$, but it is simpler). The 
tensor square of this Frobenius manifold computed according to Theorem 
describes the quantum cohomology of ${\bf CP}^1 \times {\bf CP}^1$.
\medskip
{\bf Example 4.6.} Let $M_h$ be the polynomial two-dimnsional Frobenius  
manifolds of the form
$$
F={1\over 2} t_1^2 t_2 + t_2 ^{h+1}, ~h\in {\bf Z}, ~h\geq 3.
$$
Tensor product of the form $M_{h'}\otimes M_{h''}$ is a polynomial 
4-dimensional Frobenius manifold only in the following three cases:
$M_3\otimes M_3$, $M_3\otimes M_4$, $M_3\otimes M_5$. 

In the next Lecture we will establish a relation between polynomial 
Frobenius manifolds and finite Coxeter groups. We will see that the 
manifolds $M_h$ correspond to the groups $I_2(h)$ of symmetries of 
regular $h$-gon on the plane. Particularly, for $h=3$ we obtain $I_2(3)=$ 
the Weyl group of the type $A_2$, for $h=4$ $I_2(4)=$ the Weyl group
of the type $B_2$. Their tensor products also correspond to certain 
finite Coxeter groups. Namely,
$$
\eqalignno{
M_{A_2}\otimes M_{A_2} &= M_{D_4}
& {(4.101)}
\cr
M_{A_2}\otimes M_{B_2} &= M_{F_4}
&{(4.102)}
\cr
M_{A_2}\otimes M_{I_2(5)} &=M_{H_4}
&{(4.103)}\cr}
$$
the notations for finite Coxeter groups as in [Bou]; see also the next 
Lecture). Besides these there are only two more cases where tensor 
products of
two polynomial Frobenius manifolds is again a polynomial Frobenius
manifold.
They correspond to the following Coxeter groups
$$
\eqalignno{
M_{A_2}\otimes M_{A_3} &=M_{E_6}
&{(4.104)}
\cr
M_{A_2}\otimes M_{A_4} &= M_{E_8}.
&{(4.105)}
\cr}
$$

More generally, in the singularity theory our operation of tensor 
products of the Frobenius structures on the parameter space of versal 
deformation of an isolated quasihomogeneous singularity corresponds to 
the operation of the direct sum of singularities. Denoting $M_{f(x)}$
the Frobenius structure on the parameter space of versal deformation of 
the singularity of a function $f(x)$ we obtain
$$
M_{f(x)+g(y)} = M_{f(x)} \otimes M_{g(y)}.
$$
Indeed, according to Deligne (see in [AGV]) the variation operator of the 
direct sum of the singularities is the tensor product of the variation 
operators of the summands. From this point of view the identifications 
(4.101), (4.104), (4.105) become obvious. The equalities (4.102) and
(4.103) seem not to 
admit simple explanation within the framework of the singularity theory.
However, they are in the agreement with the embeddings of Frobenius
manifolds obtained by folding of Dynkin diagrams explained in the next
Lecture (I am thankful to J.-B.Zuber for bringing my attention to this
point).

\vfill\eject
\def\res{\mathop{\rm res}}
\def\rr{{\rm right}}
\def\ll{{\rm left}}
\def\arg{{\rm arg}\,}
\def\nbh{{neighborhood }}
\def\tr{{{\rm tr}\,}}
\def\diag{{{\rm diag}\,}}
\def\L{{{\cal L}_E}}
\def\C{{\bf C}}
\def\U{{\cal U}}
\def\E{{\cal E}}
\centerline{Lecture 5}
\medskip
\centerline{\bf Monodromy group and mirror construction for 
semisimple Frobenius manifolds}
\medskip
We will introduce a new metric [Du5, Du7] on an open subset of a Frobenius
manifold 
$M$. The inverse of this metric will be a symmetric bilinear form on the 
cotangent bundle $T^*M$ defined everywhere.
\smallskip
{\bf Definition 5.1.} The {\it intersection form} of the Frobenius
manifold 
$M$ is the bilinear form on $T^*M$ defined by the formula
$$
(\omega_1, \omega_2) := i_{E(t)} (\omega_1\cdot \omega_2), ~~\omega_1, \, 
\omega_2\in T_t^*M.
\eqno(5.1)
$$

In the r.h.s. the product of one-forms $T_t^*M\times T_t^*M \to T_t^*M$ is 
defined using the algebra structure on $T_tM$ and the isomorphism
$$
<~,~> : T_tM \to T_t^*M.
$$
\medskip
In the flat coordinates the components of the intersection form are given 
by the formula
$$
\eqalign{
g^{\alpha\beta}(t) := \left( dt^\alpha, dt^\beta\right) &= E^\epsilon (t) 
c_\epsilon^{\alpha\beta}(t)\cr
 &= (d+1-q_\alpha-q_\beta) F^{\alpha\beta}(t) +A^{\alpha\beta}.\cr}
\eqno(5.2)
$$
Here
$$
\eqalign{
c_\epsilon^{\alpha\beta}(t) &=
\eta^{\alpha\gamma}c_{\gamma\epsilon}^\beta(t)\cr
F^{\alpha\beta}(t) &= \eta^{\alpha\lambda}\eta^{\beta\mu}{\partial^2 
F(t)\over \partial t^\lambda \partial t^\mu}\cr
A^{\alpha\beta} &=  \eta^{\alpha\lambda}\eta^{\beta\mu}A_{\lambda\mu}\cr}
$$
where the constant matrix $A_{\lambda\mu}$ was defined in (WDVV3). 

From (5.2) one obtains
$$
g^{\alpha\beta}(t) = t^1\, \eta^{\alpha\beta} + \tilde g^{\alpha\beta} 
\left( t^2, \dots, t^n\right)
$$
with
$$
 \tilde g^{\alpha\beta}
\left( t^2, \dots, t^n\right)= \sum_{\epsilon=2}^n E^\epsilon(t) 
c_\epsilon^{\alpha\beta}(t).
$$
So the bilinear form does not degenerate identically.
\smallskip
{\bf Definition 5.2.} The locus $\Sigma\subset M$ 
$$
\Sigma =\left\{ t \in M \, | \, \det \left( g^{\alpha\beta}(t)\right) 
=0\right\}
\eqno(5.3)
$$
is called {\it discriminant} of the Frobenius manifold $M$.
\medskip
{\bf Exercise 5.1.} Prove that the discriminant is specified by the
equation
$$
\det \U(t) =0
$$
where $\U(t)$ is the operator of multiplication by the Euler vector field.
\medskip
The inverse
$$
\left( g_{\alpha\beta}\right) = \left( g^{\alpha\beta}\right)^{-1}
\eqno(5.4)
$$
defines a metric on the open subset $M\setminus \Sigma$.
\smallskip
{\bf Lemma 5.1.} {\it 1). The Christoffel coefficients of the Levi-Civita 
connection for the metric (5.4) in the flat coordinates $t^\alpha$ are 
uniquely determined from the equation
$$
\Gamma_\gamma^{\alpha\beta}:= -g^{\alpha\epsilon} 
\Gamma_{\epsilon\gamma}^\beta = \left( {d+1\over 2} -q_\beta\right) 
c^{\alpha\beta}_\gamma.
\eqno(5.5)
$$

2). The metric (5.4) on $M\setminus \Sigma$ is flat.}

Proof can be found in [Du7].
\medskip
For brevity we will call the bilinear form $g^{\alpha\beta}(t)$ on $T^*_tM$
{\it contravariant metric} and the expressions 
$\Gamma_\gamma^{\alpha\beta}:= -g^{\alpha\epsilon}
\Gamma_{\epsilon\gamma}^\beta$ {\it contravariant Christoffel 
coefficients} of the Levi-Civita connection for the metric. 

We make a digression about linear pencils of contravariant metrics.

Let $\left( g_1^{ij}(x), {\Gamma_1}_k^{ij}(x)\right)$ and  $\left( 
g_2^{ij}(x), {\Gamma_2}_k^{ij}(x)\right)$ be two contravariant metrics 
invertible on an open subset of a manifold $M$ together with the 
correspondent  contravariant Christoffel coefficients.
\smallskip
{\bf Definition 5.3.} We say that the two contravariant metrics form a 
{\it linear quasihomogeneous pencil} of the charge $d$ if

1). For any $\lambda \in \C$ the metric
$$
g_1^{ij}(x) -\lambda\,g_2^{ij}(x)
$$
does not degenerate on an open subset in $M$.

2). The functions
$$
 {\Gamma_1}_k^{ij}(x)-\lambda\,  {\Gamma_2}_k^{ij}(x)
$$
are the contravariant Christoffel coefficients of the metric (5.4).

3). There exists a function $\varphi(x)$ on $M$ such that the vector fields
$$
E^i(x) := g_1^{ij}(x){\partial\varphi\over\partial x^j}, ~~e^i(x) := 
g_2^{ij}(x){\partial\varphi\over\partial x^j}
\eqno(5.6)
$$
satisfy the following properties
$$
\left[ e, \, E\right] =e
\eqno(5.7)
$$
$$
\L g_1^{ij}(x) =(d-1) \, g_1^{ij}(x), ~\L g_2^{ij}(x) = (d-2)\, 
g_2^{ij}(x), ~ {\cal L}_e g_1^{ij}(x) = g_2^{ij}(x), ~{\cal L}_e 
g_2^{ij}(x) =0.
\eqno(5.8)
$$
\medskip
{\bf Theorem 5.1.} {\it The intersection form of a Frobenius manifold 
together with the flat metric $<~,~>$ form a flat pencil of the charge $d$.}

Proof can be derived from Lemma 5.1 (see [Du7]). The function $\varphi(t)$ 
equals $\varphi = t_1 =\eta_{1\epsilon}t^\epsilon$.
\medskip
It can be shown [Du8] that, vice versa, a manifold with  a flat pencil of 
contravariant metrics satisfying certain assumptions about eigenvalues
of the linear operator $\nabla E$
carries a natural Frobenius structure such that, in 
the flat coordinates for $g_2^{ij}$, the metric $g_1^{ij}$ has the form 
(5.2) (cf. [Du7], [DZ1]).
\smallskip
{\bf Definition 5.4.} A function $x=x(t)$ is called {\it flat coordinate} 
of a metric if the differential $dx$ is covariantly constant w.r.t. the 
Levi-Civita connection for the metric.
\medskip
The flat coordinates of the intersection form on a Frobenius manifold are 
determined from the system of linear differential equations
$$
g^{\alpha\epsilon} \partial_\beta \xi_\epsilon + \sum_\epsilon \left(
{1\over 2} -\mu_\epsilon \right) c_\beta^{\alpha\epsilon}\xi_\epsilon =0
\eqno(5.9)
$$
where $\xi_\beta =\partial_\beta x$.
\smallskip
{\bf Definition 5.5.} The equations (5.9) are called {\it Gauss - Manin 
system of the Frobenius manifold}.
\medskip
{\bf Exercise 5.2.} Prove that the flat coordinates of the linear pencil
$g^{\alpha\beta}(t)-\lambda\,\eta^{\alpha\beta}$ have the form
$$
x(t^1-\lambda, t^2, \dots, t^n)
\eqno(5.10)
$$
where $x(t^1, t^2, \dots, t^n)$ are flat coordinates of the intersection 
form. Prove that the gradients  $\xi^\alpha =\eta^{\alpha\beta}\partial_\beta
x(t^1-\lambda, t^2, \dots, t^n)$ satisfy the system of equations
$$
\left(\U-\lambda\right) \partial_\beta \xi + C_\beta \left({1\over 2} 
+\mu\right) \xi =0
\eqno(5.11)
$$
$$
\left(\U-\lambda\right)\partial_\lambda \xi = \left({1\over 2}
+\mu\right) \xi.
\eqno(5.12)
$$

This is an extension of the Gauss - Manin system (5.9) onto $M\times 
\C_\lambda$. The second equation (5.12) has rational coefficients in 
$\lambda$. As above, compatibility of the full system will imply 
isomonodromicity of the Fuchsian system (5.12).
\medskip
{\bf Digression.} One can see by a straightforward computation that also 
the system
$$
\left(\U-\lambda\right)\partial_\beta \phi +C_\beta \mu\,\phi=0
\eqno(5.13a)
$$
$$
\left(\U-\lambda\right)\partial_\lambda \phi =\mu\,\phi
\eqno(5.13b)
$$
is compatible. We will use it to reduce WDVV for $n=3$, $d\neq 0$ to a 
particular case of Painlev\'e-VI equation (semisimplicity is assumed).
For $n=3$ the matrix $\mu$ degenerates
$$
\mu=\diag (\mu_1, \, 0, \, -\mu_1).
$$
So, the equations (5.13b) for the vector-function $\phi=(\phi_1, \phi_2, 
\phi_3)^T$ splits into a $2\times 2$ subsystem for $\chi =(\phi_1, \phi_3)^T$
and a quadrature for $\phi_2$
$$
{d\chi\over d\lambda} = -\mu_1 A(\lambda)\,\chi
\eqno(5.14)
$$
$$
A(\lambda) ={A_1\over \lambda-u_1} + {A_2\over \lambda-u_2} + {A_3\over 
\lambda-u_3}.
$$
Here $u_1$, $u_2$, $u_3$ are the eigenvalues of $U(t)$ (i.e., the 
canonical coordinates), the $2\times 2$-matrices have the form
$$
A(\lambda) =\mu_1 \left(\matrix{ v_1^1(\lambda; t) & -v_3^1(\lambda; t)\cr
v_1^3 (\lambda; t) & - v_3^3(\lambda; t) \cr}\right)
$$
where the matrix $\left( v^\alpha_\beta(\lambda; t)\right) 
:=\left(\U(t)-\lambda\right)^{-1}$,
$$
A_i =\left(\matrix{ \psi_{i1}\psi_{i3} & -\psi_{i3}^2\cr
\psi_{i1}^2 & - \psi_{i1}\psi_{i3}\cr}\right), ~~i=1, \, 2, \, 3.
\eqno(5.16)
$$
Clearly, the matrices satisfy the conditions
$$
\det A_i = \tr A_i =0, ~i=1, \, 2, \, 3
\eqno(5.17a)
$$
$$
A_1+A_2+A_3 =\left(\matrix{1 & 0\cr 0 & -1\cr}\right).
\eqno(5.17b)
$$
Following [JM] we introduce coordinates $p$, $q$, $k$ on the space of 
matrices $A_1$, $A_2$, $A_3$ satisfying (5.17). The coordinate $q$ is the 
root of the linear equation
$$
\left[ A(q)\right]_{12}=0;
\eqno(5.18a)
$$
the coordinate $p$ is the value
$$
p=\left[ A(q)\right]_{11}.
\eqno(5.18b)
$$
\def\half{{1\over 2}}
\def\deli{{\partial_i}}
Explicitly,
$$
\eqalign{q &= \left( g^{11} g^{22} - {g^{12}}^2\right)/g^{11}\cr
p &= \mu_1 {g^{11} g^{22}\over {g^{12}}^3 + g^{11} g^{12} g^{13} -
g^{11} g^{12} g^{22} - {g^{11}}^2 g^{23}}\cr}
\eqno(5.19)
$$
The entries of the matrices $A_i$ can be expressed via
the coordinates $p$, $q$, $k$
as follows
$$\eqalign{\psi_{i1}\psi_{i3} &=
-{q-u_i\over 2\mu_1^2 P'(u_i)} \left[ P(q) p^2 + 2 \mu_1 {P(q)\over q-u_i} p
+ \mu_1^2 (q + 2 u_i - \sum u_j)\right]
\cr
\psi_{i3}^2 &= -k {q-u_i\over P'(u_i)}
\cr
\psi_{i1}^2 &= -k^{-1} {q-u_i\over 4\mu_1^4 P'(u_i)}
\left[ P(q) p^2 + 2 \mu_1 {P(q)\over q-u_i} p
+ \mu_1^2 (q + 2 u_i - \sum u_j)\right]^2
\cr}
\eqno(5.20)
$$
where the polynomial $P(\lambda)$ has the form
$$
P(\lambda) := (\lambda-u_1)(\lambda-u_2)(\lambda-u_3).
\eqno(5.21)
$$
Compatibility of the system (5.13) implies
$$
\eqalign{\deli q &= {P(q)\over P'(u_i)}\left[ 2p+{1\over q-u_i}\right]\cr
\deli p &= -{P'(q) p^2 +\left(2q+u_i - \sum u_j\right) p + \mu_1 (1-\mu_1 )
\over P'(u_i)}\cr}
\eqno(5.22)
$$
and it gives a quadrature for the function $\log k$  
$$\deli \log k = (2\mu_1 -1) {q-u_i\over P'(u_i)}.
\eqno(5.23)
$$
Eliminating $p$ from the system we obtain a second order differential 
equation
for the function $q = q(u_1, u_2, u_3)$
$$
\deli^2 q = \half {P'(q)\over P(q)} \left(\deli q\right)^2
-\left[ \half {P''(u_i)\over P'(u_i)} +{1\over q-u_i}\right] \deli q
$$
$$+ \half {P(q)\over \left( P'(u_i)\right)^2}
\left[ (2\mu_1 -1)^2 +{P'(u_i)\over (q-u_i)^2}\right] , ~~ i=1, \, 2, \, 3.
$$
The system (5.22) is invariant w.r.t. transformations of the 
form
$$\eqalign{u_i &\mapsto a u_i + b\cr
q &\mapsto a q +
b.\cr}
$$
Introducing the invariant variables
$$\eqalign{x &={u_3-u_1\over u_2-u_1}\cr
y &= {q\over u_2 - u_1} - {u_1\over u_2-u_1}\cr}
$$
we obtain for the function $y = y(x)$ the following particular
Painlev\'e-VI equation
$$y'' = \half \left[ {1\over y} + {1\over y-1} + {1\over y-x}\right] (y')^2
- \left[ {1\over x}+{1\over x-1}+{1\over y-x}\right] y'
$$
$$+\half {y(y-1)(y-x)\over x^2(x-1)^2}\left[ (2\mu_1 -1)^2 +{x(x-1)\over 
(y-x)^2}
\right].
\eqno(PVI(\mu))
$$

Conversely, for a solution $y(x)$ of the equation $PVI(\mu)$ we 
construct
functions $q = q(u_1,u_2,u_3)$ and $p = p(u_1,u_2,u_3)$ putting
$$
\eqalign{q &= (u_2-u_1)y\left({u_3-u_1\over u_2-u_1}\right)+u_1\cr
p &= \half {P'(u_3)\over P(q)} y'\left( {u_3-u_1\over u_2-u_1}\right)
-\half {1\over q-u_3}.\cr}
$$
Then we compute the quadrature (5.23) determining the function $k$
(this provides us with one more arbitrary integration constant).
After this we are able to compute the matrix $\left(\psi_{i\alpha}
(u)\right)$  from the equations (5.20) and
$$
\left(\psi_{12}, \psi_{22}, \psi_{32}\right)
= \pm i\left( \psi_{21}\psi_{33} -\psi_{23} \psi_{31},
\psi_{13}\psi_{31}- \psi_{11}\psi_{33}, \psi_{11}\psi_{23}
-\psi_{13}\psi_{21}\right).
$$
The last step is in reconstructing the flat coordinates $t=t(u)$ and the 
tensor
$c_{\alpha\beta\gamma}$ using the formulae (3.21) and (3.17).
\medskip
{\bf Example 5.1.} Applying the above procedure to the three polynomial
solutions (1.22) - (1.24) of WDVV we obtain the following three algebraic
solutions
of $PVI(\mu)$ with $\mu = -1/4$, $-1/3$, $-2/5$ resp. [Du7,DM] represented
in
a parametric form
$$
\eqalign{
y & = {(s-1)^2(1+3s)(9s^2-5)^2\over(1+s)(25-207 s^2+1539 s^4+243 s^6)},
\cr
x &={(s-1)^3(1+3s)\over(s+1)^3(1-3s)}\cr}
\eqno(5.24)
$$ 
$$
\eqalign{
y &={(2-s)^2(1+s)\over(2+s)(5s^4-10s^2+9)}
\cr 
x &={(2-s)^2(1+s)\over(2+s)^2(1-s)}
\cr}
\eqno(5.25)
$$
$$
\eqalign{
y &={(s-1)^2(1+3s)^2(-1+4 s+s^2)(7-108 s^2+314 s^4-588 s^6+119 s^8)^2\over
 (1+s)^3(-1+3s) P(s^2)}
\cr
x &={{{{\left( -1 + s \right) }^5}\,{{\left( 1 + 3\,s \right) }^3}\,
     \left( -1 + 4\,s + {s^2} \right) }\over 
   {{{\left( 1 + s \right) }^5}\,{{\left( -1 + 3\,s \right) }^3}\,
     \left( -1 - 4\,s + {s^2} \right) }}
\cr}
\eqno(5.26a)$$
where
$$
\eqalign{
P(z)=&49-2133 z+34308 z^2-259044 z^3+16422878 z^4-7616646 z^5+
13758708 z^6\cr
&+5963724 z^7-719271 z^8+42483 z^9.\cr}
\eqno(5.26b)
$$

Some other particular solutions of Painlev\'e-VI in a relation with
Frobenius manifolds were constructed in [Se].
\bigskip
Let us return to the intersection form of a Frobenius manifold.
Due to Lemma 5.1 in a \nbh of a point $t_0 \in M\setminus \Sigma$ one can 
choose $n$ independent flat coordinates $x^1(t)$, \dots, $x^n(t)$ of the 
intersection form. In these coordinates the matrix
$$
g^{ab} =\left( dx^a, dx^b\right) = {\partial x^a\over \partial t^\alpha} 
{\partial x^b\over \partial t^\beta} g^{\alpha\beta}(t)
\eqno(5.27)
$$
becomes constant, and the Christoffel coefficients vanish. The flat 
coordinates are determined uniquely up to shifts and orthogonal 
transformations $\in O\left( n, g^{ab}\right)$.
\smallskip
{\bf Exercise 5.3.} Show that fo $d\neq 1$ the flat coordinates $x(t)$ can 
be chosen in such a way that
$$
\L x(t) = {1-d\over 2} x(t).
$$
\medskip
So, for $d\neq 1$, the flat coordinates of the intersection form satisfying
the quasihomogeneity condition of Exercise 5.3 are determined uniquely up
to a transformation from $ O\left( n, 
g^{ab}\right)$.

The solutions of the Gauss - Manin system can be continued analyticaly 
along any path in $M\setminus \Sigma$. We obtain a multivalued {\it period
map}
$$
t\mapsto \left( x^1(t), \dots, x^n(t)\right)
\eqno(5.28)
$$
defined on $M\setminus\Sigma$ (cf. the end of Lecture 2 above). The
multivaluedness of the period map is described by 
a representation of the fundamental group of the complement to the 
discriminant
$$
\pi_1\left( M\setminus \Sigma; t_0\right) \to  O\left( n, g^{ab}\right)
\eqno(5.29)
$$
(for $d=1$ instead of the orthogonal group we obtain a representation 
into the group of affine isometries of the metric $g^{ab}$).
\smallskip
{\bf Definition 5.6.} The image $W(M)$ of the representation (5.29)
is called {\it monodromy group} of the Frobenius manifold.
\medskip
Our main aim now is to compute 
the monodromy group of a semisimple Frobenius manifold in terms of
the Stokes matrix of the manifold.
 
In the semisimple case doing the gauge transform
$$
\phi =\Psi\,\xi
\eqno(5.30)
$$
we rewrite the extended Gauss - Manin system (5.11), (5.12) in the form
$$
\left(U-\lambda\right){d\phi\over d\lambda}=\left(\half +V\right) \, \phi
\eqno(5.31a)
$$
or, equivalently,
$$
{d\phi\over d\lambda}=\sum_{i=1}^n {B_i\over \lambda-u_i}\phi
\eqno(5.31b)
$$
with
$$
B_i = -E_i \left( \half +V\right)
\eqno(5.31c)
$$
where $E_i$ is the matrix unity (3.31), 
$$
\deli \phi = - {B_i\over \lambda-u_i}\phi
+V_i\phi, ~i=1, \dots, n.
\eqno(5.32)
$$
We obtain a Fuchsian system (5.31b) with the matrix residues $B_1$, \dots,
$B_n$
of a particular form (5.31c). Compatibility of (5.31a) with (5.31b)
provides 
isomonodromicity of the dependence of the coefficients of the system on 
the position of the poles $u_1$, \dots, $u_n$. We will now relate the
structure of the monodromy of the Fuchsian system (5.31b) to the Stokes
matrix
of the Frobenius manifold. 
\smallskip
{\bf Lemma 5.2.} {\it Let $\phi^{(1)}$, $\phi^{(2)}$ be two solutions of 
the system (5.31). Then the bilinear form
$$
\left( \phi^{(1)}, \phi^{(2)}\right) := { \phi^{(1)}}^T 
\left(U-\lambda\right)\, \phi^{(2)}
\eqno(5.33)
$$
does not depend neither on $\lambda$ nor on $u_1$, \dots, $u_n$.}

Proof can be obtained by straightforward differentiation.
\medskip
{\bf Remark 5.1.} We remember that the solutions $\phi =(\phi_1\dots, 
\phi_n)^T$ of the system (5.31) correspond to flat coordinates $x(t)$
of the intersection form
$$
\phi_i =\sum \psi_{i\alpha}\eta^{\alpha\beta} \partial_\beta 
x\left(t^1-\lambda, t^2, \dots, t^n\right).
\eqno(5.34)
$$
If $\phi^{(1)}$ corresponds to $x_1(t)$, $\phi^{(2)}$ to $x_2(t)$ then 
the bilinear form (5.33) equals
$$
\left( \phi^{(1)}, \phi^{(2)}\right) =(dx_1, dx_2) -\lambda <dx_1, dx_2>.
\eqno(5.35)
$$
\medskip
We will now construct, essentially following [BJL], a particular system of
solutions of the Fuchsian system (5.31b).
Let us choose an argument $\varphi$ in such a way that
$$
\arg (u_i-u_j) \neq {\pi\over 2} +\varphi \,({\rm mod}\, 2 \pi) ~{\rm for
~ 
any}~i\neq j.
\eqno(5.36)
$$
We make $n$ distinct parallel branchcuts $L_1$, \dots, $L_n$ on the 
complex $\lambda$-plane of the form
$$
L_j =\left\{ \lambda = u_j +i\rho e{-i\varphi}, ~~\rho\geq 0\right\} ,~ 
j=1, \dots, n.
\eqno(5.37)
$$
Each branchcut has positive and negative sides
$$
L_j^+ =\left\{ \lambda\,|\, \arg (u_j-\lambda) =-{\pi\over 2} -\varphi 
+0\right\} , ~~L_j^- =\left\{ \lambda\,|\, \arg (u_j-\lambda) =-{\pi\over 
2} -\varphi+2\pi
-0\right\} .
$$
On the complement
$$
\C\setminus \cup_j L_j
\eqno(5.38)
$$
the single-valued functions $\sqrt{u_1-\lambda}$, \dots, $\sqrt{u_n-\lambda}$
are well-defined. We specify them uniquely requiring that
$$
{\rm on}~ L_j^+ ~~\arg \sqrt{u_j-\lambda} = -{\pi\over 4} -{\phi\over 2} +0.
\eqno(5.39)
$$

Let us choose small loops $\gamma_1$, \dots, $\gamma_n$ going around the 
points $u_1$, \dots, $u_n$ in the counterclockwise direction. Let 
$R_1^*$, \dots, $R_n^*$ be the monodromy transformations in the space of 
solutions of (5.31b) corresponding to the loops  $\gamma_1$, \dots,
$\gamma_n$.
\smallskip
{\bf Lemma 5.3.} {\it 1). There exist unique solutions
$\phi^{(1)}(\lambda)$, 
\dots, $\phi^{(n)}(\lambda)$ of (5.31b) analytic in (5.38) such that
$$
R_j^* \phi^{(j)} =-\phi^{(j)}, ~~j=1, \dots, n
\eqno(5.40)
$$
$$
\phi_a^{(j)}(\lambda) ={\delta_a^j\over\sqrt{u_j-\lambda}} 
+O(\sqrt{u_j-\lambda}), ~~\lambda\to u_j.
\eqno(5.41)
$$

2). Introduce a symmetric matrix $G=\left( G_{ij}\right)$
$$
G_{ij} =\left( \phi^{(i)}, \phi^{(j)}\right).
\eqno(5.42)
$$
The monodromy transformations $R_j^*$ are the reflections
$$
R_j^* \phi^{(i)} = \phi^{(i)}- 2\, G_{ij}  \phi^{(j)}, ~~i,\, j=1, \dots, n
\eqno(5.43)
$$
in the hyperplanes orthogonal to $ \phi^{(j)}$ w.r.t. the bilinear form
(5.33).}

{\bf Proof.}  The matrix residue $B_j$ in (5.31b) has one eigenvalue
$-1/2$ and
$n-1$ eigenvalues $0$. 
So one can construct a fundamental system of solutions
$\phi^{(j)}(\lambda)$, $r_2(\lambda)$, \dots, $r_n(\lambda)$ such that
$R_j^* \phi^{(j)} = - \phi^{(j)}$, $R_j^*r_k = r_k$, $k=2, \dots, n$. That
means that 
the last $n-1$ solutions are analytic at $\lambda=u_j$. The solution
$\phi^{(j)}(\lambda)$ is determined uniquely up to a nonzero factor.
 From this it easily follows the first part of
Lemma.

To prove the second part let us represent $\phi^{(i)}(\lambda)$ as a 
linear combination of $\phi^{(j)}(\lambda)$ and of the solutions analytic 
at $\lambda = u_j$
$$
\phi^{(i)}(\lambda) =C_{ij} \phi^{(j)}(\lambda) + r_{ij}(\lambda).
$$
Here $C_{ij}$ is some constant, the solution $r_{ij}(\lambda)$ is 
analytic at $\lambda=u_j$. Computing the bilinear form (5.42) and using 
Lemma 5.2 we obtain
$$
G_{ij}=\lim_{\lambda\to u_j} \sum_{a=1}^n (u_a-\lambda) 
\phi_a^{(i)}(\lambda) \phi_a^{(j)}(\lambda) =C_{ij}.
$$
We obtain
$$
C_{ij} = G_{ij} ~{\rm for}~ i\neq j.
$$
Similar computation gives
$$
G_{ii}=1.
$$

We obtain a representation
$$
\phi^{(i)}(\lambda) =G_{ij} \phi^{(j)}(\lambda) + r_{ij}(\lambda).
\eqno(5.44)
$$
So
$$
\eqalign{
R_j^* \phi^{(i)} & = -G_{ij} \phi^{(j)}(\lambda) +r_{ij}(\lambda)\cr
 &= \phi^{(i)}(\lambda) -2 G_{ij} \phi^{(j)}(\lambda) \cr
 &= \phi^{(i)}(\lambda) -2 {\left( \phi^{(i)}, \phi^{(j)}\right)\over \left(
\phi^{(j)}, \phi^{(j)}\right)} 
\phi^{(j)}(\lambda). \cr}
$$
Lemma is proved.
\medskip
We will now establish, using the technique of [BJL],  a simple relation
between the matrix $G$ (5.42) for 
the system (5.31b) and the Stokes matrix of the operator (3.30).

Let us assume that the angle $\varphi$ is chosen in such a way that the 
order of the rays $L_1$, \dots, $L_n$ on the complex $\lambda$-plane
corresponds to the order of the complex numbers $u_1$, \dots, $u_n$ in 
the following sense: looking along the ray $L_j$ from the endpoint 
$\lambda=u_j$ we must see $L_{j-1}$ as the nearest ray on the left
and $L_{j+1}$ as the nearest one on the right, $2\leq j\leq n-1$.
\smallskip
{\bf Lemma 5.4.} {\it The oriented line $\ell = \ell_+ \cup (-\ell_-)$
$$
\ell_+ =\left\{ z\,|\, \arg z=\varphi\right\}
$$
is admissible for the operator (3.30). The correspondent Stokes matrix
$S$ is upper triangular. It satisfies the relation
$$
S+S^T = 2\, G
\eqno(5.45)
$$
where $G$ is the matrix (5.42).}

{\bf Proof.} Admissibility is obvious from (5.36). Let us construct the 
fundamental matrices $Y_\rr(z)$, $Y_\ll(z)$ having the needed asymptotic 
development (5.26) in the half-planes $\Pi_\rr$, $\Pi_\ll$. We will 
construct them taking an appropriate inverse Laplace transform of the 
solutions $\phi^{(j)}(\lambda)$ defined in Lemma 5.3. Put
$$
Y_{aj}(z) =-{\sqrt{z}\over 2 \sqrt{\pi}} 
\int_{C_j}\phi^{(j)}_a(\lambda)\,e^{\lambda\,z} d\lambda, ~~a,\,j =1, 
\dots, n.
\eqno(5.46)
$$
Here $C_j$ is an infinite contour coming from infinity along the positive 
side of the branchcut $L_j$, then encircling the point $\lambda = u_j$ 
and, after, returning to infinity along the negative side of the 
branchcut $L_j$. Since $\lambda=\infty$ is a regular singularity of the 
system (5.31b), the solutions $\phi^{(j)}(\lambda)$ grow at
$\lambda\to\infty$
not faster than a certain power of $|\lambda |$. We conclude that the 
integral converges absolutely for $z\in\Pi_\ll$. Using (5.31b) and
integrating
by parts we prove that the matrix $Y(z)=\left( Y_a^j(z)\right)$ satisfies 
the equation (3.30). To obtain the asymptotic development of this solution 
as $|z|\to\infty$ we can, due to Watson Lemma [WW], integrate the terms of 
the convergent expansion (5.41) of the solution $\phi^{(j)}(\lambda)$ near 
$\lambda=u_j$. Doing so we easily see that the solution $Y_\ll(z):=Y(z)$ 
has the needed asymptotic development 
$$
Y_{aj}(z) \sim \left( \delta_{aj} +O\left({1\over z}\right)\right) \, e^{z
u_j}
$$
as $|z|\to\infty$, $z\in\Pi_\ll$.

Let us now construct the fundamental matrix $Y_\rr(z)$. We are to choose 
the system of the opposite branchcuts
$$
L_j' =\left\{ \lambda=u_j -i\rho e^{-i\varphi}, ~\rho\geq 0\right\}, ~~
j=1, \dots, n
\eqno(5.47)
$$
to construct the correspondent solutions $\phi^{(j)'}(\lambda)$ and to define
$$
{Y_\rr}_{aj}(z) =
-{\sqrt{z}\over 2 \sqrt{\pi}}
\int_{C_j'}\phi^{(j)'}_a(\lambda)\,e^{\lambda\,z} d\lambda, ~~a,\,j =1,
\dots, n.
\eqno(5.48)
$$
Here the contour $C_j'$ goes around the branchcut $L_j'$. As above we prove
that the solution $Y_\rr(z) := \left({Y_\rr}_a^j(z)\right)$ of (3.30) has 
the needed asymptotic development in $\Pi_\rr$ as $|z|\to\infty$. It 
remains to establish a relation between the integrals (5.46) and (5.48).
To 
continue analytically $Y_\ll(z)$ through $\ell_+$ in the clockwise direction
into $\Pi_\rr$ we are to rotate the branchcuts $L_j$ in the counterclockwise
direction until they take the places of $L_j'$, $j=1, \dots, n$. For 
$j=1$ such a deformation does not meet obstructions. So
$$
\phi^{(1)'} =\phi^{(1)}.
$$
To deform $L_2$ to $L_2'$ we are to pass through the branchcut $L_1$. 
This is equivalent to the action of monodromy transformation $R_1^*$. So
$$
\phi^{(2)'}=R_1^* \phi^{(2)}.
$$
Continuing this process we obtain that
$$
\phi^{(k)'}=R_1^* R_2^* \dots R_{k-1}^* \phi^{(k)}, ~k=2, \dots, n.
$$
Using the computation in the proof of Coxeter identity (see [Bou])
we obtain
$$
\phi^{(k)} =2G_{k1}\phi^{(1)'}+2 G_{k2} \phi^{(2)'} +\dots
+ 2G_{k\, k-1} \phi^{(k-1)'} + \phi^{(k)'}, ~k=1, \dots, n.
$$
Lemma is proved.
\medskip
{\bf Corollary 5.1.} {\it If
$$
\det (S+S^T)\neq 0
\eqno(5.49)
$$
then the functions $\phi^{(1)}(\lambda)$, \dots, $\phi^{(n)}(\lambda)$ form
a basis of the space of solutions of (5.31b).}
\medskip
{\bf Exercise 5.4.} Prove that the Stokes matrix of quantum cohomology
of a manifold $X$ of an even complex dimension (assuming semisimplicity
of the quantum cohomology) satisfies the nondegeneracy condition (5.49).
\medskip
In the rest of this Lecture I will assume that $d\neq 1$ and that the 
Stokes matrix satisfies the nondegeneracy condition (5.49).

All the above constructions of the basis $\phi^{(1)}(\lambda)$, \dots, 
$\phi^{(n)}(\lambda)$ were done for a given fixed point $(u_1, \dots, u_n)$
of the Frobenius manifold. Since the solutions $\phi^{(j)}(\lambda)$
are determined uniquely, they become locally well-defined analytic functions
of $(u_1, \dots, u_n)$.
\smallskip
{\bf Lemma 5.5.} {\it They satisfy also the equations (5.32).}

Proof. Let us consider the vector-function
$$
\tilde\phi^{(j)} := \partial_i\phi^{(j)} +{B_i\over 
\lambda-u_i}\phi^{(j)}-V_i\phi^{(j)}
$$
for some $i$ between $1$ and $n$. Because of compatibility of (5.31) and
(5.32) 
the vector-function $\tilde\phi^{(j)}$ satisfies (5.31). It is easy to see 
that this solution is regular near the points $\lambda=u_1$, \dots, 
$\lambda=u_n$. Hence $\tilde\phi^{(j)}=0$. Lemma is proved.
\medskip
Let $\tilde M =\tilde{Fr}(e,\mu,R,S,C; u^0)$ be the Frobenius structure
on the universal covering of $\C^n\setminus\diag$ defined by the given 
monodromy data $(e, \mu, R, S, C)$ with $\mu_1\neq -1/2$, $\det 
(S+S^T)\neq 0$. The discriminant $\tilde \sigma$ of this Frobenius manifold
consists of the lifts of the coordinate hyperplanes $u_i=0$, $i=1, \dots, n$.

Let $\E$ be $n$-dimensional linear space equipped with a symmetric 
nondegenerate bilinear form $(~,~)$ on the dual space $\E^*$ having in 
some basis $e^1$, \dots, $e^n$ the Gram matrix
$$
(e^i, e^j) = \left( S+S^T\right)_{ij}.
\eqno(5.50)
$$
For any $1\leq i\leq n$ denote
$$
R_i :\E \to \E
\eqno(5.51)
$$
the transformation dual to the reflection $R^*_i: \E^*\to\E^*$ in the 
hyperplanes orthogonal to $e^i$:
$$
R_i^*(x) = x-(x,e^i) \,e^i.
\eqno(5.52)
$$
\smallskip
{\bf Theorem 5.2.} {\it The image of the monodromy representation
in the group $O\left(\E, (~,~)\right)$ of orthogonal transformations of 
the space $\E$
$$
\pi_1(\tilde M; u^0) \to O\left(\E, (~,~)\right)
\eqno(5.53)
$$
is the group generated by the reflections $R_1$, \dots, $R_n$.}

Proof. According to Lemma 5.3, locally we have a basis
$\phi^{(1)}(\lambda; 
u)$, \dots, $\phi^{(n)}(\lambda;u)$ of solutions of (5.31), (5.32). The
formula
$$
x_j (\lambda; u) ={2 \sqrt{2}\over 1-d} \sum_a (u_a-\lambda) 
\psi_{a1}(u)\phi_a^{(j)}(\lambda; u)
\eqno(5.54)
$$
gives flat coordinates of the linear pencil $(~, ~)-\lambda\, <~,~>$. Due 
to Lemma 5.4 we have
$$
(dx_i, dx_j) -\lambda\, <dx_i, dx_j> =\left(S+S^T\right)_{ij}.
\eqno(5.55)
$$
We obtain a locally well-defined isometry (the period map of the Frobenius
manifold, cf. (2.101))
$$
\tilde M\setminus \tilde \Sigma \to \E,
\eqno(5.56a)
$$
$$
u\mapsto \left( x_1(u), \dots, x_n(u)\right)
\eqno(5.56b)
$$
where
$$
x_j(u):= x_j(\lambda; u)|_{\lambda=0}.
$$
The monodromy around the branch $u_i=0$ of $\tilde\Sigma$ in the given 
chart of the universal covering of $\C^n\setminus\diag$ is equivalent to 
the monodromy of the vector-function
$$
\left( x_1(\lambda; u), \dots, x_n(\lambda; u)\right)
$$
corresponding to a small loop around $\lambda=u_i$ in the 
$\lambda$-plane. We obtain the transformation
$$
x_k(u) \mapsto x_k(u) -\left( S+S^T\right)_{ki} x_i(u), ~~k=1, \dots, n
\eqno(5.57)
$$
where we identify the coordinates in $\E$ with the dual basis in $\E^*$. 
That means that, locally, the monodromy group is generated by the 
reflections (5.57).

What happens with the analytic continuation into another chart of $\tilde 
M$? We arrive in another chart when some of the Stokes rays (4.34) passes 
through the admissible line $\ell$. Simultaneously, two of the branchcuts 
$L_1$, \dots, $L_n$ pass one through another one. The Stokes matrix 
changes according to the rule (4.68). It is sufficient to understand what 
happens with the flat coordinates $x_1(\lambda; u)$, \dots, $x_n(\lambda; u)$
with an elementary transformation (4.74) of the braid group.
\smallskip
{\bf Lemma 5.6.} {\it The elementary braid $\sigma_i$ permuting the points 
$\lambda=u_i$ and $\lambda=u_{i+1}$ in the complex $\lambda$-plane gives 
the following transformation of the solutions $\phi^{(1)}(\lambda)$, 
\dots, $\phi^{(n)}(\lambda)$
$$
\eqalign{
\sigma_i\left( \phi^{(k)}\right) &=\phi^{(k)}, ~~ k\neq i, \, i+1
\cr
\sigma_i\left(\phi^{(i)}\right) &= \phi^{(i+1)}
\cr
\sigma_i\left(\phi^{(i+1)}\right) &= R_{i+1}^* \phi^{(i)} = \phi^{(i)} - 
S_{i, i+1}\phi^{(i+1)}\cr}
\eqno(5.58)
$$
}

We leave the proof as an exercise to the reader.
\medskip
We obtain, after the analytic continuation along the braid $\sigma_i$, 
that the new monodromy transformations in $\E^*$ are reflections in the 
hyperplanes orthogonal to the vectors
$$
\left( e^1, \dots, e^{i-1}, e^{i+1}, R_{i+1}^*(e^i), e^{i+2}, \dots, 
e^n\right).
$$
But reflection w.r.t. the hyperplane orthogonal to $R_{i+1}^*(e^i)$ is 
equal to
$R_{i+1}^* R_i^* R_{i+1}^*$. This transformation belongs to the group 
generated by $R_1^*$, \dots, $R_n^*$. Theorem is proved.
\medskip
To complete our description of an arbitrary semisimple Frobenius manifold
in terms of an appropriate ``singularity theory'' we are to construct an
analogue of versal deformation. This can be done (at least, under the
nondegeneracy assumption (5.49)) in the following way ([Du7], Appendix I).
We
will construct
a family of functions $\lambda(p; u)$, $u=(u_1, \dots, u_n)$ of complex
variable $p$ defined in an open domain ${\cal D}$ of a Riemann surface
${\cal R}$
realized as branched covering over complex plane with finite number of
sheets.
The Riemann surface may depend on $u$. However, the projection of 
the domain ${\cal D}$ on the complex plane will be fixed. These functions
depend on complex pairwise distinct parameters $u_1$, \dots, $u_n$
belonging to a sufficiently small domain $\Omega \subset \C^n$. 

The first main property is that $\lambda(p; u)$ as function of $p\in {\cal
D}$ has critical values just $u_1$, \dots, $u_n$. The correspondent
critical points must not be degenerate. The second condition we require
from the function $\lambda (p; u)$ is that, for any two points
$p_i^{(1,2)}\in {\cal D}$ with the same critical value $u_i$, we must have
$$
\lambda''(p_i^{(1)}; u) = \lambda''(p_i^{(2)}; u).
$$
Here the prime denotes the $p$-derivative.
\smallskip
{\bf Definition 5.7.} The function $\lambda (p; u)$ on ${\cal D}\times
\Omega$ satisfying the above two properties is called {\it superpotential}
of some domain $M_\Omega$ in the Frobenius manifold $M$ if:

1). The canonical coordinates $(u_1, \dots, u_n)$ map $M_\Omega$ to
$\Omega\subset \C^n$.

2). For any critical points $p_1$, \dots, $p_n \in {\cal D}$
of $\lambda(p;u)$ with the critical values $u_1$, \dots, $u_n$ resp.
the following expressions for the flat metric $<~,~>$ on $T_tM$, the
intersection form $(~,~)$ (outside the discriminant $\Sigma$), and the
multiplication of tangent vectors hold true
$$
\left< \partial', \partial''\right>_t =-\sum_{i=1}^n \res_{p=p_i}
{\partial' \left( \lambda(p; u(t))dp\right) 
\partial'' \left( \lambda(p; u(t))dp\right)
\over
d\lambda(p; u(t))}
\eqno(5.59)
$$
$$
\left( \partial', \partial''\right)_t =-\sum_{i=1}^n \res_{p=p_i}
{\partial' \left( \log\lambda(p; u(t))dp\right)
\partial'' \left( \log\lambda(p; u(t))dp\right)
\over
d\log\lambda(p; u(t))}
\eqno(5.60)
$$
$$
\left< \partial'\cdot \partial'', \partial'''\right>_t =-\sum_{i=1}^n
\res_{p=p_i}
{\partial' \left( \lambda(p; u(t))dp\right)
\partial'' \left( \lambda(p; u(t))dp\right)
\partial''' \left( \lambda(p; u(t))dp\right)
\over
dp\, d\lambda(p; u(t))}
\eqno(5.61)
$$
In these formulae $\partial'$, $\partial''$, $\partial'''$ are any three
vector fields on $M$,
$$
d\lambda(p;u) := {\partial \lambda(p;u)\over \partial p} dp,
~~
d\log\lambda(p;u) := {\partial \log\lambda(p;u)\over \partial p} dp.
$$

3). For some 1-cycles $Z_1$, \dots, $Z_n$ in ${\cal D}$ the integrals
$$
\tilde t_j (u;z) ={1\over \sqrt{z}} \int_{Z_j} e^{z\, \lambda(p;u)} dp,
~~j=1, \dots, n
\eqno(5.62)
$$
converge and give a system of independent flat coordinates of the
connection $\tilde \nabla$. 
\medskip
{\bf Example 5.2.} For the polynomial Frobenius manifold corresponding to
$A_n$ singularity the versal deformation
$$
\lambda = p^{n+1} + a_n p^{n-1} +\dots, + a_1
$$
gives the needed superpotential. The variables $u_1$, \dots, $u_n$ are the
critical values of this function. Locally one can express the coefficients
of the polynomial as single-valued functions of $u_j$. For other simple
singularities the versal deformation is a family of polynomials of two
variables. However, one can reduce double integrals for the residues
(1.20), (1.21)
and for the oscilatory integrals (2.96) to one-dimensional residues and 
single integrals of the above
form. For the $D_n$ case this was done in [DVV]. (the superpotential
becomes
a rational function). For the case of $E_6$ singularity the superpotentil
is algebraic. It was found in [LW].
\medskip
{\bf Example 5.3.} For the case when the Riemann surfaces ${\cal R}$ can
be
compactified at infinity in such a way that there is exactly one
branch point on the Riemann surface over any of the critical values $u_j$
then the Frobenius manifold
can be identified with a Hurwitz space of branched coverings.
The Frobenius structure on the Hurwitz spaces was constructed in [Du1,
Du2]
(see also [Du7]).
In [Kr] the method of [Du1, Du2] was extended to produce also certain
algebro-geometric solutions satisfying WDVV1 and WDVV2 but not WDVV3.
\medskip
We will now construct a superpotential for semisimple Frobenius manifolds
satisfying nodegeneracy assumption (5.49). 

Let $u^0=(u_1^0, \dots, u_n^0)$, $u_i^0\neq u_j^0$ for $i\neq j$, be any
point of $M$ (written in the canonical coordinates). We choose the
branchcuts $L_1^0$, \dots, $L_n^0$ as above. For $M\ni u$ sufficiently
close to $u^0$ we will choose the branchcuts $L_1$, \dots, $L_n$
coinciding with $L_1^0$, \dots, $L_n^0$ outside some small neigborhoods
of the points $u_1^0$, \dots, $u_n^0$ resp. This allows us to construct
solutions $\phi^{(1)}(\lambda; u)$, \dots, $\phi^{(n)}(\lambda;u)$
as in Lemma 5.3. Denote
$$
G^{ij} =\left( \phi^{(i)}, \phi^{(j)}\right).
$$
Let $\left( G_{ij}\right)$ be the inverse matrix. Let us consider the
solution (cf. [BJL2])
$$
\phi (\lambda;u) =\sum_{i,\,j=1}^n G_{ij}\phi^{(j)}(\lambda;u).
\eqno(5.63)
$$
\smallskip
{\bf Lemma 5.7.} {\it The solution $\phi(\lambda; u) =\left(
\phi_a(\lambda;
u)\right)$ for $\lambda\to u_j$ has the behaviour
$$
\phi_a(\lambda; u) ={\delta_{aj}\over\sqrt{u_j-\lambda}}+O(1), ~a,\,j=1,
\dots, n.
\eqno(5.64)
$$
}

Proof follows from (5.44).
\medskip
Denote $p=p(\lambda; u)$ the correspondent flat coordinate of the
intersection form
$$
p(\lambda; u) ={\sqrt{2}\over 1-d} \sum_{a=1}^n (u_a-\lambda) \psi_{a1}(u)
\phi_a(\lambda; u)
\eqno(5.65)$$
i.e., 
$${\partial p(\lambda;u)\over \partial u_a} ={1\over \sqrt{2}}
\psi_{a1}(u) \phi_1(\lambda;u).
\eqno(5.66)
$$
It is analytic in the domain
$$
\C \setminus \cup_j L_j.
\eqno(5.67)
$$
For $\lambda\to u_j$ it behaves as follows
$$
p(\lambda; u) = p_j + \sqrt{2} \psi_{j1} \sqrt{u_j-\lambda} 
+O(u_j-\lambda)
\eqno(5.68)
$$
where $p_j = p(u_j; u)$. For $\lambda \to \infty$ the function $p(\lambda;
u)$ has a regular singularity. Hence $dp(\lambda; u)/d\lambda$ has at most
finite number of zeroes $r_1$, \dots, $r_N$ in (5.67). Without loss of
generality we may assume that all these zeroes are simple and that they do
not belong to the branchcuts $L_j$.

Let ${\cal D}_0$ be the image of the domain
$$
\lambda\in\C \setminus \cup_j L^0_j
$$
w.r.t. the map $p(\lambda; u^0)$. Particularly, the two sides of the
branchcut $L_j^0$  open to produce a smooth boundary curve of ${\cal D}_0$
passing through $p_j^0$. Denote
$$
\zeta_j^0 =p(r_j; u^0), ~j=1, \dots, N.
$$

Let us consider the inverse $\lambda=\lambda(p; u^0)$ to the function
$p(\lambda; u^0)$. It lives on a certain branched covering $\hat{\cal
D}_0$ of the domain ${\cal D}_0$ obtained by cutting ${\cal D}_0$ along
some paths going from $\zeta_1^0$, \dots, $\zeta_N^0$ to infinity and by
subsequent glueing of a finite number of copies of ${\cal D}_0$ with the
cuts. Near a point of the boundary of $\hat{\cal D}_0$ passing through
$p_j^0$ we have
$$
\lambda = u_j^0 -{1\over 2 \psi_{j1}^2 (u^0)} (p-p_j^0)^2 + O(p-p_j^0)^3.
$$
Thus we can analytically continue $\lambda(p; u^0)$ through the boundary
of ${\cal D}_0$ near a domain ${\cal D}$ containing $p_j^0$ as its
internal point. We can repeat this construction for all $u$ sufficiently
close to $u^0$ (actually, it is sufficient to require the points $u_j$ not
to intersect the branchcuts $L_i^0$). In this way we will produce a family
of Riemann surfaces with the branchpoints $\zeta_1$, \dots, $\zeta_N$. 
We may also assume that the image
of (5.67) w.r.t. the map $p(\lambda; u)$  for any $u$ close to $u^0$
belongs
to the projection of the domain ${\cal D}$ onto the complex $p$-plane.
This completes the construction of the family of functions $\lambda(p;
u)$.

The cycles $Z_j$ we need to compute the integrals (5.62) have the form
$$
Z_j = p(C_j; u)
\eqno(5.69)
$$
(more precisely, an arbitrary lift of this cycle on the Riemann surface)
where $C_j$ was defined in (5.46).
\smallskip
{\bf Theorem 5.3.} {\it The function $\lambda(p; u)$ is a superpotential
of the Frobenius manifold for a sufficiently small neighborhood of of the
point $u^0$.
}

Proof. By the construction the function $\lambda(p; u)$ has critical
values $u_1$, \dots, $u_n$. For any critical point $p_j\in {\cal D}$ on
the Riemann
surface with the critical value $u_j$ we obtained
$$
\lambda = u_j -{1\over 2 \psi_{j1}^2 (u)} (p-p_j)^2 + O(p-p_j)^3.
$$
So the second derivatives of $\lambda(p; u)$ do not depend on the choice
of the critical point.

Let us prove the formulae (5.59) - (5.61). We take $\partial'=\partial_a$,
$\partial''= \partial_b$, $\partial'''=\partial_c$ the vector fields along
the canonical coordinates. Then
$$
\partial_a \left( \lambda (p; u)dp \right) =
\left[ \delta_{aj} + O(p-p_j)\right] dp, ~p\to p_j,
$$
$$
d\lambda (p; u) = -\left[ {p-p_j\over \psi_{j1}^2} + O(p-p_j)^2\right] \,
dp, ~p\to p_j.
$$
So
$$
\res_{p=p_j} {\partial_a \left( \lambda(p; u) dp \right) \partial_b
\left( \lambda (p; u) dp\right)
\over
d\, \lambda (p; u)} =-\psi_{j1}^2 \delta_{aj}\delta_{bj}.
$$
Thus the formula (5.59) gives
$$
<\partial_a, \partial_b> =\delta_{ab} \, \psi_{a1}^2.
$$
This coincides with (3.15).

Similarly, (5.60) for $u_a\neq 0$ gives
$$
(\partial_a, \partial_b) =\delta_{ab}{\psi_{a1}^2 \over u_a}.
$$
This coincides with the definition of the intersection form written in the
canonical coordinates. Finally, the last formula (5.61) gives
$$
<\partial_a\cdot\partial_b, \partial_c>
=\delta_{ab}\delta_{ac}\psi_{a1}^2.
$$
This is equivalent to (3.15) together with the definition of the canonical
coordinates $\partial_a\cdot\partial_b =\delta_{ab}\partial_a$.

To prove (5.62) we use that the integrals
$$
\tilde t_j = -\sqrt{z} \int_{C_j} p(\lambda; u) e^{z\lambda}\,d\lambda,
~j=, \dots, n
\eqno(5.70)
$$
give flat coordinates of the deformed connection. Let us first prove
their independence. Indeed, the Jacobi matrix
$$
Y_{aj}(z; u) ={1\over \psi_{a1}} {\partial \tilde t_j\over \partial u_a}
= -\sqrt{z} \int_{C_j} \phi_a(\lambda; u) e^{z\lambda} d\lambda
$$
coincides with the fundamental matrix $Y^\ll(z; u)$ (up to a factor
$2\sqrt{\pi}$) due to Lemma 5.7 and Lemma 5.4. The final step of the
derivtion
is integration by parts and change of the integration variable $\lambda\to
p$:
$$
\eqalign{
\tilde t_j &= -{1\over \sqrt{z}} \int_{C_j}p(\lambda; u) \, de^{\lambda z}
\cr
&= {1\over \sqrt{z}} \int_{C_j}e^{\lambda z} {dp(\lambda; u)\over
d\lambda}\, d\lambda
\cr
&= {1\over\sqrt{z}}\int_{Z_j} e^{z\, \lambda(p; u)}dp.
\cr}
$$
Theorem is proved.

\medskip
{\bf Example 5.4.} Let the reflections $R_1^*$, \dots, $R_n^*$ generate a
finite group $W$ acting in the Euclidean space $\E$ of dimension $n$.
Recall [Bou]
that finite groups generated by reflections (5.43) are called {\it Coxeter
groups}. Let us assume the group $W$ to be an irreducible one. We will
construct, following [Du6], a Frobenius manifold $M_W$ with the
monodromy
group $W$.

The underlying  manifold of $M_W$ will be the orbit space
$$
M_W = \E/W.
\eqno(5.71)
$$
The coordinte ring of $M_W$ is, by definition, the ring of $W$-invariant
polynomials 
$$
\C [x_1, \dots, x_n]^W
$$ 
where $x_1$, \dots, $x_n$ are
Euclidean coordinates on $\E$. Due to Chevalley theorem [Bour] $M_W$ has a
natural structure of a graded affine algebraic variety:
$$
\C [x_1, \dots, x_n]^W \simeq \C [y^1, \dots, y^n]
\eqno(5.72)
$$
where
$$
y^i = y^i(x_1, \dots, x_n), ~~i=1, \dots, n
$$
are certain homogeneous $W$-invariant polynomials of the degrees
$$
d_i := \deg y^i(x) = m_i+1,
\eqno(5.73)
$$
$m_1$, \dots, $m_n$ are the exponents of the Coxeter group. The basic
invariant polynomials determine a coordinate system on the orbit space
$M_W$. The are determined uniquely up to invertible transformations of the
form
$$
y^i(x) \mapsto {y^i}'\left( y^1(x), \dots, y^n(x)\right), ~~i=1, \dots, n
$$
with quasihomogeneous polynomials ${y^i}'\left( y^1, \dots,
y^n\right)$ of the same degree $d_i$.

We will construct a polynomial Frobenius structure on $M_W$. That means
that the structure functions $c_{\alpha\beta}^\gamma$ will be elements
of the ring (5.72). Important ingredients of this contsruction will be
the Arnold's construction of convolution of invariants [Ar1, Gi1] and also
the flat coordinates on the orbit space $M_W$ discovered by K.Saito {\it 
et al.} in [Sai1, SYS].

Let $y^1(x)$ be the invariant polynomial of the maximal degree $h=d_1$.
The number $h$ is called {\it the Coxeter number} of the group $W$.
We define the unity vector field
$$
e:={\partial\over \partial y^1}
\eqno(5.74)
$$
and the Euler vector field
$$
E:= {1\over h} \sum _a x_a {\partial \over \partial x_a}.
\eqno(5.75)
$$
The unity vector field is defined up to a constant factor.

The construction of the metric $<~, ~>$ and of the multiplication law of
tangent vectors is more complicated. Let $(~,~)$ denote the $W$-invariant
Euclidean metric on the space $\E$. We will use the orthonormal
coordinates $x_1$, \dots, $x_n$ w.r.t. this metric.

Let us define
a bilinear symmetric form on
$T^*M_W$. In the coordinates $y^1$,\dots, $y^n$ it has the matrix
$$
(dy^i, dy^j) = \sum_{a=1}^n {\partial y^i\over \partial x_a} {\partial
y^j\over \partial x_a} = g^{ij}(y)
\eqno(5.76)
$$
for some polynomials $g^{ij}(y)$, $y=(y^1, \dots, y^n)$ (these exist
due to Chevalley theorem). The matrix $\left( g^{ij}(y)\right)$ is
invertible on $M_W\setminus \Sigma$ where the discriminant $\Sigma$
consists of all singular orbits.
\smallskip
{\bf Theorem 5.4}. {\it There exists a unique, up to an equivalence,
polynomial
Frobenius structure on the space of orbits of a finite Coxeter group with
the unity vector field (5.74), the Euler vector field (5.75), and the
intersection form (5.76).}

Sketch of the proof. We put
$$
<~,~> :={\cal L}_e (~,~)
\eqno(5.77)
$$
(cf. (5.8)). This gives {\it Saito metric} on the orbit space $M_W$.
According to [Sai1, SYS] this metric is flat, and there exists a
distinguished
system of basic homogeneous $W$-invariant polynomials $t^1(x)$, \dots,
$t^n(x)$ such that
$$
\eta^{\alpha\beta} :=<dt^\alpha, dt^\beta>
$$
is a constant nondegenerate matrix. Our main observation is that the
metrics $(~,~)$ and $<~,~>$ form a flat pencil (see above the definition).
This allows us to reconstruct the Frobenius structure inverting the
formula (5.2) (in the present case $A^{\alpha\beta}=0$ in (5.2)). Namely,
we define a $W$-invariant homogeneous polynomial $F$
of the degree $2h+2$ from the equations
$$
\eta^{\alpha\lambda}\eta^{\beta\mu}{\partial^2F\over \partial
t^\lambda\partial t^\mu} = {h \, (dt^\alpha, dt^\beta)\over \deg t^\alpha
+ \deg t^\beta -2}, ~~\alpha, \, \beta = 1, \dots, n
\eqno(5.78)
$$
(cf. (5.2)). Such a polynomial exists and it satisfies WDVV. Theorem is
proved.
\medskip
Observe that for the Frobenius structure on $M_W$ 
$$
d=1-{2\over h}, ~~q_\alpha = 1- {\deg t^\alpha \over h}, ~\alpha = 1,
\dots, n.
\eqno(5.79)
$$
\smallskip
{\bf Exercise 5.5.} Prove that the monodromy group of the Frobenius
manifold
$M_W$ is isomorphic to $W$. (Hint: prove that the flat coordinates
of the intersection form coincide with the Euclidean coordinates in the
space $W$.)
\medskip
{\bf Exercise 5.6.} Prove that the Frobenius manifolds $M_W$ satisfy the
semisimplicity condition.
\medskip
Particularly, for $n=2$ the polynomial Frobenius manifold corresponding to
the group $I_2(k)$ of symmetries of regular $k$-gon has the form
$$
F= {1\over 2} t_1^2 t_2 + t^{k+1}.
$$
For $n=3$ there are three irreducible finite Coxeter groups $W(A_3)$,
$W(B_3)$ and $W(H_3)$. They are the groups of symmetries of regular
tetrahedron, octahedron and icosahedron resp. The corresponding
polynomial Frobenius manifolds have the form (1.22), (1.23), and (1.24)
resp.
We give here aalso the list of all our polynomial Frobenius manifolds
in the dimension 4.
Group $W(A_4)$.
$$F = {1\over 2} t_1^2 \,t_4 + t_1\, t_2\, t_3 + {1\over 2} t_2^3 +
{1\over 3}
t_3^4 + 6 t_2 t_3^2 t_4 + 9 t_2^2 t_4^2 + 24 t_3^2 t_4^3 +{216\over 5}
t_4^6.$$
\medskip
Group $W(B_4)$.
$$F = {1\over 2} t_1^2 \,t_4 + t_1\, t_2\, t_3
+{{{ t_2}}^3} + {{{ t_2}\,{{{ t_3}}^3}}\over 3} +
  3\,{{{ t_2}}^2}\,{ t_3}\,{ t_4} +
  {{{{{ t_3}}^4}\,{ t_4}}\over 4}$$
$$ +  3\,{ t_2}\,{{{ t_3}}^2}\,{{{ t_4}}^2} +
  6\,{{{ t_2}}^2}\,{{{ t_4}}^3}
+ {{{ t_3}}^3}\,{{{ t_4}}^3} +
  {{18\,{{{ t_3}}^2}\,{{{ t_4}}^5}}\over 5} +
  {{18\,{{{ t_4}}^9}}\over 7}.$$
\medskip
Group $W(D_4)$.
$$F = {1\over 2} t_1^2 \,t_4 + t_1\, t_2\, t_3 +t_2^3 t_4+ t_3^3 t_4
+ 6 t_2 t_3 t_4^3 +{54\over 35} t_4^7.$$
\medskip
Group $W(F_4)$.
$$F = {1\over 2} t_1^2 \,t_4 + t_1\, t_2\, t_3 + { t_2^3\, t_4\over 18}
+ {3\, t_3^4\, t_4\over 4} + {t_2 \,t_3^2 \,t_4^3\over 2}$$
$$+ {t_2^2\, t_4^5\over 60} + {t_3^2\, t_4^7\over 28} +
{t_4^{13}\over 2^4 \cdot 3^2\cdot 11\cdot 13} .$$
\medskip
Group $W(H_4)$ 
$$F = {{{ t_1}\,{ t_2}\,{ t_3}}} +
  {{{{{ t_1}}^2}\,{ t_4}}\over 2} +
  {{{2\,{{ t_2}}^3}\,{ t_4}}\over {3}} +
  {{{{{ t_3}}^5}\,{ t_4}}\over {240}} +
  {{{ t_2}\,{{{ t_3}}^3}\,{{{ t_4}}^3}}\over {18}} +
  {{{{{ t_2}}^2}\,{ t_3}\,{{{ t_4}}^5}}\over {15}} +
  {{{{{ t_3}}^4}\,{{{ t_4}}^7}}\over {2^3\cdot 3^3\cdot 5}} +
  $$
$${{{ t_2}\,{{{ t_3}}^2}\,{{{ t_4}}^9}}\over {2\cdot 3^4\cdot 5}} +
  {{8\,{{{ t_2}}^2}\,{{{ t_4}}^{11}}}\over {3^4\cdot 5^2\cdot 11}} +
  {{{{{ t_3}}^3}\,{{{ t_4}}^{13}}}\over {2^2\cdot 3^6\cdot 5^2}} +
  {{2\,{{{ t_3}}^2}\,{{{ t_4}}^{19}}}\over {3^8\cdot 5^3\cdot 19}} +
  {{32\,{{{ t_4}}^{31}}}\over {3^{13}\cdot 5^6\cdot 29\cdot 31}}.$$

As it was shown in [Bl] these are all semisimple polynomial solutions
of WDVV for $n=4$ satisfying the conditions
$$
0<q_\alpha \leq d <1, ~\alpha = 2,\, 3,\, 4.
$$
\medskip
Other examples of polynomial solutions of WDVV associated with
finite Coxeter groups can be found in [Zu].
\smallskip
{\bf Remark 5.2}. There are certain inclusions between the polynomial
Frobenius
manifolds of the form
$$M_W := {\bf C}^n/W$$
(the orbit spaces)
for a finite Coxeter group $W$ acting in $n$-dimensional
Euclidean space. These inclusions correspond to the operation
of folding of Dynkin graphs [AGV]. As it is shown in [Ya], if the
Dynkin graph of a Coxeter group $W'$ is obtained by folding   
of the Dynkin graph of another Coxeter group $W$ then the corresponding
orbit space $M_{W'}$ is a (graded) linear subspace in $M_W$ w.r.t. the
Saito
linear structure. From our construction we immediately conclude
that the inclusion   
$$M_{W'}\subset M_W$$  
is also an embedding of Frobenius manifolds. We obtain the following
list of embeddings (they were obtained in [Zu] by
a straightforward computation)
$$
\eqalign{
M_{B_n} &\subset M_{A_{2n-1}}
\cr
M_{I_2(k)} &\subset M_{A_{k-1}}
\cr
M_{F_4} &\subset M_{E_6}
\cr
M_{H_3} &\subset M_{D_6}
\cr
M_{H_4} &\subset M_{E_8}.
\cr}
$$
(The group $W_{G_2}$ coincides with $W_{I_2(6)}$ and, therefore, $M_{G_2}
\subset M_{A_5}$.) The inclusions mean that, for example,
$$
F_{E_8} (t_1, 0, t_3, 0, 0, t_6, 0, t_8) =F_{H_4} (t_1, t_3, t_6, t_8).
$$
\medskip
{\bf Conjecture.} {\it Any irreducible semisimple polynomial Frobenius
manifold is equivalent to $M_W$ for some finite irreducible Coxeter
group $W$.}
\medskip
{\bf Remark 5.3.} According to our construction, the Frobenius structure
depends not only on the monodromy group $W$ but also on class of
equivalence of the ordered system of generating reflections $R_1^*$, 
\dots, $R_N^*$. The equivalence is established by simultaneous
conjugations of the generators by $(~,~)$-orthogonal transformations
and by the following action of the braid group
$$
\eqalign{
& \sigma_i(R_k^*) = R_k^*, ~k\neq i, \, i+1
\cr
& \sigma_i(R_i^*) = R_{i+1}^*
\cr
& \sigma_i (R_{i+1}^*) = R_{i+1}^* R_{i}^* R_{i+1}^*.
\cr}
\eqno(5.80)
$$
Here, as above,  $\sigma_i$ are the standard generators of the braid group
${\cal B}_n$. Any such class of equivalence is determined by the orbit of
the Stokes matrix $S=(S_{ij})$
$$
S_{ii}=1, ~S_{ij}=(e_i^*, e_j^*) ~{\rm for}~i<j
\eqno(5.81)
$$
w.r.t. the ${\cal B}_n$-action (5.80). Here $e_i^*$ is the basis of
normals
to the mirrors of the reflections normalized by the condition
$$
(e_i^*, e_i^*)=2
$$
for any $i$. For example, for an algebraic Frobenius manifold the orbit of
the given Stokes matrix $S$ must be finite. For the first nontrivial case
$n=3$ the classification of finite orbits of the action of ${\cal B}_3$
on the space of $3\times 3$ Stokes matrices satisfying the nondegeneracy
condition (5.49) was obtained in [DM]. Namely, there are only five finite
orbits. Three of them correspond to standard system of generating
reflections in the groups $W(A_3)$, $W(B_3)$, $W(H_3)$ of symmetries of
regular tetrahedron, octahedron and icosahedron respectively. Recall the
construction of a standard system of generating reflections in the group
of
symmetries of a regular polyhedron. Let $O$ be the center of the
polyhedron, $M$ the center of its face, $A$ a vertex of the face, $H$ the
center of an edge of the face having $A$ as an endpoint. Then the
reflections w.r.t. the planes $OMA$, $OMH$ and $OAH$ generate the group of
symmetries of the polyhedron [Cox]. Reordering of these generators give
the
same equivalence class.

One can repeat this construction with
cube
(the reciprocal of octahedron) or with dodecahedron (the reciprocal
of icosahedron) just to obtain the same system of generators in $W(B_3)$
and in $W(H_3)$ resp. Now we are able to describe the remaining two finite
orbits of the action of ${\cal B}_3$. The correspondent mirrors of the
reflections are obtained by applying the above construction to the great
icosahedron and great dodecahedron. The description of these regular
Kepler - Poinsot star-polyhedra one can find in the Coxeter book [Cox]. As
above, their reciprocals give the same equivalence class. All these
regular star-polyhedra have icosahedral symmetry. Thus, we obtain three
classes of triples of generating reflections in the group $W(H_3)$.

The above classification was applied in [DM] to the problem of
classification of algebraic solutions of $PVI(\mu)$. One can show that the
standard systems of generators in $W(A_3)$, $W(B_3)$, $W(H_3)$ correspond
to the polynomial solutions (1.22) - (1.24) of WDVV (thus, to the
algebraic solutions
(5.24) - (5.26) of $PVI(\mu)$). The last two finite orbits give algebraic
(non-polynomial) Frobenius manifolds.

The problem of classification for any $n$ of finite orbits of the action
of the braid group ${\cal B}_n$ in the space of $n\times n$ Stokes
matrices remains open. The solution of this problem could be useful to
prove the above Conjecture. For $n=4$ one can prove that all finite orbits
of irreducible Stokes matrices satisfying the nondegeneracy condition
(5.49)
correspond to a system of generating reflection of a finite Coxeter group
acting in ${\bf R}^4$. In the groups $W(A_4)$ and $W(B_4)$ there is only
one equivalence class of systems of generating reflections, namely, the
standard one. In the groups $W(D_4)$ and $W(F_4)$ there are two classes.
Finally, in the group $W(H_4)$ there are 10 classes of systems of
generating reflections. One of them corresponds to the standard system of
generators in the group of symmetries of the regular 600-cell, 6 others
to 4-dimensional regular star-polyhedra considered modulo reciprocity
(see the definitions in [Cox]) but the remaining 3 classes do not have
clear
geometrical meaning. Of the full list of 16 finite orbits we give here
Stokes matrices of representatives  in the finite ${\cal B}_4$-orbits
corresponding to only nonstandard systems of generators
$$
\eqalign{
D_4: ~S &=\left(\matrix{1 & -1 & 0 & 1 \cr
0 & 1 & -1 & -1 \cr
0 & 0 & 1 & 1\cr
0 & 0 & 0 & 1 \cr}
\right)
\cr
F_4: ~ S &= \left(\matrix{1 & -1 & 0 & \sqrt{2} \cr
0 & 1 & -\sqrt{2} & -\sqrt{2} \cr
0 & 0 & 1 & 1 \cr
0 & 0 & 0 & 1 \cr} \right)
\cr
H_4: ~S &= \left( \matrix{ 1 & -1 & 0 & {1+\sqrt{5}\over 2} \cr
0 & 1 & -1 & - {1+\sqrt{5}\over 2} \cr
0 & 0 & 1 &  {-1+\sqrt{5}\over 2} \cr
0 & 0 & 0 & 1 \cr} \right)
\cr
S &= \left( \matrix{ 1 & -1 & 0 &  {1+\sqrt{5}\over 2} \cr
0 & 1 & - {1+\sqrt{5}\over 2} & - {1+\sqrt{5}\over 2} \cr
0 & 0 & 1 & 1 \cr
0 & 0 & 0 & 1 \cr} \right)
\cr
S & =  \left( \matrix{ 1 & -1 & 0 &  {1-\sqrt{5}\over 2} \cr
0 & 1 & - {1-\sqrt{5}\over 2} & - {1-\sqrt{5}\over 2} \cr
0 & 0 & 1 & 1 \cr
0 & 0 & 0 & 1 \cr} \right)
\cr}
$$
\medskip
The construction of Theorem 5.4 was generalized in [DZ1] to produce
quasipolynomial solutions of WDVV, i.e., solutions with $d=1$ of the form
$$
F(t^1, \dots, t^n) = {\rm cubic}\, + f(t^2, \dots, t^{n-1}, \exp t^n)
$$
with a polynomial $f$. The monodromy group of these Frobenius manifolds
are certain extensions of affine Weyl groups. Particularly, the monodromy
group of quantum cohomology of ${\bf CP}^1$ is given by this construction
(see [Du7]).
\medskip
{\bf Example 5.5.} Let us compute the monodromy group of the quantum
cohomology of ${\bf CP}^2$. The Stokes matrix $S$ (4.97) of this Frobenius
manifold
satisfies the nondegeneracy condition (5.49). The basic reflections in the
monodromy group in the basis of the flat coordinates
$(x,y,z)$ corresponding to
the basis $\Phi^\ll=\left(\Phi_1^\ll, \Phi_2^\ll, \Phi_3^\ll\right)$  of
the
solutions (4.96) have the matrices
$$
R_1^* =\left(\matrix{-1 & -3 & 3 \cr 0 & 1 & 0 \cr 0 & 0 & 1 \cr}\right),
~R_2^* = \left(\matrix{ 1 & 0 & 0 \cr -3 & -1 & 3 \cr 0 & 0 & 1
\cr}\right),
~R_3^* =\left(\matrix{1 & 0 & 0 \cr 0 & 1 & 0 \cr 3 & 3 & -1 \cr}\right).
$$
The full monodromy group of the Frobenius manifold we obtain adding the
monodromy transformation corresponding to the only nontrivial loop
$$
t_2\mapsto t_2+2 \pi i.
$$
This corresponds to the rotation
$$
z\mapsto z\, e^{2 \pi i\over 3}.
$$
Using the explicit formulae (4.96) and the identity (4.92) we immediately
obtain
the needed transformation
$$
\Phi^\ll \left( z\, e^{2 \pi i\over 3} \right) =\Phi^\ll (z)\, T
$$
or, equivalently,
$$
(x,y,z)\mapsto (x,y,z)\,T
\eqno(5.82)
$$
with
$$
T=\left(\matrix{ 0 & -1 & 0 \cr 0 & 0 & 1\cr -1 & -3 & 3 \cr}\right).
\eqno(5.83)
$$
For the matrix
$$
T_0 := T\, R_1^* =\left(\matrix{0 & -1 & 0 \cr 0 & 0 & 1 \cr 1 & 0 &
0\cr}\right)
\eqno(5.84)
$$
we have the identities
$$
T_0^3 = -1, ~~R_2^* = T_0^{-1} R_1^* T_0, ~~ R_3^* = T_0^{-1} R_2^* T_0.
\eqno(5.85)
$$
So, we introduce the new system of generators $A$, $B$, $C$ in the full
monodromy group putting
$$
A=R_1^*, ~~B=T_0^4 = -T_0, ~~ C=T_0^3 = -1.
\eqno(5.86)
$$
All the transformations of the group preserve the integer lattice in ${\bf
R}^3$. They also preserve the indefinite quadratic form with the Gram
matrix $S+S^T$
$$
q(x,y,z) = 2\left( x^2 + y^2 + z^2 + 3 x y - 3 x z - 3 y z\right).
\eqno(5.87)
$$
The group acts discretely on the complexification of the cone
$q(x,y,z)>0$.

Introducing the coordinates $r$, $\tau$, $\bar\tau$
$$
\eqalign{
x &= {i r\over 2} {2 \tau \bar\tau -3 ( \tau +\bar\tau) + 2\over  \tau -
\bar\tau}
\cr
y & = {i r \over 2} {2  \tau \bar\tau +  \tau +\bar\tau -2 \over  \tau
-\bar\tau}
\cr
z &= {i r\over 2} {2  \tau \bar\tau - ( \tau +\bar\tau) -2\over  \tau
-\bar\tau}
\cr}
\eqno(5.88)
$$
we obtain the action of the generating transformations 
$$
\eqalign{
A: &\left( r\mapsto r, ~ \tau \mapsto -{1\over \tau}, ~ \bar\tau \mapsto
-{1\over \bar\tau}\right)
\cr
B: & \left( r\mapsto r, ~ \tau \mapsto {1\over 1-\tau}, ~ \bar\tau \mapsto
{1\over 1-\bar\tau}\right)
\cr
C: & \left( r \mapsto -r, ~ \tau \mapsto \tau, ~ \bar\tau \mapsto \bar\tau
\right).
\cr}
\eqno(5.89)
$$
We proved
\smallskip
{\bf Theorem 5.5.} {\it The monodromy group of quantum cohomology of
${\bf CP}^2$
is isomorphic to $PSL_2({\bf Z})\times \{ \pm \}$.}
\medskip
It would be interesting to develop an appropriate theory of invariants
for this action of the modular group. This could help to obtain analytic
formulae for quantum cohomology of ${\bf CP}^2$.

\vfill\eject
\centerline{\bf References}
\medskip
\item{[Ama]} Amari S.-I., Differential-Geometric Methods in Statistics,
Springer Lecture Notes in Statistics {\bf 28} 1985.
\medskip
\item{[Ar1]} Arnol'd V.I., Wave front evolution and equivariant Morse lemma,
{\sl Comm. Pure Appl. Math.} {\bf 29} (1976) 557 - 582.
\medskip
\item{[Ar2]} Arnol'd V.I., Singularities of Caustics and Wave Fronts,
Kluwer Acad. Publ., Dordrecht - Boston - London, 1990.
\medskip
\item{[AGV]} Arnol'd V.I., Gusein-Zade S.M., and Varchenko A.N.,
Singularities
of Differentiable Maps, volumes I, II, Birkh\"auser, Boston-Basel-Berlin,
1988.
\medskip
\item{[At]} Atiyah M.F., Topological quantum field theories, {\sl
Publ. Math.
I.H.E.S.} {\bf 68} (1988) 175.
\medskip
\item{[BJL1]}
Balser W., Jurkat W.B., and Lutz D.A., Birkhoff invariants
and Stokes multipliers for meromorphic linear differential equations,
{\sl J. Math. Anal. Appl.} {\bf 71} (1979), 48-94.
\medskip
\item{[BJL2]} Balser W., Jurkat W.B., and Lutz D.A., On the reduction of
connection
problems for differential equations with an irregular singular point
to ones with only regular singularities, {\sl SIAM J. Math. Anal.} {\bf 12}
(1981) 691 - 721.
\medskip
\item{[BS]} Barannikov S., Kontsevich M., Frobenius manifolds and 
formality of Lie algebras of polynomial vector fields, alg-geom/9710032.
\medskip
\item{[Beh]} Behrend K., Gromov - Witten invariants in algebraic geometry,
{\sl Inv. Math.} {\bf 127} (1997) 601 - 627.
\medskip
\item{[Beau]} Beauville A., Quantum cohomology of complete intersections, 
alg-geom/9501008. 
\medskip
\item{[Bi]} Birman J., Braids, Links, and Mapping Class Groups,
Princeton Univ. Press, Princeton, 1974.
\medskip
\item{[Bl]} Blanco Cede\~{n}o A., On polynomial solutions of equations of 
associativity, Preprint ICTP IC/97/152.
\medskip
\item{[BV]} Blok B. and Varchenko A., Topological conformal 
field theories
and the flat coordinates, {\sl Int. J. Mod. Phys.} {\bf A7} (1992) 1467.
\medskip
\item{[Bou]} Bourbaki N., Groupes et Alg\`ebres de Lie, Chapitres 4, 5 et 6,
Masson, Paris-New York-Barcelone-Milan-Mexico-Rio de Janeiro, 1981.
\medskip
\item{[CL]} Coddington, E. A., Levinson, N., Theory of ordinary
differential equations, New York, McGraw-Hill 1955.
\medskip
\item{[COGP]} Candelas P., de la Ossa X.C., Green P.S., and Parkes L.,
A pair of Calabi - Yau manifolds as an exactly soluble superconformal   
theory, {\sl Nucl. Phys.} {\bf B359} (1991), 21-74.
\medskip
\item{[CV1]} Cecotti S. and Vafa C., Topological-antitopological fusion,
{\sl Nucl. Phys.} {\bf B367} (1991) 359-461.
\medskip
\item{[CV2]} Cecotti S. and  Vafa C., On classification of $N=2$ 
supersymmetric
theories, {\sl Comm. Math. Phys.} {\bf 158} (1993), 569-644.
\medskip
\item{[Cox]} Coxeter, H.S.M., Regular polytopes, New York, Macmillan
1963.
\medskip
\item{[DI]} Di Francesco P., and Itzykson C., Quantum intersection rings,
Preprint SPhT 94/111, September 1994.
\medskip
\item{[Dij1]} Dijkgraaf R., Intersection theory, integrable 
hierarchies and
topological field theory, Preprint IASSNS-HEP-91/91, December 1991.
\medskip
\item{[Dij2]} Dijkgraaf R., Notes on topological string theory and 2D
quantum gravity,
Preprint PUPT-1217, IASSNS-HEP-90/80, November 1990.
\medskip
\item{[DVV]} Dijkgraaf R., E.Verlinde, and H.Verlinde, Topological
strings in $d<1$, {\sl Nucl. Phys.}
{\bf B 352} (1991) 59.
\medskip
\item{[DW]} Dijkgraaf  R., and Witten E., Mean field theory, topological
field theory, and multimatrix models,
{\sl Nucl. Phys.} {\bf B 342} 
(1990) 486-522.
\medskip
\item{[Du1]} Dubrovin B., Differential geometry of moduli spaces and its
application to soliton equations and to topological field theory,
Preprint No.117,  Scuola Normale Superiore, Pisa (1991).
\medskip
\item{[Du2]} Dubrovin B., Hamiltonian formalism of Whitham-type hierarchies
and topological Landau - Ginsburg models, {\sl Comm. Math. Phys.}
{\bf 145} (1992) 195 - 207.
\medskip
\item{[Du3]} Dubrovin B., Integrable systems in topological field theory,
{\sl Nucl. Phys.} {\bf B 379} (1992) 627 - 689.
\medskip
\item{[Du4]} Dubrovin B., Geometry and integrability of
topological-antitopological fusion, {\sl Comm. Math. Phys.}{\bf 152}
(1993), 539-564.
\medskip
\item{[Du5]} Dubrovin B., Integrable systems and classification of 
2-dimensional topological field theories,
In \lq\lq Integrable Systems", Proceedings
of Luminy 1991 conference dedicated to the memory of J.-L. Verdier.
Eds. O.Babelon, O.Cartier, Y.Kosmann-Schwarbach, Birkh\"auser, 1993.
\medskip
\item{[Du6]} Dubrovin B., Differential geometry of the space of orbits
of a Coxeter group, Preprint SISSA-29/93/FM (February 1993).
\medskip
\item{[Du7]} Dubrovin B., Geometry of 2D topological field theories,
In: ``Integrable Systems and Quantum Groups'', Eds. M.Francaviglia, S.Greco,
Springer Lecture Notes in Math. {\bf 1620} (1996) 120 - 348.
\medskip
\item{[Du8]} Dubrovin B., Flat pencils of metrics and Frobenius manifolds,
math.DG/9803106, to appear in Proceedings of 1997
Taniguchi Symposium ``Integrable Systems and Algebraic Geometry''.
\medskip
\item{[DM]} Dubrovin B., Mazzocco M., Monodromy of certain Painlev\'e-VI
transcendents and reflection groups, Preprint SISSA 149/97/FM.
\medskip
\item{[DZ1]} Dubrovin B., Zhang Y., Extended affine Weyl groups and 
Frobenius manifolds, {\sl Compositio Math.} {\bf 111} (1998) 167-219.
\medskip
\item{[DZ2]} Dubrovin B., Zhang Y., Bihamiltonian hierarchies in 2D 
topological field theory at one-loop approximation, Preprint 
SISSA 152/97/FM, hep-th/9712232, to appear in {\sl Comm. Math. Phys.}.
\medskip
\item{[DM]} Duval A., Mitschi C., Matrices de Stokes et groupe de Galois
des equations hypergeometriques confluentes generalisees, {\sl Pacific J.
Math.} {\bf 138} (1989) 25 - 56.
\medskip
\item{[EKYY]} Eguchi T., Kanno H., Yamada Y., and 
Yang S.-K., Topological
strings, flat coordinates and gravitational descendants, {\sl Phys. 
Lett.}{\bf
B305} (1993), 235-241;
\medskip
\item{[EYY1]} Eguchi T., Yamada Y., and Yang S.-K., Topological field 
theories 
and the period integrals, {\sl Mod. Phys. Lett.} {\bf A8} (1993), 1627-1638.
\medskip
\item{[EYY2]} Eguchi T., Yamada Y., and Yang S.-K., On the genus expansion
in the topological string theory, Preprint UTHEP-275 (May 1994).
\medskip
\item{[Ga]} Gantmakher F.R., The theory of matrices, New York, Chelsea
Pub. Co., 1960.
\medskip
\item{[Gi1]} Givental A.B., Convolution of invariants of groups generated
by reflections, and connections with simple singularities of functions,
{\sl Funct. Anal.} {\bf 14}  (1980) 81 - 89.
\medskip
\item{[Gi2]} Givental A.B., Equivariant Gromov - Witten invariants, 
alg-geom/9603021.
\medskip
\item{[Gi3]} Givental A.B., Stationary phase integrals, quantum Toda 
lattice, flag manifolds, and the mirror conjecture, alg-geom/9612001.
\medskip
\item{[Gi4]} Givental A.B., A mirror theorem for toric complete 
intersections,
alg-geom/9701016.
\medskip
\item{[Gi5]} Givental A.B., Elliptic Gromov - Witten invariants
and the generalized mirror conjecture,
math.AG/9803053.
\medskip
\item{[Gr]} Gromov M., Pseudo-holomorphic curves in symplectic
manifolds, {\sl Invent. Math.} {\bf 82} (1985), 307.
\medskip
\item{[Hi1]} Hitchin N., 
Poncelet polygons and the Painlev\'e equations, In: ``Geometry and
Analysis", Tata Inst. of Fundamental Rsearch, Bombay, 1995, pp. 151-185.
\medskip
\item{[Hi2]} Hitchin N., Frobenius Manifolds (Notes by David Calderbank),
Preprint, 1996.
\medskip
\item{[Ho]} Hori K., Constraints for topological strings in $D\geq 1$,
{\sl Nucl. Phys.} {\bf B 439} (1995) 395 - 420.
\medskip
\item{[IN]} Its A.R.,  and Novokshenov V.Yu., {\it The Isomonodromic 
Deformation
Method in the Theory of Painlev\'e Equations}, {\sl Lecture Notes in
Mathematics}
1191, Springer-Verlag, Berlin 1986.
\medskip
\item{[JM]} Jimbo M., and Miwa T., Monodromy preserving deformations
of linear ordinary differential equations with rational coefficients. II.
{\sl Physica} {\bf 2D} (1981) 407 - 448.
\medskip
\item{[Ka]} Kaufmann R., The geometry of moduli spaces of pointed curves,
the tensor product in the theory of Frobenius manifolds and the explicit
K\"unneth formula in quantum cohomology, Ph.D. thesis, MPI, Bonn, 1997.
\medskip
\item{[Ko1]} Kontsevich M., Intersection theory on the moduli space of
curves, {\sl Funct. Anal.} {\bf 25} (1991) 50.
\medskip
\item{[Ko2]} Kontsevich M., Intersection theory on the moduli space of curves
and the matrix Airy function,
{\sl Comm. Math. Phys.} {\bf 147} (1992) 1-23.
\medskip
\item{[Ko3]} Kontsevich M.,
Enumeration of rational curves via torus action, {\sl Comm. Math. Phys.}
{\bf 164} (1994) 525 - 562.
\medskip
\item{[KM1]} Kontsevich M., Manin Yu.I., Gromov - Witten classes, quantum
cohomology and enumerative geometry, {\sl Comm. Math. Phys.} {\bf 164}
(1994) 525 - 562.
\medskip
\item{[KM2]} Kontsevich M., Manin Yu.I., Quantum cohomology of a product
(with Appendix by R.Kaufmann), {\sl Inv. Math.} {\bf 124} (1996) 313 -
339. 
\medskip
\item{[Kr]} Krichever I.M.,
 The $\tau$-function of the universal Whitham hierarchy, matrix models
and topological field theories, {\sl Comm. Pure Appl. Math.} {\bf 47}
(1994) 437 - 475.
\medskip
\item{[LW]} Lerche W., Warner N., Exceptional SW Geometry from ALE 
Fibrations,  Preprint hep-th/9608183.
\medskip
\item{[LLY]} Lian B.H., Liu K., Yau S.-T., Mirror principle I, 
hep-th/9712011.
\medskip
\item{[Loo]} Looijenga E., A period mapping for certain semi-universal
deformations, {\sl Compositio Math.} {\bf 30} (1975) 299 - 316.
\medskip
\item{[Los]} Losev A., ``Hodge strings'' and elements of K.Saito's theory 
of the primitive forms, hep-th/9801179.
\medskip
\item{[Lu]} Luke Y. L., Mathematical functions and their approximations, 
New York Academic Press 1975.
\medskip
\item{[Ma]} Malgrange B., \'Equations Diff\'erentielles \`a
Coefficients Polynomiaux, Birkh\"auser, 1991.
\medskip
\item{[MS].} McDuff D., Salamon D., J-holomorphic curves and quantum
cohomology,
Providence, RI., American Mathematical Society, 1994.
\medskip
\item{[Man1]} Manin Yu.I., Frobenius manifolds, quantum cohomology, and
moduli spaces, Preprint MPI 96-113.
\medskip
\item{[Man2]} Manin Yu.I., Sixth Painlev\'e equation, universal elliptic
curve, and mirror of $P^2$, alg-geom/9605010.
\medskip
\item{[Man3]} Manin Yu.I., Three constructions of Frobenius manifolds:
a comparative study, Pre\-print math.AG/9801006.
\medskip
\item{[MM]} Manin Yu.I., Merkulov S.A., Semisimple Frobenius (super)manifolds
and quantum cohomology of $P^r$, alg-geom/9702014.
\medskip 
\item{[Mi]} Miwa T., 
Painlev\'e property of monodromy presereving
equations and the analyticity of $\tau$-functions, {\sl Publ. RIMS}
{\bf 17} (1981), 703-721.
\medskip
\item{[Pi]} Piunikhin S., Quantum and Floer cohomology have the same ring 
structure,
Preprint MIT (March 1994).
\medskip
\item{[RT]} Ruan Y., Tian G., A mathematical theory of quantum cohomology,
{\sl Math. Res. Lett.} {\bf 1} (1994), 269-278.
\medskip
\item{[Sab]} Sabbah C., Frobenius manifolds: isomonodromic deformations
and infinitesimal period mappings, Preprint, 1996.
\medskip
\item{[Sad]} Sadov V., On equivalence of Floer's and quantum 
cohomology,     
Preprint HUTP-93/A027.
\medskip
\item{[Sai1]} Saito K., On a linear structure of a quotient variety 
by a finite
reflection group, Preprint RIMS-288 (1979).
\medskip
\item{[Sai2]} Saito K., Period mapping associated to a primitive form,
{\sl Publ. RIMS} {\bf 19} (1983) 1231 - 1264.
\medskip
\item{[SYS]} Saito K., Yano T., and Sekeguchi J., On a certain generator
system  
of the ring of invariants of a finite reflection group,
{\sl Comm. in Algebra} {\bf 8(4)} (1980) 373 - 408.
\medskip
\item{[Sat]} Satake I., Flat structure of the simple elliptic
singularity of type $\tilde E_6$ and Jacobi form, hep-th/9307009.
\medskip
\item{[Se]} Segert J., Frobenius manifolds from Yang - Mills instantons,
dg-ga/9710031.
\medskip
\item{[Ta]} Takahashi Atsushi, Primitive forms, topological LG models
coupled to gravity and mirror symmetry, math.AG/9802059.
\medskip
\item{[TX]} Tian Gang, Xu Geng, On the semisimplicity of the quantum 
cohomology algebra of complete intersections, alg-geom/9611035.
\medskip
\item{[Wa]} Wasow W., Asymptotic expansions for ordinary differential
equations, Wiley, New York, 1965.
\medskip
\item{[WW]} Whittaker E.T., and Watson G.N., A Course of Modern Analysis,
N.Y. AMS Press, 1979.
\medskip
\item{[Wi1]} Witten E., On the structure of the topological phase of
two-dimensional
gravity,
{\sl Nucl. Phys.} {\bf B 340} (1990) 281-332.
\medskip
\item{[Wi2]} Witten E., Two-dimensional gravity and intersection theory on
moduli
space,  
{\sl Surv. Diff. Geom.} {\bf 1}  (1991) 243-210.
\medskip
\item{[Wi3]} Witten E., Lectures on mirror symmetry, In: [Yau].
\medskip
\item{[Yan]} Yano T., Free deformation for isolated singularity,
{\sl Sci. Rep.   
Saitama Univ.} {\bf A 9} (1980) 61 - 70.
\item{[Yau]} Yau S.-T., ed., Essays on Mirror Manifolds, International 
Press Co., Hong Kong, 1992.
\medskip
\item{[Zu1]} Zuber J.-B., On Dubrovin's
topological field theories, 
{\sl Mod. Phys. Lett.} {\bf A 9} (1994) 749 - 760.
\medskip
\item{[Zu2]} Zuber J.-B., Graphs and reflection groups, {\sl Comm.
Math. Phys.} {\bf 179} (1996) 265 - 294.
\medskip
\item{[Zu3]} Zuber J.-B., Generalized Dynkin diagrams and root systems and
their folding, hep-th/9707046.
\vfill\eject\end